\documentclass[10pt,a4paper]{article}

\usepackage{amsmath,amssymb,amsthm,bbm,graphicx,xcolor,xargs,sectsty}
\usepackage[square,compress,comma, numbers,sort]{natbib}
\usepackage{amsfonts,mathtools}
 \usepackage[colorlinks=true, linkcolor=blue,citecolor=green]{hyperref}

\usepackage[shortlabels]{enumitem}
\usepackage[title,titletoc,header]{appendix}
\usepackage{mathrsfs}   

\def\proj{ \mathfrak{p}}
\def\projt{ \mathfrak{p}_t}

\newcommand\spaceD{\mathrm{D}}

\makeatletter
\def\@fnsymbol#1{\ensuremath{\ifcase#1\or *\or ** \or \ddagger\or
		\mathsection\or \mathparagraph\or \|\or **\or \dagger\dagger
		\or \ddagger\ddagger \else\@ctrerr\fi}}
\makeatother

\usepackage{cleveref}

\DeclarePairedDelimiter{\parens}()
\def\Pspace{ \parens{\Omega,\mathscr{F},\mathbb{P}}}
\def\Dspace{ \parens{\spaceD,\mathscr{D}}}
\def\DTspace{ \parens{\spaceD,d_\spaceD}}

\def\Dtau{d_\spaceD}

\def\LK{\lambda_{\TO}} 
\def\borel#1{\mathscr{B}(#1)}
\def\borelD{\borel{\spaceD}}
\def\borelR{\mathscr{B}(\R)}

\def\AA{\mathscr{D}}

\def\borelV{\borel{\spaceD}}
\def\dV{d_\spaceD}

\newcommand\VV{\spaceD}
\def\AAV{\AA}
\def\TTTK{K(\TO)}

\crefname{hypothesis}{Assumption}{Assumptions}
\crefname{equation}{}{}
\crefname{enumi}{}{}

\def\ve{\varepsilon}

\numberwithin{equation}{section}
\subsubsectionfont{\itshape}

\allowdisplaybreaks[4]

\newcommand{\QED}{\hfill $\Box$}

\newcommand{\abs}[1]{\left\lvert #1 \right\rvert}

\DeclarePairedDelimiterXPP\pk[1]{\mathbb{P}}\{ \}{}{ #1}
\DeclarePairedDelimiterXPP\E[1]{\mathbb{E}}\{ \}{}{	#1}

\newcommand{\nelem}[1]{{Lemma \ref{#1}}}
\newcommand{\neprop}[1]{{Proposition \ref{#1}}}
\newcommand{\netheo}[1]{{Theorem \ref{#1}}}

\def\IF{\infty}

\newcommand{\BQN}{\begin{eqnarray}}
	\newcommand{\EQN}{\end{eqnarray}}
\newcommand{\BQNY}{\begin{eqnarray*}}
	\newcommand{\EQNY}{\end{eqnarray*}}
\def\ldot{, \ldots,}

\newcommand{\limit}[1]{\lim_{#1 \to   \infty}}
\def\akn#1{\begin{align} #1 \end{align} }

\def\bkny#1{ \begin{eqnarray*} #1 \end{eqnarray*}}
\def\bkn#1{ \begin{eqnarray} #1 \end{eqnarray}}
\def\bqn#1{ \begin{eqnarray} #1 \end{eqnarray}}
\def\bqny#1{ \begin{eqnarray*} #1 \end{eqnarray*}}

\newcommand{\kb}[1]{\boldsymbol{#1}}
\newcommand{\vk}[1]{\kb{#1}}

\newcommand{\COM}[1]{}

\definecolor{c20}{rgb}{0.,0.7,0.}
\definecolor{c30}{rgb}{0.,0.,1.}
\definecolor{c40}{rgb}{1,0.1,0.7}
\definecolor{c50}{rgb}{1,0,0}
\definecolor{c60}{rgb}{1,0.9,0.1}
\definecolor{c70}{rgb}{0.50,1.00,0.00}

\def\mb#1{{#1}}

\def\tE#1{{\textcolor{c30}{#1}}}
\def\eEE#1{{\textcolor{c50}{#1}}}
\def\vEE#1{{\textcolor{c20}{#1}}}

\def\eEE#1{{#1}}
\def\vEE#1{{#1}}
\def\cEE#1{{#1}}
\def\cEV#1{{\textcolor{c50}{#1}}}
\def\cEV#1{{#1}}
\def\cEG#1{{\textcolor{c50}{#1}}}
\def\cEG#1{{#1}}

\def\tE#1{#1}

\newcommand{\zhora}[2]{{#1}\marginpar{\footnotesize #2}}

\newcommand{\Nset}{\mathbb{N}}
\newcommand{\Rset}{\mathbb{R}} 
\newcommand{\Zset}{\mathbb{Z}}

\def\dualpha{\,\alpha r^{-\alpha-1} d r\,}

\newcommand{\pr}{\mathbb{P}}

\newcommand{\convprob}{\stackrel{P}{\longrightarrow}} 
\newcommand{\convvague}{\stackrel{v,\mcb}{\longrightarrow}}
\newcommand{\convvagueO}{\stackrel{v, \mcb_0}{\longrightarrow}}
\newcommand{\convweak}{\stackrel{w}{\Longrightarrow}}

\newcommandx{\Norm}[1]{\left|#1\right|}

\newcommandx{\norm}[2][1=]{|\!|#2|\!|_{#1}}
\newcommandx{\normE}[1]{|\!|#1|\!|}

\newcommandx{\lpnorm}[3][1=,3=]{\left\|#2\right\|_{#1}^{#3}}
\newcommandx{\supnorm}[3][1=,3=]{\left\|#2\right\|_{#1}^{#3}}

\newcommand\exc{\mathcal{E}^q}

\newcommand\nualpha{v_\alpha}
\newcommand\tailmeasure{\boldsymbol{\nu}}

\newcommand\ind[1]{\mathbbm{1}{\left(#1\right)}}

\def\boldsymbol#1{#1}
\def\bszero{\boldsymbol{0}}

\newcommand{\bsy}{f}

\newcommand{\bsX}{\boldsymbol{X}}
\newcommand{\bsY}{\boldsymbol{Y}}
\newcommand{\bsZ}{\boldsymbol{Z}}

\def\bTheta{\boldsymbol\Theta}

\def\bszero{\boldsymbol{0}}

\def\bY{\boldsymbol{Y}}
\def\bX{\boldsymbol{X}}

\def\mcb{\mathcal{B}}

\def\cadlag{c\`adl\`ag}
\def\Cadlag{C\`adl\`ag}

\def\eg{e.g.,}
\def\ie{i.e.,}

\newtheorem{theorem}{Theorem}[section]
\newtheorem{lemma}[theorem]{Lemma}
\newtheorem{corollary}[theorem]{Corollary}
\newtheorem{proposition}[theorem]{Proposition}

\newtheorem{definition}[theorem]{Definition}

\newtheorem{remark}[theorem]{Remark}
\newtheorem{example}[theorem]{Example}
\newtheorem{Condition}[theorem]{Condition}

\newcommand{\prooftheo}[1]{ \textsc{Proof of Theorem} \ref{#1} }
\newcommand{\proofprop}[1]{\textsc{Proof of Proposition} \ref{#1}}
\newcommand{\prooflem}[1]{\textsc{Proof of Lemma} \ref{#1}}
\newcommand{\proofkorr}[1]{\textsc{Proof of Corollary} \ref{#1}}
\def\tPr{ \bTheta^{[h]} }

\newcommand{\BEL}{\begin{lemma}}
	\newcommand{\EEL}{\end{lemma}}

\newcommand{\BRM}{\begin{remark}}
	\newcommand{\ERM}{\end{remark}}

\newcommand{\BD}{\begin{definition}}
	\newcommand{\ED}{\end{definition}}

\def\RL{\Rset^l}

\def\spaceD{\DDD}
\def\DDD{\mathrm{D}}

\def\TTT{\mathrm{T}}
\def\TO{\mathcal{Q}}
\def\TT{\TTT}
\def\R{\Rset}

\def\SSY{ \mathcal{S}^q}
\def\intTO{\int_{\TO}}

\def\RVD{ \mathscr{R}_\alpha(\mcb_0, \tailmeasure )}

\def\RP{\Rset_{>}}
\def\RPP{\Rset_{\ge}}

\def\TO{\mathcal{Q}}

\def\mM{\mathcal{M}^+}
\def\mS{\mathcal{M}_\alpha}
\def\mSD{\mathcal{M}_\alpha(\AA)}

\def\AAK{\AA_{\TO}}

\def\zD{0_\spaceD}
\def\zD{0}

\def\epsx{  (\ve z)^{-1}}
\def\epsilon{\ve}
\def\RA{\mathcal{R}_\alpha}
\def\RNA{\cEG{\mathcal{R}_{1/\alpha}} }
\def\RAA{\mathscr{R}_\alpha}
\def\RVBA{ \mathcal{R}(g, \mcb, \tailmeasure)}
\def\RVB{ \RA(g, \mcb, \tailmeasure)}
\def\RVC{ \widetilde{\RA}(a_n, \mcb, \tailmeasure)}
\def\MBB{\mM(\mcb)}

\def\MBV{\mM(\mcb)}

\def\DOT{ \!\cdot \!}
\def\DOT{\bullet}
\def\DOT{\cdotp}
\def\DOT{\boldsymbol{\cdot}}
\makeatletter
\newcommand*\bigcdot{\mathpalette\bigcdot@{.55}}
\newcommand*\bigcdot@[2]{\mathbin{\vcenter{\hbox{\scalebox{#2}{$\m@th#1\bullet$}}}}}
\makeatother
\def\DOT{\! \bigcdot \!}

\def\famHMAP{\mathcal{H}(\spaceD)}
\begin{document}

\title{Tail Measures and Regular  Variation}
\vspace{ -2 cm}

\author{Martin Bladt\footnote{University of Lausanne, Switzerland }, Enkelejd Hashorva$^*$ 
	and Georgiy Shevchenko\footnote{Taras Shevchenko National University of Kyiv,  Ukraine }}

\COM{
\address{Enkelejd Hashorva, Department of Actuarial Science, 
	University of Lausanne,\\
	UNIL-Dorigny, 1015 Lausanne, Switzerland
}

\author{Enkelejd  Hashorva}
\address{Enkelejd Hashorva, Department of Actuarial Science, 
	University of Lausanne,\\
	UNIL-Dorigny, 1015 Lausanne, Switzerland
}

\author{Georgiy Shevchenko}
\address{Georgiy Shevchenko,  Department of Probability Theory, Statistics and Actuarial Mathematics, Taras Shevchenko National University of Kyiv, Volodymyrska 64, Kiev 01601, Ukraine \\
	Tel/fax.: +38044-259-03-92
}	

}
\maketitle 
\begin{abstract}  
A general framework for the study of regular variation (RV) is that of Polish star-shaped metric spaces, while recent developments in \cite{kulik:soulier:2020} have discussed RV with  respect to a  properly localised boundedness $\mcb$. Along the lines of the latter approach, we discuss the RV of Borel measures and random processes on a general Polish metric spaces $(\spaceD, d_\spaceD)$. Tail measures introduced in \cite{owada:samorodnitsky:2012} appear naturally as limiting measures of regularly varying time series. 
	We define tail measures on the measurable space $(\spaceD, \AA)$ indexed by $\famHMAP$, a countable family of 1-homogeneous coordinate maps,  and show some tractable instances for the investigation of RV when $\mcb$ is determined by $\famHMAP$. This allows us to study the regular variation of  \cadlag\ processes on $D(\R^l, \R^d)$ retrieving in particular results obtained in \cite{PH2020} for RV of stationary  \cadlag\ processes on the real line removing $l=1$ therein. \cEG{Further, we discuss potential applications and open questions}. 
\end{abstract}

\cadlag\ processes;  max-stable processes; tail processes; spectral tail processes;  weak convergence; 


\tableofcontents

\section{Introduction}
Let ${\bsX}(t), t\in \TTT$ be an $\R^d$-valued random process with $\TTT$ a  non-emtpy set (here $d,l,k$ are  reserved for positive integers). For given  $t_1 \ldot t_k \in \TTT$ and $A \in \borel{\R^{dk}} $   it is of interest for many investigations  to determine  the asymptotic behaviour as $n\to \IF$ of $$p_{t_1 \ldot t_k}(a_n\DOT A)=\pk* { \left({\bsX(t_1)}{},\dots,{\bsX(t_k)}{}\right) \in  a_n \DOT A }$$ for some  positive scaling constants $a_n,n\ge 1$.  
 Such an  investigation is reasonable if $a_n \DOT A, n\ge 1$ are Borel absorbing events, i.e., the outer multiplication $\DOT	$ satisfies $a_n \DOT A \in \borel{\R^{dk}}$ and  $\limit{n}p_{t_1 \ldot t_k}(a_n \DOT  A)=0$.  Throughout this paper $\borel{\VV}$ stands for the Borel $\sigma$-field  on the topological space  $\VV$.\\
Considering for simplicity the canonical scaling, i.e., $ c \DOT A:= \{ c \DOT  a, a\in A\}$ for all $c\in (0,\IF)$, where $\DOT$ is the usual product on $\Rset$, it  is   natural to require that the  $A$ is  separated from the origin   (denoted   by 0) of $\R^{dk}$, i.e.,  $A$  is included in the complement of
a neighbourhood of $\bszero$ in the {usual} topology. 
\cEE{For such $A$,  the  rate of convergence to 0 of}
 $p_{t_1\ldot t_k}(a_n \DOT A)$  is the main topic in the theory of RV of random vectors.   Indeed the  RV of functions,  random processes and Borel measures is  important in various research fields  and is not confined to probabilistic applications, see e.g., \cite{MR4029007} and the references therein. \\ 
 The problem at hand can be regarded as a   scaling approximation discussed  for instance in \cite{MR3331244} in terms of Kendall's theorem and is investigated in  the framework  of   RV  of measures, in finite or infinite dimensional spaces, see e.g., \cite{MR1925445,MR2271424,hult:lindskog:2005,hult:lindskog:2006,RDTM:2008,MR3561100,SegersEx17,kulik:soulier:2020,MR4036028,PH2020,Planic}.\\
  As, for instance, in  \cite{Davis,BojanS,owada:samorodnitsky:2012,kulik:soulier:2020}, we say that $\bsX$ is finite dimensional regularly varying if there exist positive $a_n$'s  such that for all 
$t_1 \ldot  t_k\in\TTT, k\ge 1,$ there exists a non-null  measure $\nu_{t_1,\dots,t_k}$ on
$\borel{\Rset^{dk}}$ satisfying 
\bkny{  \lim_{n\to\infty}
	n\pk* { \left( {\bsX(t_1)}{},\dots, {\bsX(t_k)}{}\right) \in a_n \DOT A }
	= \nu_{t_1,\dots,t_k} (A) < \infty
}
for all $\nu_{t_1,\dots,t_k}$-continuity $A\in \borel{\R^{dk}}$ separated from $\bszero$. 
The measure $\nu_{t_1,\dots,t_k}$ is called the exponent measure of $(X(t_1),\dots,X(t_k))$. If the outer multiplication $\DOT$ is the usual product, it is well-known that the exponent measures are  $-\alpha$-homogeneous, \ie\  there exists $\alpha>0$ (not depending on $t_i$'s) such that 
\bkn{\label{exaAH} 
	\nu_{t_1,\dots,t_k}(z \DOT A)=z^{-\alpha} \nu_{t_1,\dots,t_k}(A), \quad \forall z\in (0,\IF), \forall ( t_1,\dots,t_k )\in \TTT^k, \forall k\ge 1.
} 
 RV of Borel measures on some Polish metric space $\DTspace$ is investigated in \cite{hult:lindskog:2005,hult:lindskog:2006,MR3271332,SegersEx17,kulik:soulier:2020}. 
  The recent manuscripts \cite{kulik:soulier:2020,PH2020} treat  RV  of measures and processes in terms of a given \eEE{properly localised boundedness $\mcb$} on $\spaceD$ following the ideas in \cite{MR3979308}. In \cite{kulik:soulier:2020}  several weak  conditions are formulated with respect to the scaling and the topology of  $\spaceD$, see  \cite{kulik:soulier:2020}[Appendix B: (M1)-(M3), (B1-B3)].  
We highlight next some key developments and findings: 
\begin{enumerate} [{\bf F}1]
	\item \label{it:1} All investigations  in the literature, e.g., 	\cite{MR1825160,hult:lindskog:2005,hult:lindskog,RDTM:2008,MR3238572} consider RV of random processes with  compact parameter space   $\TTT\subset \R^l, l\in \Nset$. \cEV{Moreover, RV of Borel measures on star-shaped Polish metric spaces are considered. Surprisingly, the non-compact case $\TTT=\R^l$, which is of great interest for the investigation of time series,  has been  investigated only recently in \cite{PH2020} for stationary \vEE{stochastically continuous} \cadlag\ random processes when $l=1$;}
	\item \label{it:2} \cEG{The recent manuscripts \cite{kulik:soulier:2020, PH2020} develop the theory of RV with respect to a properly localised boundedness $\mcb$. This new approach  has several advantages including the unification of RV and hidden RV; }	
	\item \label{it:3}     RV  of \cEG{stationary} time series can be characterised by  tail and spectral tail processes, see 
	\cite{BojanS, Hrovje,WS,basrak:planinic:2019, kulik:soulier:2020}.
 See	also \cite{SegersEx17, klem,MR4153591} for non-stationary time series where also local tail processes play a crucial role for the characterisation of RV;
	\item \label{it:4}  Characterisation \cEG{of  RV of stationary time series}  in terms of tail measures is first investigated in  \cite{owada:samorodnitsky:2012} and further discussed in  \cite{Hrovje,klem,PH2020};
\item \label{it:5} \vEE{There are different definitions of RV useful in various applications, which in view of  \cite{hult:lindskog:2006}[Thm 3.1] are equivalent for star-shaped Polish metric spaces;}	
\end{enumerate}
\Cref{it:1}--\Cref{it:5}  \cEG{and recent applications developed in \cite{PH2020}} \mb{motivate the following two {topics}, which constitute the backbone of the present contribution}:
\begin{enumerate}[{\bf T}1]
	\item \label{them:2}  RV  of processes \cEG{(not necessarily stationary)} with non-compact parameter space $\TTT$, or in general  RV  of Borel measures  \mb{in non-star-shaped Polish metric spaces \eEE{with respect to some properly localised boundedness $\mcb$}};
	\item \label{them:1} Basic properties of   tail measures in general measure spaces  and their relationship with RV;
	\item \label{them:3} \cEG{Relation between RV and local tail processes};
	\item \label{them:4} \cEG{Potential applications of RV to \cadlag\ processes (random fields) with non-compact $\TTT$};
	\item \label{them:5} \cEG{Discussion on possible different definitions of RV relevant for applications}.
	\end{enumerate} 
RV  of  \vEE{stochastically continuous} stationary \cadlag\ processes defined on the real line is  recently investigated in \cite{PH2020}.
For the case of locally compact $\TTT=\R^l$  
the  corresponding functional metric spaces (we denote them by $(\spaceD, \Dtau)$ below) are not radially monotone (or star-shaped, see \cite{SegersEx17}), which is the case when $\TTT$ is compact. 
Specifically, for a hypercube $\TTT\subset\R^l, l\ge 1$ and $D(\TTT, \R^d)$ the space of generalized \cadlag\ functions $\TTT \mapsto \R^d$ (see e.g., \cite{MR0402847,JJ} for definitions), a metric $\Dtau$  can be chosen so that $D(\TTT, \R^d)$ is Polish and
$$\Dtau(c  \DOT f, \zD)= c \Dtau(f,\zD), \ \ \forall c>0,\forall f\in D(\TTT, \R^d),$$
where  $\zD$ is the zero function. Consequently, $\Dtau(cf, \zD)$ is strictly monotone for all $c>0$ and fixed \zhora{$f\neq 0$}{}; this is referred to as the radial monotonicity property and has been a key assumption in the treatment of RV of measures in Polish metric spaces, e.g., \cite{hult:lindskog:2005,Segers}.\\
  When  $\TTT=\R^l$, in view of \netheo{sam},\Cref{tev:6} in Appendix, radial monotonicity does not hold. That property  is crucial for the proof of \cite{hult:lindskog:2006}[Thm 3.1]. Therefore when dealing with  
  $D(\R^l,\R^d)$ the equivalence of different definitions of RV of Borel measures  does not follow from the aforementioned theorem, but can be nonetheless confirmed as shown in \nelem{chans}.
  
Following \cite{kulik:soulier:2020}, where 
a boundedness along with the chosen group action plays a crucial role, we discuss first  RV of Borel measures on general Polish metric spaces.
From \cite{owada:samorodnitsky:2012,Hrovje,klem,PH2020} it is known that for particular Polish spaces the limit measure in the definition of RV is a tail measure,  \eEE{which is essentially characterised by the following   properties:  
\begin{enumerate}[{P}1)]
	\item \label{P1} $-\alpha$-homogeneity as described by \eqref{exaAH};
	\item \label{P2}  countable indexing by 1-homogeneous maps.
\end{enumerate}
}
 In abstract settings, \Cref{P1}   is introduced   under the assumption that $(\spaceD, \DOT \ , \AA, \RP, \cdot)$ is a measurable cone with 
 $\AA$ being a $\sigma$-field  on $\spaceD$,   \ie\ the outer multiplication (we prefer here the formulation as a pairing)  $(z,f) \mapsto z \DOT     f \in \spaceD$ for $z\in \RP$ and $f\in \spaceD$   is a group action from the multiplicative group $(\RP, \DOT \,)$ on $\spaceD$ and  is jointly measurable.\\
 {Hereafter},    $Z$ is a $\spaceD$-valued   random element defined on a complete probability space $\Pspace$. The cone measurability assumption and the Fubini-Tonelli theorem yield that
 \bqn{ \tailmeasure_Z (A) = \E*{ \int_0^\IF \ind{   z \DOT \vk Z \in  A}  \alpha z^{- \alpha- 1} dz  }, 
 	\quad A \in \AA
 	\label{repNU}
 }
is a  non-negative  measure  on $\AA$ for all $\alpha>0$. It follows   that $\tailmeasure=\tailmeasure_Z$ satisfies
\begin{enumerate}[$M0$)]
	\item \label{L:P1}  	$\tailmeasure( t\DOT A) = t^{-\alpha} \tailmeasure(A), \ \ \forall t \in (0,\IF), \forall  A\in \AA$
\end{enumerate}
and therefore  $\tailmeasure$ is called $-\alpha$-homogeneous, with $\alpha$ referred to as its  index.\\ 
\eEE{If the space $\spaceD$ is a countable product of measurable spaces, then \Cref
{P2} can be introduced with respect to a given  $1$-homogeneous positive definite (point-separating) measurable map as in \cite{klem}. \vEE{The essential feature of such a map is that it defines a  cone on $\spaceD$ and moreover its support can be countably generated by Borel sets of $\R$.} In this paper we do not restrict to such measurable spaces and therefore the countable indexing is introduced below in \Cref{L:P3} 
with respect to a countable family $\famHMAP$ of  $1$-homogeneous measurable  maps, see Definition \ref{defEL}.}\\  
A crucial  consequence of both \Cref{P1}-\Cref{P2} is that the  introduced tail measures $\tailmeasure$ are $\sigma$-finite. Moreover, the family $\famHMAP$
allows us to introduce the     families of local tail/ spectral tail processes. The latter are utilised to show that   tail measures $\nu$  possess  a stochastic representer $Z$ such that $\nu=\nu_Z$ as defined in 
\eqref{repNU}.\\
\eEE{As in \cite{kulik:soulier:2020}, RV of general measures $\nu$ on $(\spaceD, \borelD)$  is discussed in \Cref{sec:regvM} with respect to some properly localised  boundedness $\mcb$ on $\spaceD$. We show that tractable instances for the investigation of  RV of $\nu$ arise if $\mcb$ can be characterised by $\famHMAP$ as in \Cref{item:condB} below, which is in particular the case for some common boundedness on the space of  general \cadlag\ processes, or on $\tilde l \setminus \tilde 0$ defined in \cite{MR3877549}[p. 877].\\
 	In \netheo{thmA} we relate the RV of \cadlag\ processes on $D(\R^l, \R^d)$ with the RV of their restrictions on $D(K,\R^d)$ for $K$ a given hypercube on $\R^l$. 
  Moreover, we present necessary and sufficient conditions for RV of \cadlag\ processes in \netheo{Th:RV}. 
  \cEG{Our findings show that RV of \cadlag\ processes can be investigated without imposing the stationarity assumption.}}
 
Besides being more complicated, the non-stationary case is also inevitably less tractable than the stationary one.  Despite those limitations, numerous  interesting results still continue to hold for non-stationary \cadlag\ processes, \vEE{including the equivalence of different definitions of RV and Breiman's lemma \cEG{(see e.g., \cite{AAP1579} for some extensions)}  with  its ramifications,}  see \Cref{sec:disc}. \\
\eEE{Another conclusion of this paper is that tractable cases arise for general Polish metric spaces  if the boundedness can be intrinsically  related to $\famHMAP$.} 
\cEG{The importance of our results is illustrated also by the wide range of potential applications and open problems  discussed in \Cref{openq}. Besides, our  results for local spectral tail processes, their relationship with tail measures and RV are of certain theoretical importance.}

 Below is a short summary of  some new aspects of this contribution: 
\begin{enumerate}[i)]
	\item We introduce tail measures, families of local tail/spectral tail  processes   for general measure spaces indexed by  $\famHMAP$; 
	\item  \label{lgS:b}All the results on the families of local tail and local spectral tail processes are new for the   settings of this paper. \neprop{pnu} is  new also  for the simpler cases   $\spaceD= D(\R^l, \R^d)$ or $\spaceD=D(\Zset^l, \R^l)$ and all $l\ge 1$;
	\item \netheo{prophope},  \netheo{thmA}  and the characterisation of  the limit measure  $\nu$ in \nelem{rusl}  are new also   for stationary $X$ taking values in $\spaceD$ as in \Cref{lgS:b}  above, whereas 
	\nelem{chans} is new if $\spaceD=D(\R^l, \R^d), l\ge 1$. Further \netheo{Th:RV} presents new results for $X$ with \cadlag\ sample paths  also when $X$ is stationary and $l>1$;
		\item \cEG{Our applications in \Cref{openq} include novel results for the  tail behaviour of supremum of regularly varying \cadlag\ processes. }
\end{enumerate} 
 
The paper is organized as follows:  \Cref{sec:prim} introduces notation and exhibits some  preliminary results  \cEV{concluding with} our  main assumptions.     Tail  measures, local tail/ spectral tail processes and stochastic representers   are discussed in \Cref{sec:HM}, whereas the  RV  of Borel measures and random processes is treated in \Cref{sec:regvarinD}.   \Cref{sec:disc} is dedicated to discussions and some extensions. 
\cEG{Potential applications, results for max-stable and $\alpha$-stable processes as well as open problems are presented in \Cref{openq}.}    All proofs are relegated to \Cref{sec:proofs}. In Appendix  we \cEV{review} \cEE{some properties} of general \cadlag\ functions  
and then display the mapping theorem.

\section{Preliminaries}
\label{sec:prim}
We present  first  several definitions and notation related to a given metric space. Then we continue with properties of a properly localised boundedness $\mcb$ 
and conclude this section with the formulation of the main assumptions.

\subsection{Measurable cones and the family of maps $\famHMAP$}
\label{fun:space}
 
Let $(\spaceD, \Dtau) $ be a metric space with corresponding Borel $\sigma$-field $\borel{\spaceD}$ and let  $\AA$ be another generic $\sigma$-field  on $\spaceD$.  In order to define a homogeneous measure on $\AA $ that satisfies \Cref{L:P1} we shall assume that a  pairing  
$$(z,f) \mapsto z \DOT f\in \spaceD, \  \   \   f\in \spaceD,z\in \RP=(0,\IF) $$
 (thus $\spaceD$ is a cone for the outer multiplication  $\DOT\ $) is a group action of the product group $(\RP, \cdot  \, )$ on $\spaceD$. This simply means  
   $$1 \DOT f =f, \  \   \ (z_1z_2)\DOT f=z_1\DOT (z_2 \DOT f) \in \spaceD, \   \forall f\in \spaceD,\forall z_1,z_2\in \RP .$$ 
\BD \cEE{ 
	We shall call  $(\spaceD, \DOT\ , \AA, \RP, \cdot)$ a measurable cone, if $\spaceD$ is non-empty and the  
corresponding group action  $(z,f) \mapsto z \DOT f, z\in \RP, f\in \spaceD$  of $(\RP, \cdot \ )$ on $\spaceD$ is   $\borel{\RP} \times \AA/ \AA$ measurable.
}
\label{defMC}
\ED
 In some cases  $\spaceD$ possesses   a  
zero element $0_\spaceD$, \ie\ 
$$z \DOT 0_\spaceD= 0_\spaceD, \ \  
\forall z \in \RPP=[0,\IF).$$
In the following we shall   write $0$ instead of $0_\spaceD$; abusing slightly the notation $0$  will also denote the origin of $\R^m,m\in \Nset$.
  
 Hereafter  $\TO= \{t_i, i\in \Nset\}$ is a non-empty subset of a given parameter space \cEE{$\TTT$}.

\BD \mb{We introduce the maps} $\norm[t]{\cdot}: \spaceD\mapsto [0,\IF] ,t\in \TO	 $, \mb{which} satisfy 
$$\norm[t]{z\DOT f}=z\norm[t]{f}, \  \forall f\in \spaceD,\forall z\in \RP$$
 and are   $\AA/\borel{[0,\IF]}$-measurable. \eEE{Suppose further that for all $t\in \TO$ there exists $ f\in \spaceD$ such that $\norm[t]{f} \in (0,\IF)$ and denote by $\famHMAP$ the family of the maps $\norm[t]{\cdot}, t\in \TO$.}
\label{defEL}
\ED 
In the sequel  we shall assume that  $\famHMAP$ is non-empty. Next,  given $f\in \spaceD$  and $K \subset \TTT$, we define
$$
f_{K}^* = \max_{ t \in K \tE{\cap \TO} } \norm[t]{f}.
$$ 
If $K\cap \TO=\emptyset$, interpret $f_{K}^*$ as  0 and write simply $f^*$ if $K=\TO$.

\def\Hh{\mathfrak{H}}
\cEE{Hereafter  $\Hh$ shall  denote the class of all   maps $\Gamma: \spaceD  \mapsto \R$  and all    maps $\Gamma: \spaceD  \mapsto [0,\IF]$
 which are $\AA /\borel{\R}$ and $\AA/\borel{[0,\IF]}$ measurable, respectively.  Write $\Hh_\lambda,  \lambda\ge 0$ for the class of  maps 
$\Gamma\in \Hh$ such that for all $f \in \spaceD$ and some $c>  0$,  $\Gamma(c \DOT f) = c^\lambda \Gamma( f)$.\\ 
}

\subsection{Boundedness on Polish spaces and $\mcb$-boundedly finite   measures}
\label{sec:boundV}
Consider   a non-empty set $\VV$ equipped with a $\sigma$-field $\AAV $.    
\BD A measure $\tailmeasure$ on $\AAV$  is a countably additive set-function  $\AAV \mapsto [0,\IF]$ with $\tailmeasure(\emptyset)=0$.   
We call $\tailmeasure$ non-trivial  if  $\tailmeasure(A)\in (0,\IF)$ for some $A\in \AA$ and denote the set of non-trivial measures on 
 $\mathscr{V}$ by $\mM(\AAV)$.   If $\AAV=\borelV$,  then  $\nu$ is called Borel.
 \ED
 
 Suppose next that $(\VV, \dV)$ is a Polish metric space and set  $\AAV= \borel{\VV}.$ 
Write   $\overline{A}$ and $\partial A$ for the closure and the topological frontier (boundary) of a non-empty set $A \subset \VV$, respectively. Note that $\partial A=\overline{A} \setminus int(A)$, where $int(A)$ is the interior of $A$. If $\tailmeasure\in \mM(\AAV)$, then the events (\ie\ the elements of $\AAV$) of interest are $A \in \AAV'$, where $\AAV'$ consists of  all events   such that $\tailmeasure(A) < \IF$. 
Since  $\AAV'$ is in general too large, reducing it to a   countably generated set is of great advantage for dealing with properties of $\tailmeasure$.  This motivates the concept of the properly localised boundedness which is quite general and not restricted to Polish spaces; our definitions below    are essentially taken from \cite{kulik:soulier:2020}[Appendix B], see also 
\cite{kallenberg:2017,MR3979308}.

\BD  A non-empty class $\mcb =\{ A:  A \subset \VV\}$   is called a properly localised boundedness on $\VV$ if
\begin{enumerate}[{B}1)] 
	\item \label{B3} $\mcb$ is closed with respect to finite unions and the  subsets of elements of $\mcb$ belong to $\mcb$; 
	\item  \label{shall} 	There exist  open sets  $O_n\in \mcb,n\in \Nset$ such that $\overline O_n \subset O_{n+1},n\in \Nset$  and $\bigcup_{n=1}^\IF O_n=\VV$. Moreover for all $A \in \mcb$ we have  $A \subset O_n$ for some $n\in \Nset$.
\end{enumerate} 
\ED

\BRM \vEE{A properly localised boundedness $\mcb$ contains the compact sets of $\spaceD$, see \cite{kulik:soulier:2020}[Rem B.1.2]. Moreover, all metrically bounded sets of $\spaceD$ form a localised boundedness and also the converse holds, namely if $\mcb$ is a properly localised boundedness, then there exists a metric $d'$ on $\spaceD$ for which $(\spaceD, d')$ is complete and $A \in \mcb  \iff$   $A$ is metrically bounded for $d'$, see \cite{kulik:soulier:2020}[B.1.3], \cite{MR3979308}[Rem 2.7].}
\ERM
\cEG{Throughout the following  $\mcb$ denotes a  properly localised boundedness on $\spaceD$.} 
\BD A  Borel measure  $\nu $  on $\AAV$ that satisfies $\nu(A)<\infty$ for all $ A \in   \mcb \cap \AAV$
is called $\mcb$-boundedly finite. If further 
$\nu$ is non-trivial, then  we write  $\nu \in \MBV$.
\ED 
  
If  $F$ is a closed subset of $\VV$  we set  $\VV_F= \VV\setminus F$ (assumed to be non-empty), which is again  a  Polish space.  Write $\mcb_F$ for the collections of subsets of $\VV_F$  with elements $B$ such that 
$$ \dV(x,F)= \inf_{f\in F} \dV(x,f) > \ve, \quad \forall x\in B$$
for some  $\ve>0,$ which may depend on $B$.   We can equip $\VV_F$  with a metric $d_{\VV_F}$, 
which induces  the trace topology on $\VV_F$ and the elements of $\mcb_F$ are metrically bounded. One instance is the metric given in \cite{kulik:soulier:2020}
[Eq.\ (B.1.4)]. In view of \cite{kulik:soulier:2020}[Example B.1.6] $\mcb_F $ is a  properly localised boundedness on  $ \VV_F$. 
In the particular case $F=\{a\}$ we write simply $\mcb_a$ and $ \VV_a$, respectively.\\
The boundedness $\mcb_F$, for $F$ being further a cone, appears in connection with hidden regular variation, see \eg\  \cite{MR4036028}, whereas $\mcb_0$ is the common boundedness used in the definition of RV, see \eg\ \cite{hult:lindskog:2005,SegersEx17} and references therein.  \\

Hereafter  the support of  $H\in \Hh$ is denoted by $\operatorname{supp}(H)$  
	defined by $\operatorname{supp}(H)= \overline{ H^{-1}((-\IF,0) \cap (0,\IF] )}$.
Suppose next that   $(\VV,\, \DOT\, , \AAV, \RP, \cdot)$  is a measurable cone 
and consider the following 
restrictions for a given  properly localised boundedness $\mcb$    on $\VV$:

\begin{enumerate}[$B$3)]
	\item  \label{item:boundProd}  
	{For all $A\in \mcb$ and all $z\in \RP$ we have $z\DOT A \in \mcb$;}
\end{enumerate}
\begin{enumerate}[$B$4)]
	\item \label{item:leben}
	There exists an open set $A\in \mcb$ such that $ z\DOT A \subset A$ for all $z> 1$. Suppose further that   
	$t \DOT \overline{A} \subset s\DOT A, \forall t> s >0$ and $\cap_{s\ge 1} (s \DOT A) $ equals  the empty set $\emptyset$;
\end{enumerate}
\begin{enumerate}[$B$5)]
	\item \label{item:condB}
\eEE{If $\famHMAP$ is a family of  maps as defined in Definition \ref{defEL},}
 then    	$A \in \mcb$ if and only if there exists some index set 
 $K_A \subset \TTT$ and $\ve_A>0$ such that $\forall f\in A$ we \zhora{have}{}  
	$$f^*_{K_A} = \cEE{\sup}_{t\in  K_A \cap \TO} \norm[t]{f}> \ve_A.$$
\end{enumerate}
Given   $\Gamma\in \Hh$ and a measure $\tailmeasure$ on $\AAV$, write
$$
\tailmeasure[\Gamma]= \int_{\VV} \Gamma(f) \tailmeasure(df).
$$ 
\BRM If $\tailmeasure \in \MBV$, then   $\tailmeasure$ is uniquely defined by 
$\nu[\Gamma]$ for all $\Gamma:\VV  \mapsto \R$ bounded continuous supported on $\mcb$. 
Moreover, $\nu$  is $-\alpha$-homogeneous, provided that 
$\tailmeasure[\Gamma_z]  = z^{-\alpha}\tailmeasure[\Gamma]$ for all bounded continuous $\Gamma \in \Hh$ with support in $\mcb$, 
		with $\Gamma_z(v)= \Gamma(z\, \DOT\, v), v\in \VV$ 
	    and \Cref{item:boundProd} holds. For more details we refer to \cite{kulik:soulier:2020}[Appendix B].
\label{lepuruxhi} 
\ERM
 \subsection{Main assumptions}  
Below we write $D(K, \R^d)$ for the space of functions $f: K\mapsto \R^d$. If $K=\R^l, K=(0,\IF)^l$ or $K$ is a hypercube of $\R^l$, then $f \in D(K,\R^d)$ is assumed to be a \cadlag\ function, 
see e.g., \cite{MR625374,JJ} for the definition in the less common case    $l>1$.   
Next, we  formulate the following set of assumptions: 
\begin{enumerate} [{$A$}1)] 
	\item  \label{A:1} 	$(\spaceD, \DOT\ ,  \AA,   \RP, \cdot )  $ is a  measurable cone, 
	 $\TO=\{t_i, i\in \Nset\}$ is a subset of some parameter space $\TTT$  and the family of maps $\famHMAP$  exists;  
	\item   \label{A:5}   $(\spaceD,\Dtau)$ is  a Polish  space with  a properly localised boundedness $\mcb$. Further,   
	\vEE{$\norm[t]{\cdot}$}'s 
		are finite  
		and  \Cref{A:1},\Cref{item:boundProd} hold;  
		\item   \label{A:2}   
		Let \mb{$\spaceD=D(\TTT, \R^d)$  with $\TTT=\R^l$ or $\TTT=\Zset^l$ equipped with the} Skorohod $J_1$ topology and the corresponding metric $\Dtau$ which turns it into a Polish space.   		
Set   $\norm[t]{f}=\norm{f(t)}, f\in \spaceD, t\in \TTT,$  where 
		$\norm{\cdot}: \R^d \mapsto [0,\IF)$ is a norm on $\R^d$.  
		 Here $\TO$ is a countable dense subset of $\TTT$,  the  pairing  $(z,f) \mapsto z \,\DOT f= zf$ with $(zf)(t)=  z f(t),t\in \TTT$ is the canonical  pairing. 
\end{enumerate} 

Under \Cref{A:2} the assumption \Cref{A:5} holds for $\spaceD=D(\R^l, \R^d)$,  which follows from \netheo{sam},\Cref{tev:0}-\Cref{tev:3}. 
  Moreover by  \netheo{sam},\Cref{tev:2}, the Borel $\sigma$-field  $\borelD$ agrees with $\AAK=\sigma(\proj_t, t\subset \TO)$. Consequently, 
  $\norm[t]{\cdot}, t\in \TO$ are $\borelD/\borel{\R}$ measurable and $1$-homogeneous and thus $\famHMAP$ exists.\\
Consider next the boundedness $\mcb_0$ defined on $\spaceD_0=\spaceD \setminus \{0\}$ with $\spaceD=D(\R^l, \R^d)$ and $0\in \spaceD$ the zero function. 
In view of \netheo{sam},\Cref{{tev:Mreis}} $A \in \mcb_0$  
if and only if there exists  a  hypercube  $K_A\subset \R^l$  and some $\ve_A>0$ such that 
\bqn{   \sup_{t\in K_A} \norm[t]{f} = f^*_{K_A} > \ve_A , \forall f\in A.
	\label{heart}
}

\BRM \vEE{Eq. \eqref{heart} shows that $\mcb_0$ satisfies \Cref{item:condB}. This is also the case if $\spaceD=D(\Zset^l, \R^d),$  see 
	\cite{kulik:soulier:2020}[p.\ 105] for $l=1$. Note that other properly localised boundedness satisfying \Cref{item:condB} exist, for instance   $\mcb_{\tilde 0}$ on the space $\tilde\spaceD_{\tilde 0}= \tilde l \setminus \{ \tilde 0\}$ defined in \cite{MR3877549}[p. 877] by metrically bounded sets therein.}
\ERM 

\section{Tail  measures}
\label{sec:HM} 
Tail  measures introduced in \cite{owada:samorodnitsky:2012} play a crucial role in the study of RV, see  e.g., \cite{owada:samorodnitsky:2012, Hrovje, klem,kulik:soulier:2020, PH2020}. \eEE{In the literature so far the main emphasis has been on shift-invariant tail measures and tail measures defined on  product spaces.}  In this section we shall assume   that  \cEG{\Cref{A:1} holds} and fix some $\alpha>0$.

\subsection{Definition and basic properties}
\label{sec:A} 
 If $\tailmeasure$ is a $-\alpha$-homogeneous measure on $\AA$ and $A\in \AA$ satisfies   $z\DOT A=A$ for some positive $z\not=1$,  then  		
\begin{equation}\label{lemOhne}
	\tailmeasure(A)\in \{0,\IF\}.
\end{equation}
By the $1$-homogeneity of the  maps $\norm[t]{\cdot}$,  \eqref{lemOhne}  \vEE{implies  that    $F_*$ defined by }
$$
F_*=   \{f\in \spaceD : f^*_{\TO} = 0 \} ,\ \ f^*_\TO=\cEE{\sup}_{t\in \TO} \norm[t]{f}
 $$
satisfies $\tailmeasure(F_*) 
\in \{0,\IF\}$. 
Of particular interest are measures $\nu$ such that 
\begin{enumerate}[{M}1)]
	\item \label{L:P3} { $\tailmeasure(F_*)=0$},
\end{enumerate}
since this property is crucial for establishing their $\sigma$-finiteness. 
Next, we define    tail measures on $\AA$, supported by  the findings of \cite{owada:samorodnitsky:2012}, in which  tail measures on the product $\sigma$-field of $\spaceD=(\R^d)^\TT$ are introduced. See also \cite{Hrovje, kulik:soulier:2020, klem} for special product spaces containing a zero element 0 (we do not assume existence of 0 here) and \cite{PH2020} for  $\spaceD=D(\R, \R^d)$.  

\BD A  measure $\nu$ on $\AA$ that satisfies \Cref{L:P1},\Cref{L:P3} is called a tail measure (write $\nu \in \mSD$) if     
\begin{enumerate}[{M}2)]
	\item \label{L:P2} 
	$p_h:= \tailmeasure(\{f\in\spaceD:\norm[h]{f}>1\})  \in [0,\IF),   \tE{\forall h\in\TO}$, with $p_{h_0}\tE{\in (0, \IF)}$ for some $h_0\in \TO$.
 \end{enumerate}
\ED 
  
\BRM    
	The measurability of 
$\norm[h]{\cdot}, h\in \TO$  implies   $A_h=\{ f\in \spaceD: \norm[h]{f}=1\} \in \AA, h\in \TO$.  If $\tailmeasure\in \mSD$, then 
{ by \Cref{L:P1}
	and \Cref{L:P2} 
 $\tailmeasure(A_h)=0$  for all $h \in \TO$.
}
\eEE{Consequently,  for all $x>0$
\begin{equation}\label{thus}
	 p_h(x)=\tailmeasure(\{ f\in \spaceD: \norm[h]{f} \ge x\}) = 
 x^{-\alpha}  \tailmeasure(\{ f\in \spaceD: \norm[h]{f} > 1\})= 
x^{-\alpha} p_h, \ h \in \TO
\end{equation}
and thus if $p_h=0$, then $p_h(0)=0$ follows by the countable additivity of $\tailmeasure$. Since  \Cref{L:P2}  and \eqref{lemOhne} 
imply 
\begin{equation}\label{gurt}
	\nu(\{f\in \spaceD: \norm[h]{f}=\IF \})=0, \ \forall h \in \TO,
\end{equation}
then \Cref{L:P3} is equivalent to   
\bqn{\label{eM2}
\tailmeasure\Bigl(\Bigl\{f\in \spaceD :\sup_{t\in \TO} \norm[t]{f} \in \{ 0,\IF\} \Bigr\}  \Bigr)=
\tailmeasure\Bigl(\Bigl\{f\in \spaceD :\sup_{t\in \TO: p_t>0} \norm[t]{f} \in \{0,\IF\} \Bigr\}  \Bigr) =0.
}	
}
\label{simonDiv}
\ERM
 
 \cite{owada:samorodnitsky:2012}[Prop 2.4] derives necessary and sufficient conditions for the $\sigma$-finiteness of tail measures defined on the 
product $\sigma$-field of  $\spaceD=(\R^d)^\TT$.   Our definition of tail measures implies their  $\sigma$-\eEE{finitness} and as in \cite{klem}[Prop. 2.3] we have the following result (its proof is omitted). 
 
\BEL \label{triv}  
	If $\tailmeasure\in \mSD$, then it  is  $\sigma$-finite 
 and $\tailmeasure$ is uniquely determined by its restrictions to  \eEE{\mb{$\{ f\in \spaceD : \norm[h]{f}> 1\})$}} for all $h\in \TO$.
 \EEL 
 
Recall that  $ Z$ denotes throughout this paper a $\spaceD$-valued   random element defined on a complete probability space $\Pspace$. Suppose next that  
$ \norm[h]{Z}$ are random variables (rv's) for all $h\in \TO$ and  further   
\bqn{ \label{harixh}
	\E{\norm[h_0]{Z}^\alpha} {\in (0, \IF)}, \quad 	\E{\norm[h]{Z}^\alpha}\in [0,\infty), \quad { \forall h\in \TO}, 
	\quad {\pk{ Z_{\TO}^*\not=  0  }=1}
}
for some  $h_0\in \TO$ and consider  $\tailmeasure_{\vk Z}$ defined in \eqref{repNU}. Since $\spaceD$ is  a measurable cone, then  $\tailmeasure_Z$ is the image  	 measure of $(z, f) \mapsto z \DOT f$ with respect to the product measure $\mu(df) \times \nualpha(dr)$, with 
$$\mu= \mathbb{P} \circ  Z^{-1}, \quad \nualpha(dr)= \dualpha.$$
Clearly, $\tailmeasure_Z$  satisfies \Cref{L:P1}-\Cref{L:P3} with  \mb{$p_h=\E{ \norm[h]{Z}^\alpha }< \IF$} for all $h \in \TO $ and hence 
$\tailmeasure_Z\in \mSD$.\\
Hereafter $R$ is an  $\alpha$-Pareto rv with $\pk{R> t}= t^{-\alpha}, t\ge 1$, independent of all other random elements. It can be utilised to link $Z$ and $\tailmeasure_Z$ as in \eqref{mada0} below.

If  a measure $\tailmeasure$  on $\AA$  has  representation \eqref{repNU} with  $Z$    satisfying the first two conditions in \eqref{harixh}, then  it follows   that   for all  $h \in \TO, \Gamma\in \Hh, \ve \in (0,\IF)$ (here 
$\Hh$ is the class of maps defined in \Cref{fun:space})
\bqn{
	 \int_{\spaceD}  \Gamma(f) \ind{\norm[h]{f}> \ve } \nu(df) =
	 \frac 1 {\ve^{\alpha}} \E{ \norm[h]{Z}^\alpha  \Gamma((\ve R/ \norm[h]{Z })  \DOT  Z) }
	\label{mada0}
}  
and hence 
$$p_h= \tailmeasure( \{f\in \spaceD: \norm[h]{f}> 1  \}) = \E{\norm[h]{Z}^\alpha} \in [0,\IF),\  \forall h \in \TO.$$

\subsection{Local tail and  local spectral tail processes}
 We introduce next the local tail and local spectral tail processes as in \cite{klem};    \eEE{our setup here is less restrictive compared to that  of product spaces dealt with in the aforementioned paper}. Recall that  $\tailmeasure\in \mSD$ stands for $\tailmeasure$ is  a  tail measure on $\Dspace$.
 \BD  Given $\tailmeasure\in \mSD$ and  $h\in \TO$ such that $p_h>0$,  the local  process $ Y^{[h]}_\nu$ of $\nu$ at $h$   has law  $\nu_h(A)=\nu( \{ f \in A: \norm[h]{f}> 1 \})/p_h, A \in \AA$. 
  We call $\vk \Theta^{[h]}_\nu = ( \norm[h]{\vk Y^{[h]}_\nu})^{-1} \DOT\vk Y^{[h]}_\nu$  the  local spectral tail process of $\nu$ at $h$. 
  If $p_h=0$, then  set $Y^{[h]}_\nu=R \DOT g, \Theta^{[h]}_\nu=g$ with $g\in \spaceD$ satisfying $\norm[h]{g}=1$.
\ED 
 
We shall drop the subscript $\tailmeasure$ for  local tail/  spectral tail processes, when there is no ambiguity.  
\BRM 
\label{shurdh}
{    $Y^{[h]}$ is a random element from  a probability space $(\Omega, \mathscr{F}, \mathbb{P})$ to $\Dspace $ and similarly for $\Theta^{[h]}$. Take for instance  $\Omega=\spaceD, \mathscr{F}=\AA, \mathbb{P}=\nu_h$ and define $Y^{[h]}:  f \mapsto \ind{\norm[h]{f} >1 } f, f\in \spaceD $, which in view of \Cref{A:1} is a $\mathscr{F}/\AA$ measurable map for all $h\in \TO$. The assumption on $\norm[t]{\cdot}$  and \eqref{gurt} imply that $\norm[t]{Y^{[h]}}, t\in \TO$ is a   non-negative   rv for all $h,t \in \TO$.}
\ERM

\begin{proposition}  If  $\tailmeasure\in \mSD$, then  for all   {$h\in \TO$}  
\bkn{  
	\label{kksuter} 
	\pk{ \norm[h]{  \bY^{[h]}}>1 }=	\pk{ \norm[h]{  \bTheta^{[h]}} =1  }=1
} 
and 
\eEE{if $p_h=0, p_t>0$, then  $\norm[h]{ Y ^{[t]}}=
	\norm[h]{ \Theta ^{[t]}}=0$ almost surely.} 
Further, for all {$h,\,t\in \TO$} such that $p_hp_t>0$  
\bqn{  
	p_h\E{  \norm[t]{\tPr}^\alpha \Gamma(  \bTheta^{[h] })    } &=&
	p_t\E*{ \ind{ \norm[h]{  \bTheta^{[t]}}  \not= 0} \Gamma(    \bTheta^{[t]} )}, \quad \forall \Gamma\in \mb{\Hh_0},
	\label{tsf}
}	
and for all $x>0$ 
\bqn{ \label{tiltY}
	p_h \E*{ \Gamma( x  \DOT \bY^{[h]}  ) \ind{ x\norm[t]{   \bY ^{[h]}}> 1}} 	 = p_t x^{\alpha} 
	\E*{  \Gamma(   \bY^{[t]}) \ind{ \norm[h]{  \bY^{[t]}}> x}}, \quad \forall \Gamma\in \mb{\Hh}.
} 
Moreover, \eEE{the law of $Y^{[h]}$ agrees with  that of} $R \, \DOT\, \bTheta^{[h]},h\in \TO$ and  
$Y^{[h]}, \eEE{h\in \TO: p_h>0}$ uniquely determine $\tailmeasure$.
\label{pnu}
\end{proposition}

\BRM For all $h,  t\in \TO$ such that $p_hp_t>0$
\bqn{  	p_h \E*{\ind{ \norm[t]{\vk \Theta^{[h]}}  \not= 0}  \Gamma( \vk \Theta^{[h]})  }
	&=& 
	p_t \E*{ \ind{ \norm[h]{\vk \Theta^{[t]}}  \not= 0} \Gamma(  \vk \Theta^{[t]} )}, \quad \forall  \Gamma\in \mb{\Hh_\alpha},
	\label{tsf2}
}	
 see also \cite{MR4153591}.
\label{hofi2}
\ERM

If $\nu=\nu_Z$ is given by   \eqref{repNU} with  $ Z$  satisfying \eqref{harixh}, 
then by definition   the claim in \eqref{mada0} implies  
\bqn{  
	 \int_{\spaceD}  \Gamma(f) \ind{\norm[h]{f}> 1 } \nu(df)= p_h \E{  \Gamma(Y^{[h]})} =  \E{\mb{\norm[h]{Z}^\alpha}  \Gamma(( R/ \norm[h]{Z })  \DOT  Z) }, \  \forall  h \in \TO: p_h>0, \forall \Gamma\in \Hh.
		\label{mada}
} 
If  $p_h=0$, then  \eqref{mada} still holds taking $\Gamma$ to be bounded.  Consequently, since $p_h=   \E{ \norm[h]{Z}^\alpha}< \IF$, then 
$Z$ determines the laws of $\bY^{[h]}$  and $\bTheta^{[h]}$ 
denoted by $\mathbb{P}_{\bY^{[h]}}$ and $\mathbb{P}_{\bTheta^{[h]}}$, respectively, \ie\  
\bkn{ \mathbb{P}_{\bY^{[h]}}(\cdot)= p_h^{-1}\E{ \norm[h]{Z}^\alpha \delta_{ (R /\norm[h]{Z}) \ \DOT \ Z} (\cdot)}, \quad \mathbb{P}_{\bTheta^{[h]}}(\cdot)
	= p_h^{-1}\E{ \norm[h]{Z}^\alpha \delta_{ \norm[h]{Z}^{-1} \ \DOT \ Z} (\cdot)}
	\label{lawY}
}
for all  $h\in \TO$ such that $p_h>0$, with $\delta_{x}(\cdot)$ the Dirac point measure of $x\in \R$. 

\def\equallaw{\stackrel{d}{=}}

\subsection{Stochastic representers}
\label{sec:spectral-representation}
A  measure $\tailmeasure$ on $\AA $ has a stochastic representer $Z$  satisfying  \eqref{harixh} if  $\nu$ equals $\nu_Z$ defined in \eqref{repNU}.  Hereafter $q_t\ge 0, t\in \TTT$ \eEE{satisfy $ q_t>0$ for all $t\in \TO$ such that $p_t>0$ } and we set  
	$$\SSY(Y^{[h]})= 
\int_\TO  \norm[t]{ Y^{[h]} }^\alpha q_t \lambda(dt) ,  \quad h\in \TO,$$
where \mb{$\lambda(dt)=\LK(dt)$} is the counting measure on $\TO$.  

If $q_t$'s are such that   $\sum_{t \in \TO}q _t =1$, then we shall consider a $\TO$-valued rv $N$ with probability mass function    $q_t, t\in \TO$ being independent of all other elements.\\
\BRM \label{remNA}
In view of Remark \ref{shurdh} and \cite{Kallenberg}[Cor. 5.8] it is possible to choose $Y^{[h]}, h\in \TO $ and $N$ to be defined in the same probability space $\Pspace$ such that all these are independent, which we shall assume below. Moreover, \eEE{$q_t$'s and thus $N$} can be chosen such  that $\E{p_N}< \IF$, with  $p_h=\E{\norm[h]{Z}^\alpha}< \IF$.   Recall $Y^{[h]}=Y_\nu^{[h]}$ is defined with respect to some 
 $\nu\in \mSD$.
 \ERM   
\eEE{Below, let  $\TTTK=\{K_n\subset \TO,n\in \Nset\}$ such that 
$\cup_{n\ge 1} K_n= \TO.$
} 
\vEE{Under assumption \Cref{A:2} $\TO$ is  a dense subset of $\TTT$ and we shall choose $K_n= [-n,n]^l \cap \TO, n\in \Nset$.} 
\BD A   measure $\tailmeasure$ on $\AA$  is  $\TTTK$-bounded  (compactly bounded when \Cref{A:2} holds) if    
\begin{enumerate}[{M}3)]
	\item \label{L:P4} 		$ 	\tailmeasure \left(\left\{f\in \spaceD : f_{K}^* >	1 \right\}\right) < \infty, 
	\ \forall K \in \TTTK$.
\end{enumerate}
\ED

\mb{The} next result shows that tail measures have a family of representers $Z$,  which can be utilised to define $Y^{[h]}$ and $\Theta^{[h]}$ via \eqref{lawY} and 
to give  an equivalent condition for  \Cref{L:P4}.   

\BEL
\label{Th:Rep} 
If $\tailmeasure\in \mSD$, then   
$\nu$ has  stochastic  representer $Z=Z_N$ given by 
\bqn{\label{spectr}
	Z_N= 		\frac{p_N^{1/\alpha} \DOT  \bY^{[N] }}{( \SSY(\bY^{[N]}))^{1/\alpha}},
 } 
where  the    local tail processes $Y^{[h]} , h\in \TO$  and $N, q_t,t\in \TO$ are as in Remark \ref{remNA}. 
\eEE{Further, $\norm[h]{ Z}=0$ if $p_h=0$}  and 
$\tailmeasure$ satisfies \Cref{L:P4}  if and only if   
\bqn{ 
	\E*{\sup_{t\in K \tE{\cap \TO}} \norm[t]{Z}^\alpha}< \IF, \quad \forall K \in  \TTTK. 
	\label{borsh}
} 
\end{lemma}

\def\TTTKN{\widetilde{K}(\TTT)}
\BRM  
\begin{enumerate}[(i)]
		\item  The claim that  $\nu\in \mSD$   is specified by 
	\eqref{repNU} for some representer $Z$  follows also by \cite{MolchEvans}[Thm 1]. 
\item 
 \cEE{ If $\tailmeasure=\tailmeasure_Z=\tailmeasure_{\tilde Z} \in \mSD$,  applying \eqref{mada} 
$ \forall h\in \TO: \eEE{p_h>0}$ we obtain 
 	\bqn{ \label{spd}
  p_h \E{  \Gamma( \Theta^{[h]}) }=		\E{ \norm[h]{ Z}^\alpha \Gamma((R/\norm[h]{Z}) \DOT\vk Z)} = 		\E{ \norm[h]{ Z}^\alpha \Gamma(\vk Z)} =   	\E{ \norm[h]{ \tilde  Z}^\alpha  \Gamma(\tilde  Z)}	,
 		\quad \, \forall \,  \Gamma\in \Hh_0.
}		 
 Since  by \eqref{harixh} we can choose $\eEE{q_h}>0,h\in \TO$ such that 
	$S(Z)= \sum_{h\in \TO} q_h\norm[h]{Z}^\alpha \in (0,\IF)$ 
	almost surely, then  (\ref{spd}) implies 
\bqn{
\E{ \Gamma_\alpha(\vk Z)} = 
	\sum_{h\in \TO}   	q_h\E*{ \norm[h]{Z}^\alpha \frac{\Gamma_\alpha(Z) }{S(Z)}}= \sum_{h\in \TO}  q_h 	\E*{ \norm[h]{\tilde Z}^\alpha\frac{  \Gamma_\alpha(\tilde Z)}{S(\tilde Z)}}
	 = \E{\Gamma_\alpha(\tilde Z)}, \ \forall \Gamma_\alpha \in \Hh_\alpha. 
}
Hence,  if in \nelem{Th:Rep} $\spaceD$ is a countable  product space or it is equal to $D(\R,\R^d)$ we retrieve  \vEE{ 
	\cite{klem}[Thm 2.4] and \cite{PH2020}[Thm 2.3]}, respectively.}
\item \vEE{ \nelem{Th:Rep}  
	together with \nelem{pnu} and 
\eqref{lawY}  implies \cite{klem}[Prop 2.7].}
\end{enumerate}
\ERM 
 
\begin{example} Assume that \Cref{A:2} holds.
	Denoting by  $\zD $ the zero function we have 
	$\{\zD\} \in \borel{\spaceD}$. Since $\norm[t]{f}=\norm{f(t)}, t\in \TO$ with $\norm{\cdot}$ a  norm on $\TT$  and further $\TO$ is a dense subset of $\TTT$, then 	\Cref{L:P3} is equivalent to 
	$$\tailmeasure(\{\zD\})=0.$$
	If $\tailmeasure\in \mSD$, then its representer $Z$ is a random process
	with almost surely \cadlag\ sample paths  and so are both $Y^{[h]}$ and $\Theta^{[h]}$ for all $h\in \TTT$.	\\
	A direct implication of (\ref{lawY}) is that if $\tailmeasure$ is shift-invariant (see \cite{klem} for definition), then by \eqref{mada} we have the equality in law 
	\bkn{ \label{shiftY} 
		Y^{[h]}  \equallaw  B^h Y^{[0]}, \quad p_h=p_{h_0}\in (0,\IF), \quad \forall h\in \TTT,
	}  
	where $B^h f(\cdot)= f(\cdot -h), h\in \TTT$. 
Note in passing that in this case \eqref{tiltY} reads (set below $Y= Y^{[0]}$)
	\bqn{  \label{eqR1} \E{ \Gamma( xB^h  Y) \ind{ x \norm{  Y(-h)} >1 }} = 
		x^\alpha \E{  \Gamma(   Y) \ind{\norm{  Y(h)} >x }   },\  \forall \Gamma \in \Hh,  \forall h\in \TT, \forall x>0,
	}  
	which for $\TTT=\Zset$ is stated initially in \cite{Hrovje}, see \cite{PH2020} for the case $\TTT=\R^l$ and   \cite{Planic} for other interesting properties of $Y$. \cEG{Also the converse holds, i.e., \eqref{shiftY} implies that  $\nu$ is shift-invariant. 
}	
	\label{exAA}
\end{example}

\def\intTKP{ \int_{K}}
\def\HH{\Hh}
\def\tHH{\widetilde{\HH}}

\subsection{Constructing $\tailmeasure$   from local tail processes}
\label{sec:TY}
A given tail measure defines the family of local tail processes $Y^{[h]}, h\in \TO$. We discuss in this section \mb{the inverse procedure, namely} how to calculate $\tailmeasure \in \mSD$ from the $Y^{[h]}$'s. Hereafter $q_t's$ are positive constants and we write  $\lambda$ for counting measure $\lambda_\TO$ on 
$\TO$ or  for Lebesgue measure on  $\TTT=\RL$  
and define  
\bkn{ \label{defek} \exc_{K}(f)= \intTKP \ind{ \norm[t]{ f} > 1} q_t \lambda(dt), \quad f\in \spaceD, 
}
with $K$ a non-empty subset of $\TO$ if $\lambda=\lambda_\TO$ and $K$ is a non-empty hypercube  of $\R^l$, otherwise.\\ 
The next result    extends  \cite{kulik:soulier:2020}[Thm 5.4.2].  

\BEL Let  $\tailmeasure\in \mSD $  be given.   Suppose that for some  $H\in \Hh$, there exists  $\ve_H>0$ and some non-empty  $K_H\subset \TO$ such that for all $f\in \spaceD $ satisfying $f_{K_H}^*\le \ve_H$ we have $H(f)=0$.  
If  
	\mb{$\int_{K} p_t  q_t\lambda(dt) \in (0, \IF) $} for some  $K\subset \TO$  such that 
	$ K_H \subset K$, then 
	\bkn{
		\tailmeasure[H] &=& \ve ^{-\alpha}   \intTKP
		\E*{  \frac{H(\ve \DOT  Y^{[t]})}{\exc_{K}(  Y^{[t]}) }} p_t q_t \lambda(dt),\quad  \forall \ve \in (0,\ve_H].
		\label{eq:goodschl}
	}  
\label{lem:repA}
\EEL 

\BRM \label{mbab} 
\begin{enumerate}[(i)]
	\item \label{mbab:1} Taking \mb{$H(f)= \ind{ \sup_{s \in K } \norm[s]{f}>1}$} for some non-empty  $K\subset \TO$, for $q_t$'s as chosen in \nelem{lem:repA} 
	 we obtain from \eqref{eq:goodschl} that \Cref{L:P4} is equivalent to 
\begin{equation}\label{ifo}
	\intTKP
\E*{  \frac{H( Y^{[t]})}{\exc_{K}(  Y^{[t]}) }} p_t q_t \lambda(dt) = \intTKP \E*{  \frac{1}{\exc_{K}(  Y^{[t]}) }} p_t q_t \lambda(dt) < \IF, 
\ \forall K \in \TTTK
\end{equation}
since $\pk{\norm[t]{Y^{[t]}} >1}=1, \forall t\in \TO$ implies that 
$ H( Y^{[t]})=1$ almost surely for all $t\in K$;
	\item \label{mbab:0} Under \Cref{A:2}, by \eqref{heart} an $-\alpha$-homogeneous 
Borel measure $\nu$ on $\borel{D}$ is $\mcb_0$-boundedly finite if and only if \Cref{L:P4} holds, or equivalently \eqref{ifo} is satisfied. 
\end{enumerate}
\ERM

We have seen that local tail/spectral tail processes can be defined directly \mb{through the} representer $Z$ of a tail measure $\nu=\nu_Z$. We \mb{may also} define such families without referring to $Z$ as follows. 
\BD The family of $\spaceD$-valued random elements $Y^{[h]}, h\in \TO$ is called a family of tail processes with index $\alpha>0$, 
if for given 
weights $p_h\in [0,\IF), h\in \TO$, with  $p_{h_0}>0$ for some $h_0\in \TO$   both \eqref{kksuter}   and \eqref{tiltY} are satisfied and further 
\eEE{$\pk{\norm[h]{ Y ^{[t]}}=0}=1$  if $p_h=0, p_t>0$.}  
\ED 
Similarly, $\spaceD$-valued spectral tail processes $\Theta^{[h]}, h\in \TO$ can be defined  without reference to some measure $\nu$. 
As  shown next   a family of tail processes defines uniquely a tail measure on $\spaceD$. 

\begin{lemma}
	\label{capit2}
Let  $\bY^{[h]}, h\in \TO $  and $ \Theta^{[h]}, h\in \TO$ be a family of $\spaceD$-valued tail  and spectral tail processes, respectively  with index $\alpha>0$ \eEE{and  let $q_t$'s be positive constants such that $\sum_{t\in \TO}  p_t q_t \in (0, \IF)$.}  
\begin{enumerate}[(i)]
	\item \label{item:toniRA}{If  $Y^{[h]}, h\in \TO, N$ are defined as in Remark \ref{remNA},}
then  there  exists   a unique $\tailmeasure\in \mSD$   such that its local tail processes are $Y^{[h]}, h\in \TO$;
\item \label{item:toniRC} $\tailmeasure$ defined in \Cref{item:toniRA} above satisfies \Cref{L:P4} if and only if 
\bkn{ \label{maalp} 
	\intTKP  \E*{  \frac{1  } { \exc_{K}( \bsY^{[t]})} }  p_t  q_t \lambda(dt	)< \IF, \quad \forall K \in \TTTK ,
}  
\cEV{which under \Cref{A:2} with $\TTT=\R^l$ is true also for $\lambda(dt)$ the Lebesuge measure on $\TTT$ and all compact $K \subset \TTT$ with 
	$q_t$ as in \nelem{lem:repA};}
\item \label{plusT} 	$R \DOT  \Theta^{[h]}, h\in \TO$ is a family of tail processes with index $\alpha>0$.
\end{enumerate} 
\end{lemma}

\begin{example}
	\label{exa:1}
Suppose that \Cref{A:2} holds and let $Z$  satisfy  \eqref{harixh} for all $h\in \TTT$. Assume  that $\tilde Z=B^h Z$ and $Z$  
 satisfy  \eqref{spd} for all \mb{$h\in \TTT$}, which implies  
 \mb{$\E{ \norm[h]{Z}^\alpha}=C \in (0,\IF),  \forall h\in \TTT$} and 
  $\tailmeasure_Z$ has local processes at $h\in \TTT$ given by  
  \begin{equation}\label{localSatio}
  Y^{[h]}=	B^h Y, \ \ Y=Y^{[0]}  ,
  \end{equation}
 with $Y^{[0]}$ having law  \mb{$\E{ \norm[0]{Z}^\alpha \delta_{R Z/\norm[0]{Z}} (\cdot)}/C$}.   
 	 For $\tailmeasure_{\sigma Z}$  with   representer $\sigma Z, \sigma \in \spaceD, \sigma\not=0$  its local tail processes are \mb{given by}
 $$Y^{[h]}(t)= \frac{\sigma(t)}{\sigma(h)} B^h Y(t), \quad \forall t\in \TT, \forall h: \sigma(h)\not=0.
  $$
\end{example}

\section{RV  of measures and random elements} 
\label{sec:regvarinD}
 We \mb{first discuss the} RV  of Borel measures  and $\spaceD$-valued random elements assuming \Cref{A:5}. 
\mb{Subsequently, we study in detail} the  RV  of processes with \cadlag\ paths. 

\subsection{$\mcb$-boundedly finite Borel measures} 
\label{sec:regvM}

  RV  of Borel measures on Polish metric spaces is discussed in \cite{hult:lindskog:2005,hult:lindskog:2006,MR3271332,SegersEx17,MR4036028}. 
 We follow the treatment of  RV  in \cite{kulik:soulier:2020}, where some properly localised boundedness $\mcb$ plays a crucial role. \\
Throughout this section we suppose  \mb{that} \Cref{A:5} holds. \\ 
Let  next ${ \nu_z,z>0}$ be $\mcb$-boundedly finite  measures on $\borel{\spaceD}$ and  recall our notation $\nu[H]=\int_{\spaceD} H(f) \tailmeasure(df) $.
\BD     $\nu_z$ converges $\mcb$-vaguely  to some Borel measure $\nu$  as $z\to \IF$ (denote this by $\nu_z\convvague\tailmeasure$) if 
\begin{equation}\label{equivC}
	\lim_{z\to\infty} \nu_z[H] 	=\tailmeasure[H] 
\end{equation} 
is valid for all \mb{continuous and bounded} maps $H: \spaceD \to \R$ with $supp(H)\in \mcb$.
\ED   
In the sequel   $g,g'$ are two  maps $\RPP\mapsto \RPP$ and  for some $\mu \in \mM(\mcb)$ we set 
$$\mu_z(A)=g(z)\mu( z \DOT  A), \  \ \mu_z'(A)=g'(z)\mu( z \DOT  A), \  \  A \in \borelD, z\in \RP.$$  
\BEL 
If    $\mu_z$, $\mu_z'$ converge $\mcb$-vaguely   to  
$\tailmeasure\in \mM(\mcb)$ and $\tailmeasure'\in \mM(\mcb)$, respectively, as $z\to \IF$, then  
 $\limit{z} g'(z)/g(z)=c$ and $\tailmeasure'= c \tailmeasure$ for some $c\in (0,\IF)$.
\label{lem:kulikB22}
\EEL
 \def\RVP{\mathcal{R}(g, \mcb, \nu)}
\BD \label{defCel}
 $\mu \in \MBB$ is regularly    varying  with scaling  function  $g$, if   
$\mu_z \convvague \nu\in \mM(\mcb)$, abbreviate this as $\mu \in \RVP $.
\ED
\BRM \begin{enumerate} [(i)]
	\item \label{nock:1} In view of \cite{kulik:soulier:2020}[Cor B.1.19]),  $\mu \in \RVP $ if and only if there exist  open sets $O_k\in \mcb, k\in \Nset$ satisfying \Cref{shall}  such that for all positive integers $k$ 
$$\nu(\partial O_k)=0, \  g(z) \mu(z \DOT (O_k \cap \cdot)) \convweak  \nu(O_k \cap \cdot), \  z \to \IF.$$
\item \label{nock:2} \cEG{Let  $\tailmeasure_Z\in \mSD $   and write $\mathbb{P}_Z$  for the law of $Z$. The $-\alpha$-homogeneity implies that 
$\tailmeasure\in \mathcal{R}(g, \mcb, \nu_Z)$ with $g(x)= x^\alpha$. Note  that $\tailmeasure_Z$ is the mean measure of the Poisson Point Process (PPP) $ N$ on $\spaceD$, which is defined by
$$  
N(\cdot)= \sum_{i=1}^\IF \delta_{P_i Z^{(i)}}(\cdot),
$$
with $\sum_{i=1}^\IF \delta_{P_i,Z^{(i)}}$ being a PPP on $(0,\IF) \times \spaceD$ with mean measure $v_\alpha(\cdot) \odot P_Z(\cdot)$ and $Z^{(i)}$'s being independent copies of $Z$.}
\end{enumerate}
\label{nock}
\ERM 
Write next $g\in \RA$, if  
\mb{$$\limit{z} g(z t)/g(z)= t^\alpha, \ \ \forall t>0
	$$}for a non-negative rv $W$ we write $W\in \RA$ if $1/\pk{W > t} \in \RA$.
   
Set  $H_t(f):= \norm[t]{f}, f\in \spaceD$ and recall the definition of $p_h(x)$ in \eqref{thus}. 
 
\begin{lemma} Let  $\mu \in \RVP$, where  $g$ is Lebesgue measurable. 
\begin{enumerate}[(i)]
	\item \label{item:AL1} If   $g\in \RA$ for some $\alpha>0$, then $\tailmeasure$ is $-\alpha$-homogeneous;
	\item \label{item:AL2} If  $p_{t_0}(x) \in (0,\IF)$  	
	for some $t_0\in \TO,x>0$ and further $H_{t_0}^{-1}(B) \in \mcb$
	for all Borel set $B\in \borel{[0,\IF)}$  separated from 0 satisfying also $\tailmeasure(Disc(H_{t_0}))=0$,  
	then $g\in \RA$  for some $\alpha>0$  and 
 $\tailmeasure$ is $-\alpha$-homogeneous; 
	\item \label{item:AL3} Suppose that \Cref{item:leben} holds.  If  $\mu(k\DOT A)>0$ for almost all $k> M>0$ and 
	the group action is continuous,    then 
	$g\in \RA$  for some $\alpha>0$  and 	 $\tailmeasure$ is $-\alpha$-homogeneous.
\end{enumerate}  
\label{lem:engj}
\end{lemma}

\def\DIF{\spaceD^{(\IF)}}
\def\Dn{\spaceD^{(n)}}

\subsection{$\spaceD$-valued random elements} 
\label{sec:rvRelements}
 Consider next a $\spaceD$-valued  random element $X$ defined on a complete probability space $\Pspace $. Define for some  $g: \RPP \mapsto \RPP$ 
$$\mu_z(A)= g(z)\pk{ X \in z \DOT A }, \quad A \in \AA, z>0.$$
 
\begin{definition}
	\label{hypo:regvar-in-D}
\mb{The random element} $X$ is called RV  with respect to $g$ and $\tailmeasure\in \mM(\mcb)$,  if  
 $\mu_z\convvague \tailmeasure$ as $z\to \IF$.  Abbreviate this as $X \in \RVBA$ and when  $g\in \RA$ write 
 $X\in \RVB$. 
\end{definition}

If $X \in \RVBA,$ with Lebesgue measurable $g$, under the assumptions of \nelem{lem:engj},\Cref{item:AL2} 
we have that $\norm[t_0]{X}\in \RA$ implies   $\tailmeasure$ is $-\alpha$-homogeneous.  
If further  the conditions of \nelem{lem:engj},\Cref{item:AL2} hold for all $h\in \TO$, then    
 \netheo{mapping}   yields   
 \bkn{\label{klubge}
	\lim_{z\to\infty} \frac{ \pk{ \norm[h]{\bX}>z }}{ \pk{ \norm[t_0]{\bX }> z}}=\frac{p_h}{ \cEE{p_{t_0}}} \in [0,\IF), \quad \forall h\in \TO, 
	\quad p_{t_0}\in (0,\IF).
} 

\eEE{Under   \Cref{A:5}, we present  in the next result a sufficient condition for the  RV  of $X$ when 
$\mcb$ is determined by the countable  family of maps $\famHMAP$ as in \Cref{item:condB}.
}

\begin{theorem}
	Let   $X$  be  
 such that  $\norm[t_0]{X}\in \RA$ for some $t_0\in \TO$ and  \eqref{klubge} holds. Assume   that    $\forall h\in \TO: p_h>0$ 		conditionally on $\norm[h]{\bsX}>z$, $z^{-1} \DOT \bsX$  converges weakly on  
		$(\spaceD, \Dtau)$  to~$\bsY^{[h]}$ as $z\to\infty$. Suppose further that 
		\bkn{\label{wtnH} 
			\limsup_{z\to \IF} \frac{\pk{X^*_K> \ve z }}{\pk{\norm[t_0]{X }> z}}< \IF, \quad \forall \ve>0, \forall K \in \TTTK
		} 
		and for  some $c>1$ and positive  $q_t$'s such that $\sum_{t\in \TO} \max (1,p_t)q_t< \IF$
		\bkn{   \label{eq:TT}
			\lim_{\eta\downarrow 0} \limsup_{z\to\infty} 
			\pk{ \bsX_{K}^* >c\epsilon z,\exc_{K}(\epsx \DOT\bsX) \leq \eta} = 0, \quad \forall \ve>0, \forall K \in \TTTK.
		} 
	If $Y^{[h]}, h\in \TO$ is a family of tail processes,  $\exc_{K}(\cdot)$ defined in (\ref{defek}) is almost surely continuous with respect to the law of $Y^{[h]}$ for all $h\in \TO$ and further $\mcb$ \eEE{satisfies   \Cref{item:condB}}, 
		then \eEE{there exists    a  $\mcb$-boundedly finite  Borel  tail measure $\nu\in  \mS(\AA)$  such that \Cref{L:P4} holds and  $X \in \RVB, g(t)= p_{t_0}/\pk{\norm[t_0]{X }>t}$. }
		\label{prophope}	
\end{theorem}

\def\RVDB{ \RAA(a_n, \mcb_0,\tailmeasure)}

\subsection{\Cadlag\ processes} 
\label{sec:cadlag}
In this section we assume \Cref{A:2} and consider  a $\spaceD$-valued random process $X$
\cEV{defined on a  complete probability space $\Pspace$}  and not identical to 0.    Further  below $g: \RPP\mapsto \RPP$ is a 
\cEV{Lebesgue measurable} function. Alternatively to the definition in the Introduction,  as in \cite{PH2020}, $X$ is called finite dimensional regularly varying if for all 
$t_1 \ldot  t_k\in\TTT, k\ge 1$ there exists a non-trivial Borel measure $\nu_{t_1,\dots,t_k}$ on
$\borel{\R^{dk} }$ satisfying  
\bkn{  \lim_{z\to\infty}
   g(z) \pk* { \left( z^{-1} \DOT{\bsX(t_1)}{},\dots, z^{-1} \DOT {\bsX(t_k)}{}\right) \in   A }
	= \nu_{t_1,\dots,t_k} (A) < \infty
\label{eq:fidi-rv}
}
for all $A \in \borel{\R^{dk}}$ separated from $\bszero$ in $\Rset^{dk}$ such that $\nu_{t_1,\dots,t_k}(\partial A)=0$.   
Moreover the measures $\tailmeasure_{t_1 \ldot t_k}$
are $-\alpha$-homogeneous for some $\alpha>0$ and $g\in \RA$.   If we set $\nu_{t_1,\dots,t_k}(\{0\})=0$, then 
 $\nu_{t_1,\dots,t_k}$ is a tail measure on $(\R^d)^k $ with index $\alpha$.
\vEE{Since any norm  $\norm{\cdot}$ on $\R^d$ is   continuous, 1-homogeneous  and $\norm{0} =0$,  Remark \ref{lem:simp} implies that 
$$\norm[t_0]{X}= \norm{X(t_0)} \in \RA.$$ 
}
 
In view of \cite{owada:samorodnitsky:2012}[Thm 2.1]  there exists $\tailmeasure'$ on $(\R^d)^ \TTT$ equipped with the cylindrical $\sigma$-field  such that $\nu_{t_1,\dots,t_k}$ is its projection on the corresponding subspace. From the aforementioned \mb{reference} $\tailmeasure'$ is $-\alpha$-homogeneous  and moreover \Cref{L:P3} holds for all $h\in \TT$.  Denote by  $Y^{[h]}, h\in \TTT$  and  $\Theta^{[h]}, h\in \TTT$ the 
 local  tail  and local spectral tail processes of $\tailmeasure'$, respectively.  
 
Utilising  \cite{klem}[Lem 3.5], \cite{SegersEx17}[Prop 3.1, Thm 4.1,5.1] and \cite{owada:samorodnitsky:2012}[Thm 12.1] 
 the finite  RV  of $X$  implies  \mb{that}\ \eqref{klubge} holds and further \cEV{(we use the notation of \cite{klem} below for convergence in distribution)}:

\def\equallaw{\stackrel{d}{=}}
	\begin{enumerate}[(i)]
			\item \label{item:tailprocess} 
 for	all $h$ such that $p_h=\tailmeasure'( \{ f\in \spaceD: \norm[h]{f}>1\}) >0$ 	and all $t_j\in \TTT,1 \le j\le k$ 
	\bkn{		\label{eq:tail-process-time-series-tailmeasure}
			  \lim_{z\to \infty} \mathcal{L}\left( z^{-1}  \DOT \bX(t_1) \ldot  z^{-1}   \DOT \bX( t_k) \ \Big|\ \mb{\norm[h]{\bX}}>x \right) 
			= \mathcal{L}\left( \bY^{[h]}(t_1) \ldot \bY^{[h]}(t_k) \right);
}
		\item \label{item:spectraltailprocess}  
for all $h\in \TTT$  such that 	$p_h>0$ and  all $t_j\in \TTT, 1\le j\le k$  we have
\bkn{ 
	\lim_{z\to \infty}	  \mathcal{L}\left( \frac{1}{\mb{\norm[h]{\bX}} }  \DOT  \bX(t_1) \ldot \frac{1}{\mb{\norm[h]{\bX}} }   \DOT \bX(t_k)  \Big|\ \mb{\norm[h]{\bX}}>u \right) = 
			\mathcal{L}\left( \bTheta^{[h]}( t_1) \ldot \bTheta^{[h]}(t_k) \right) .
}
	\end{enumerate}  
We focus next on  $\spaceD=D(\R^l,\R^d)$ and discuss RV on $\spaceD_0=\spaceD \setminus \{0\}$ equipped with the  boundedness  $\mcb_0$  defined in \Cref{sec:prim} via \eqref{heart}. The case  $\spaceD=D(\Zset^l,\R^d)$ and some more general product spaces are already investigated in \cite{klem}.\\   

\BD   
$X$    is called RV   with limit measure $\tailmeasure\in \mM(\mcb_0)$ 
if $g(z) \pk{ z^{-1} X \in\cdot } \convvagueO \tailmeasure$ as $z\to \IF$ for some \cEV{Lebesgue measurable} $g:\RPP\mapsto \RPP$.
\ED
 A $\mcb_0$-boundedly finite measure  $\tailmeasure$ on $\borel{\spaceD_0 } $, i.e., $\tailmeasure\in \mM(\mcb_0)$ can be uniquely extended to a measure 
 $\tailmeasure^*$ on $\AA=\borelD$ by   
\begin{equation}\label{extendNU}
\tailmeasure^*(\{0 \})=0, \  \tailmeasure^*(A)=\tailmeasure(A \cap  \{ f\in \spaceD: f \not=0 \}), \ A \in \AA.
\end{equation}

\cEV{\eEE{If $\nu \in \mM(\mcb_0)$ is $-\alpha$-homogeneous, then $\nu^*$ is also $-\alpha$-homogeneous and since $\tailmeasure^*(\{0 \})=0$ is equivalent to \Cref{L:P3} we have that $\nu^*$ is a tail measure on $\AA$ with index $\alpha>0$. Given such  $\nu$}  we shall write   for notational simplicity  \mb{$\tailmeasure$ instead of $\tailmeasure^*$} and hence    
$\nu \in \mM(\mcb_0)\cap \mS(\AA)$ means $\nu \in \mM(\mcb_0)$ and $\nu^* \in  \mS(\AA)$}.

\vEE{Since $\spaceD$ is not star-shaped, for a RV $X$ with limit measure $\tailmeasure$  we cannot apply at this point \cite{hult:lindskog:2006}[Thm 3.1] to conclude that $\nu$ is $-\alpha$-homogeneous. 
 The next result shows that $\nu^*$ is even a compactly bounded tail measure. As in the  previous section we set below 
 $p_h=\tailmeasure( \{ f\in \spaceD: \norm[h]{f}>1\})$. }\\
In the rest of this section $t_0\in \TT$ is such that 
$$p_{t_0}>0.$$

\begin{lemma} If $X$ \cEV{defined on the complete non-atomic probability space $\Pspace$} is RV  with limit measure $\tailmeasure\in \mM(\mcb_0)$, then $g\in \RA$  for some $\alpha>0$ and $\tailmeasure$   extends uniquely to a tail measure on $\AA$ with index $\alpha$. Moreover,     we can take $g(t)=\cEE{p_{t_0}}/\pk{\mb{\norm[t_0]{X}}>t}$ 
	and  $\nu=\nu_Z$ with representer $Z$ satisfying  
		$\pk*{ \vk Z \not=0 }=1$. Further the local tail processes $Y^{[h]}, h\in \TO$ of $\nu$ are all defined on $\Pspace$ and both  \eqref{borsh}, \eqref{ifo} hold for all compact $K\subset \TTT=\R^l$.
\label{rusl}
	\end{lemma} 
In view of the \Cref{rusl} we can adopt the following equivalent  definition.
\BD $X$ is regularly varying with $\nu \in \mM(\mcb_0)\cap \mS(\AA)$ (abbreviated $X \in \RVD$),  if \cEV{$\norm[t_0]{X} \in \RA$} and 
$$ 
p_{t_0}\pk{ z^{-1}  \DOT X\in \cdot  }/\pk{ \mb{\norm[t_0]{X}} > z} \convvagueO \tailmeasure(\cdot).
$$ 
 \label{defRV1}
 \ED 
 Next, we  shall utilise Remark \ref{nock},\Cref{nock:1} and the explicit structure of $\mcb_0$ described in \eqref{heart}. \\ 
\cEV{In the rest of this section  assume without loss of generality that 
	$$t_0=0\in \R^l$$ 
	and this will be the assumption also for RV of $X_U$, the restriction of $X$ on $U=[\mathbf{a},\mathbf{ b}]$ a hypercube of $\R^l$ that contains $[-1,1]^l$.}\\
Denote  by $\spaceD_U=D(U, \R^d)$ \eEE{the space of \cadlag\ functions $f: U \mapsto \R^d$ with $U$ some hypercube in $\R^l$ that contains $[-1,1]^l$,}  which is also a Polish space, see \eg\ \cite{MR2802050}[Lem 2.4].  
Define the boundedness $\mcb_0(\spaceD_U)$ with respect to the zero function of $\spaceD_U$ denoted by $0_U$; $\mcb_0(\spaceD_U)$ can be characterised by  \eqref{heart} with obvious modifications.  An analogous  result to \nelem{rusl} can be formulated for $\mu \in \mathcal{M}^+(\mcb_0(\spaceD_U))$ with $\TTT=U$ and hence we can define RV of a $\spaceD_U$-valued random element similarly to that of $\spaceD$-valued random elements. Further, we extend $\mu$  uniquely to a tail measure on $\borel{\spaceD_U}$ as above.\\
Now, if $\nu \in \mathcal{M}^+(\mcb_0)$ and thus 
$\nu^* \in  \mS(\AA)$, we can define its projection with respect to $U$ denoted by $\tailmeasure_{\lvert U}^*$ as 
the tail measure on $\borel{\spaceD_U}$ determined uniquely by  $Y^{[h]}_{U}, h\in U$, where $Y^{[h]}$'s are the local tail processes of $\nu^*$, since their restriction on $U$ denoted by $Y^{[h]}_{U}$ yields a family of tail processes on $\spaceD_U$. Write then $\tailmeasure_{\lvert U}$ for the restriction of $\tailmeasure $ on $\borel{\spaceD_U \setminus \{ 0_U\}}$.
 
Let $\proj_U: \spaceD \mapsto \spaceD_U,$ with $\proj_U(f)= f_U$ be the restriction of $f\in \spaceD$ on $U$.    In view of \netheo{sam}, \Cref{tev:8} we can find $\mathbf{a}_n, \mathbf{b}_n$ such that  $\nu^*(Disc(\proj_{U_n}))=0$ and $[-n,n]^l \subset [\mathbf{a}_n, \mathbf{b}_n]=:U_n$ for each given positive integer $n$. 
 \begin{theorem} 
 	If  $U_n,n\in \Nset$ is as above and $X \in \RVD$, then  
	$$X_{U_n} \in \mathscr{R}_\alpha(\mcb_0(\spaceD_{ U_n}),\nu^{(n)}), \ \ \nu^{(n)}=\nu_{\lvert U_n}, \ \  \forall n\in \Nset.$$
	 Conversely, if  
	$X_{U_n} \in \mathscr{R}_\alpha(\mcb_0(\spaceD_{ U_n}),\nu^{(n)})$ for all $ n\in \Nset$, then  
\begin{equation}\label{warten}
	X \in \RVD, \quad \nu_{\lvert U_n}=\nu^{(n)}, \ \forall  n\in \Nset, \ \quad  \nu \in \mM(\mcb_0)\cap \mS(\AA).
\end{equation}	
\label{thmA}
\end{theorem}

 \BRM \begin{enumerate}[(i)]
 	\item \label{enu:i} \cEV{Both \nelem{rusl} and \netheo{thmA} hold also  for   $\spaceD=D(\Zset^l, \R^d)$;} 
 	\item \label{enu:2} If there is a $\spaceD$-valued random element $Z$  such that 
 	\bqn{ \pk*{ \vk Z \not=0 }=1, \quad \E{ \norm[t]{\vk Z }^\alpha }  \in [0,\IF), \forall t\in \TTT,  \E{ \norm[0]{\vk Z }^\alpha } >0,  \ \E*{ \sup_{t\in K \cap \TO} \norm[t]{\vk Z}^\alpha}< \IF, \ \forall K \subset \R^l, 
 		\label{Za}
 	} 
 with $K\subset \TT$ compact  	and 	 $\nu^{(n)}=\nu_{Z_{n}}, \forall n\in \Nset$, where 
 \bqn{\label{zn} 
	Z_n(t)=  c_n^{1/\alpha}  Z(t) \Bigl \lvert \sup_{t\in U_n} \norm[t] { Z}>0,\quad t\in U_n,  \   c_n= \pk*{ \sup_{t\in U_n} \norm[t] { Z}>0 }>0,
 }
  then it follows from the proof of  \netheo{thmA} that \eqref{warten} holds with $\nu=\nu_Z$. Conversely, if  $\tailmeasure$ has repersenter $Z$, then 
  $\nu^{(n)}=\nu_{Z_{n}}, \forall n\in \Nset$ holds.
 \end{enumerate}
\label{enu}
 \ERM

Consider next an $\R^d$-valued   max-stable random process $X(t),t\in \TTT $ given  via its de Haan representation  (e.g., \cite{RDTM:2008,Roy})
\bqn{\label{eq1}
	X(t)=  \max_{i\ge 1} \Gamma_i^{-1/\alpha}   \DOT Z^{(i)}(t), \quad t\in \TTT.
}
Here $\Gamma_i= \sum_{k=1}^i \mathcal{E}_k$, where  $\mathcal{E}_k, k\ge 1$ are independent unit exponential rv's being independent of $ Z^{(i)}$'s which are independent copies of $Z(t),t\in \R^l$ with almost surely sample paths in $\spaceD$ satisfying \eqref{Za}.  In view of \cite{dom2016}  $X$ is max-stable. Commonly, $Z$ is referred to as  a spectral process of  $X$. Let $\nu_Z$ be the tail measure corresponding to $Z$, which is compactly-bounded by \eqref{Za}. The law of $X$ is uniquely determined by $\nu_Z$ or the local tail processes of $\nu_Z$, see \cite{klem}. \cEG{Moreover, as shown in \cite{MR3834849,MR4227144} $X$ is stationary if and only if  
	 (see also \cite{PH2020}[Thm 2.3] for the case $l=1$) 
	\bqn{ \label{tcfN} 
			\E{ \norm{Z(h)}^\alpha  F({\vk Z})} = \E{\norm{Z(0)}^\alpha F( B^h   {\vk Z})}, \quad \forall F\in \Hh_{0} , \forall h\in \TT 
		} 
	holds. It follows that \eqref{tcfN} is also equivalent with  $\nu_Z$ is shift-invariant, see also \cite{PH2020}[Thm 2.3] discussing $l=1$. 
}

\begin{corollary}\label{shikon}
If $X$ is given by \eqref{eq1} with $Z$ satisfying \eqref{Za},  then $X\in \RA(\mcb_0,\nu_Z)$. 
\end{corollary}

\begin{example}[Brown-Resnick max-stable processes] Let  
	$$\vk Z(t)= (   e^{  W_1(t)} \ldot 
	\  e^{  W_d(t)} ), \quad  W_i(t)= V_i(t) - \alpha Var(V_i(t))/2, \quad 1 \le i \le d, t\in \TT=\R^l, $$
	with $\alpha >0$, 
	$ (V_1(t) \ldot V_d(t)), t\in \TTT$ a centered $\R^d$-valued Gaussian process with almost surely continuous sample paths such that $V_i(0)=0, i\le d$ almost surely.	\cEE{In the light of \cite{MR3024389}[Cor.\ 6.1]}, Eq.\ \eqref{Za} holds, 
	and thus by Remark \ref{mbab},\Cref{mbab:0}    $\tailmeasure_Z$ is $\mcb_0$-boundedly finite on $\spaceD=C(\R^l,\R^d)$, 
	{the space of continuous functions $f:\R^l \mapsto  \R^d$ equipped with a metric that turns it into a Polish space}.  Consider the    max-stable process  $X$ with spectral process $Z$.  Corollary \ref{shikon} implies 
	$X \in \RAA(\mcb_0, \tailmeasure_Z)$. See also \cite{dombry} for the case where $\spaceD=C(K, \R)$ is considered with $K$ a compact set of $\R^l$.
	\label{ex:BR}
\end{example}

We focus next on $\spaceD=D(\R, \R^d)$ and utilise  \netheo{prophope} since  $\mcb_0$ is determined by the family of  maps $\famHMAP$.  See Remark  \ref{remskr}  and   \Cref{sec:spaceD} for \mb{the} definition of $w,w'$ and $w''$ that appear below.
\begin{theorem}
  \label{Th:RV}
  Let   $\bsX$  be   defined on a complete non-atomic probability space $\Pspace$ and let $t_0=0$.  
 The following  statements are equivalent:
  \begin{enumerate}[(i)]
  \item \label{item:RVD}   $\bsX \in \RVD$ with $\mb{\norm[t_0]{X}}\in \RA$;  
  \item \label{item:tightness}  \mb{Eq.} \eqref{eq:fidi-rv} holds  for all $t_1,\dots,t_k\in  \tE{T_0} , k\ge 1$  for some $T_0$ such that $\TTT\setminus T_0$ is countable
  and  
\bkn{      \label{eq:tightness}
      \lim_{\tE{\eta}\downarrow 0} \limsup_{z\to\infty} \frac{ \pk{  w'(\bsX, K,\tE{\eta}) > \ve  z  }}{\pk{ \mb{\norm[t_0]{\bsX}>z}}} = 0, 
\quad      \forall \epsilon>0, \forall K \in \TTTK.
}
\item \label{item:conditional} 
\mb{Eq.} \eqref{klubge} holds and  
for all $h\in  T_0$ for some $T_0$ such that $\TTT\setminus T_0$ is countable
\bkn{      \label{eq:tightness2A}
	\lim_{\tE{\eta}\downarrow 0} \limsup_{z\to\infty} \frac{ \pk{  w'(\bsX, K,\tE{\eta}) > \ve  z , \mb{\norm[h]{X}}\le   z }}{\pk{ \mb{\norm[t_0]{\bsX}}>z}} = 0, 
	\quad      \forall \epsilon>0, \forall K \in \TTTK.
}
Further if  $p_h>0$, then  		conditionally on $\mb{\norm[h]{\bsX}}>z$, $z^{-1}  \DOT \bsX$  converges weakly on   
$(\spaceD, \Dtau)$  to~$\bsY^{[h]}$ as $z\to\infty$, where   $Y^{[h]}$'s are $\spaceD$-valued random processes 
defined on $\Pspace$ being further the local tail processes of a tail measure $\tailmeasure$ on $\AA$  with index $\alpha>0$;
\item \label{item:proj}  Let  $s_k<t_k ,k\in \Nset$  be given constants  \mb{satisfying}  $-\limit{k}s_k=\limit{k}t_k=\IF$ and set $K_k=[s_k,t_k]$. There exists 
$\mcb_0(\spaceD_{K_k})$-boundedly finite Borel measures  $\tailmeasure_k$ with $\tailmeasure_k^*$ its corresponding tail measure  on $\borel{\spaceD_{K_k}}$
 with   index $\alpha>0$. Suppose that $\tailmeasure_k(\{f\in \spaceD: f(t)\not= f(t-)\})=0$ for 
$t\in \{ s_k,t_k\}$  and 
   $X_{K_k} \in \mathscr{R}_\alpha (\mcb_0(\spaceD_{K_k}), \nu_{k})$ for all $k\in \Nset$ with $\mb{\norm[t_0]{X}}\in \RA$.
  \end{enumerate}
\end{theorem}

\BRM  
  \begin{enumerate} [(i)] 
	\item \label{en:thm:2}  For   $l>1$ and    $\spaceD=D(\R^l, \R^d)$, if \Cref{item:RVD} holds, 
	then  the weak convergence in \eqref{posht} below and \netheo{sam}, \Cref{tev:4} imply  \netheo{Th:RV},\Cref{eq:tightness} where   
	\eqref{eq:tightness} is substituted by 
	\bqn{ \label{ture0}
	\lim_{\tE{\eta}\downarrow 0} \limsup_{z\to\infty} \frac{ \pk{  w'(\bsX, K,\tE{\eta}) > \ve  z , 
	\sup_{t \in [-k,k]^l \cap \TO} \norm[t]{X} >  z/k  }}{\pk{ \norm[t_0]{\bsX} >z }}& =& 0, 
	\quad       \forall K \in \TTTK,\\
\label{ture}
		\limit{m} \limsup_{z\to \IF} \frac{ \pk{ \sup_{t \in [-k,k]^l \cap \TO} \norm[t]{X} > mz }} 
		{\pk{{\norm[t_0]{\bsX}}>z}}& =&0,
}
with  
$k \in \Nset,\ve>0$  arbitrary. Conversely, 
\netheo{Th:RV},\Cref{eq:tightness} with the above modification implies \netheo{Th:RV},\Cref{item:RVD} and similarly \netheo{Th:RV},\Cref{item:conditional} therein can be modified to yield the equivalence with \netheo{Th:RV},\Cref{item:RVD}. 
	\item \label{en:thm:1} If 
instead of $D(\R^l, \R^d)$ we consider $C(\R^l, \R^d)$, then  \netheo{Th:RV}  holds with $w$ instead of $w'$.  This follows since in this case we can substitute $w''$ by $w$ in \eqref{ture0}. 
	\item \label{en:thm:3}
By  \cite{billingsley:1999}[Eq.~(12.28)] and \eqref{eq:tightness}  when $l=1$
 \bkn{      \label{eq:tightness2C}
 	\lim_{\tE{\eta}\downarrow 0} \limsup_{z\to\infty} \frac{ \pk{  w''(\bsX, K,\tE{\eta}) > \ve  z  }}{\pk{ \mb{\norm[t_0]{\bsX}}>z}} = 0, 
 	\quad      \forall \epsilon>0, \forall K \in \TTTK.
 } 
It follows using \cite{billingsley:1999}[Eq.~(12.32)] and 
\cite{hult:lindskog:2005}[Thm 10] that     \eqref{eq:tightness} can be substituted by 
\mb{\eqref{eq:tightness2C}}.
  \end{enumerate}  
\label{remskr}
\ERM

\begin{example}[Random scaling]   
Under \Cref{A:2}	let $ \tailmeasure_Z$ be a compactly-bounded tail measure on $\AA$.  
Let  $R$ be  an $\alpha$-Pareto rv  independent of $Z$  and  set   	$X(t)= R Z (t) , t\in \TTT$.  Utilising \cite{RDTM:2008}[Lem 2.3 (2)], since $\E{ \sup_{t \in K} \norm[t]{Z}^\alpha}\in (0,\IF)$ for all  compact  $K\in \R^l$,  it follows from Remark \ref{enu},\Cref{enu:2} and Remark \ref{remskr},\Cref{en:thm:2}   
that  $ X \in \RVD $. We note   that this example is discussed in \cite{RDTM:2008} for compact $\TTT$.  
\label{ex:RS}
\end{example}

\begin{example}[Scaled \& shifted processes] 
	\label{ex:scaled} 
 Suppose that  $X \in \RVD$ is   a $\spaceD$-valued random element, $Y^{[h]}, h\in \TTT$ are  the local processes of a \eEE{$\TTTK$-bounded} Borel tail measure  $\tailmeasure$ on $\spaceD$. Let   $X_{f,\sigma}(t)= \sigma(t)X(t) + f(t), t\in \TTT, $ with  $f,\sigma\in \spaceD$ such that 
	 $\sigma \in \spaceD$ is continuous and   $\sigma(t)\not=0$ for all $t\in \TTT$. 
 Note that if $\TTT=\R^l$, then  {$\lim_{\delta \to 0} w'(f,K,\delta)=0,\forall K\in \TTTK$,   see 	\netheo{sam},\Cref{tev:0}}.
	  Using Remark \ref{remskr},\Cref{en:thm:2} 
	we have  $X\in \RAA(\mcb_0 , \tailmeasure_{\sigma} ),$ where the tail measure
	\begin{equation}\label{umdenken}
	\tailmeasure_{\sigma}(A)= \tailmeasure(\{f\in \spaceD: \sigma f \in A \}), \ \ A \in \AA
\end{equation}
  has local tail processes given by  
$Y_{\sigma}^{[h]}(t)=  Y ^{[h]}(t)\sigma(t)/\sigma(h)$ for all $ h,t\in \TTT.$ 
 \end{example}

 \section{Discussions}
\label{sec:disc}
We shall consider first another common definition of RV, {in terms of sequences, see \cite{MR3331244} and also \cite{bingham2020sequential} for a recent full account.} The second part of our discussions is dedicated to  RV  under transformations and then we conclude with a short section on stationary \cadlag\ processes.

\subsection{An \mb{alternative} definition of  RV}
\label{sec:disc2}

Suppose   that \Cref{A:5}  
\mb{holds} and let  in the following  $a_n>0, n\in \Nset$ be a  non-decreasing sequence of constants such that  
$$\limit{n} a_{[nt]}/a_n= t^\alpha, \ \ \forall t>0,$$
where $[x]$ denotes the integer part of $x$. For such constants we write $a_n \in \RA$.
Another   common and  less restrictive definition of  RV  (see e.g., \cite{kulik:soulier:2020}[Thm B.2.1]) is the following: 

\BD  $\mu \in \MBB$ is regularly varying if for    $a_n \in \RNA$   
$$\mu_n(A)=n\mu( a_n \DOT  A) ,  \quad A \in \borelD$$ 
converges $\mcb$-vaguely    to  some   $\tailmeasure\in \MBB$  as $n\to \IF$, 
abbreviate this as $\mu \in  \RVC $.
\label{secretarin}
\ED   
If $\mu \in \RVB$  and  $g\in \RA$,   \nelem{lem:kulikB22} yields  
$\mu \in \mathcal{R}_\alpha(g_*, \mcb, \tailmeasure)$ for any  Lebesgue measurable 
$g_*: \RPP\mapsto \RPP$ such that $\limit{z} g(z)/g_*(z)=1$. Since $g\in \RA$, we can choose   $g_*\in \RA$ asymptotically non-decreasing. Taking then  $a_n= g^{-1}_*(n), n\ge 1$ with $g^{-1}$ an  asymptotic inverse of $g$,  it follows that  
$$\mu \in \RVB \implies \mu \in \RVC.$$

The inverse implication above (and thus the equivalence of both definitions of  RV)  can be  proven under \cite{kulik:soulier:2020}[(M1)-(M3),(B1)-(B3), p.\ 521/522], see  \cite{kulik:soulier:2020}[Thm B.2.2]. 

\cEV{The Definition \ref{secretarin}  can be naturally extended to $\spaceD$-valued random processes $X$, which is   abbreviate as 
	$$X \in \widetilde{\mathcal{R}_\alpha}(a_n, \mcb,\nu).$$ 
	Both definitions of RV for \cadlag\ processes are equivalent as we show next.    
\BEL If $\nu, X$  are as in \netheo{thmA}, then 	 	$X \in \RAA(\mcb_0,\nu)$ is equivalent to  $X \in \widetilde{\mathcal{R}_\alpha}(a_n, \mcb_0,\nu)$, where \cEG{$a_n$ is such that } $n \pk{ \mb{\norm[t_0]{X}} > a_n} =p_{t_0}$ for all large $  n\in \Nset$.  
\label{chans}
\EEL 
}
 
\begin{example}
	\label{examp:RD}
	We consider  the setup of  
	\cite{kulik:soulier:2020}[Prop 2.1.13] assuming  \Cref{A:2}.  
	Let   $X \in \widetilde{\mathcal{R}_\alpha}(a_n, \mcb_0,\nu) $ with $\mb{\norm[t_0]{X}}\in \RA$ for some $t_0\in \TT$  
	and 	let 
	$\Gamma:\TT\to \R^k$ be a random map  independent of $X$ defined on $\Pspace$ and let $\norm{\cdot}$ be some norm on $\R^d$. Suppose that  almost surely $\Gamma(cu)= c^\gamma \Gamma(u), \forall c>0, u \in \R^d$ for some $\gamma>0$. Assume    that \mb{$\Gamma$} is  almost surely continuous satisfying   \eqref{ture0}. If further 
	\eqref{ture} holds with $X$ substituted by $\Gamma \circ  \bsX$ 	and for some $\ve>0$
	\bkn{  \E*{ \sup_{\norm{u}\le 1} [\Gamma(u)]^{\alpha/\gamma+ \ve}}< \IF, \label{inchina}
	}
	where $\alpha>0$ is the  index of $\tailmeasure$, then 
	$\Gamma(X) \in
	\widetilde{\mathcal{R}_\alpha}( a_n^{\gamma},\mcb_0,\E{\tailmeasure \circ \mb{\Gamma}^{-1}})$, provided that $\E{\tailmeasure \circ \mb{\Gamma}^{-1}}$ is \cEV{non-trival}. \cEG{A particular instance of interest is $\Gamma(t) = A X(t), t\in \TT$ with $A$ a $k\times l$ real matrix satisfying \cite{kulik:soulier:2020}[Eq.\ (2.1.14)].}
\end{example}

\subsection{Transformations}
  We shall focus in this section on 
 $\spaceD=D(\R, \R^d)$. 
The next lemma is a restatement of \cite{MR3238572}[Lem 3.2] for our setup. 

\BEL If $X \in \widetilde{\mathcal{R}_\alpha}(a_n, \mcb_0,\nu) $ and $\sigma$ is a $\spaceD$-valued random process independent of $X$, then  
 $\sigma X \in \widetilde{\mathcal{R}_\alpha}(a_n, \mcb_0 , \E{ \tailmeasure_\sigma })$, \cEE{with $\tailmeasure_\sigma$ defined in \eqref{umdenken}}, provided that 
\bkn{ \label{sigmak} 
	\E*{ (\sup_{t \in K} \mb{\norm[t]{\sigma}})^{\alpha+ \ve }}< \IF, \quad \forall K \in \TTTK
}
for some $T_0$ such that $\TTT\setminus T_0$ is countable  \mb{and}
\bkn{ \label{sigmak2}
	\pk{ \cEE{{\sigma(t)}} \not=0}>0, \quad \forall t\in T_0.
}
\label{40bura}
\EEL  
\tE{In view of \netheo{Th:RV},\Cref{item:proj} }      \nelem{40bura}  can be extended   considering  $X_i,i\le m$ independent copies of $X$ and $\sigma_i, i\le m$ $\spaceD$-valued random processes. Then  \cite{MR3238572}[Lem 3.3] can be restated by imposing \eqref{sigmak} \mb{and \eqref{sigmak2}} on all $\sigma_i$'s. \\
{\BEL
If $X \in \widetilde{\mathcal{R}_\alpha}(a_n, \mcb_0,\nu) $, then    \cite{MR3238572}[Thm 3.3] holds also if  
in  the assumptions therein  $\lvert \sigma_j \rvert_{\IF}$ is substituted by 
 $\sup_{t\in K} \norm{\sigma_j(t)}$, for all compacts $K\in \TTTK$.
 \label{lem:hermine} 
\EEL 
}

\label{sec:LDL}

 \subsection{Stationary \cadlag\ processes}
\label{sec:stationary}
Under the settings of \Cref{sec:cadlag} assume further that $\bsX$   is  stationary. Hence 
$\mb{\norm[t_0]{X}}\in \RA$ for some  $t_0\in \TTT$  implies  $\mb{\norm[t]{X}}\in \RA$ at all $t\in \TTT$ and thus  
$$p_h=p_0\in (0,\IF), \ \ \forall h\in \TTT.$$
It follows easily using \eqref{eq:tail-process-time-series-tailmeasure} or directly by \netheo{Th:RV} and \cite{WS}[Thm 3.2] that if $\bsX \in \RVD$   holds, then the local tail processes $Y^{[h]}, h\in \TTT$ are given by 
\eqref{localSatio} which implies that the corresponding tail measure $\nu$ is shift-invariant. Moreover,  the converse holds i.e., if $\tailmeasure$ is shift-invariant, then  \eqref{localSatio} holds,  see  also \cite{klem}.  \\
With this additional knowledge on $\tailmeasure$, the counterpart of \netheo{Th:RV} for stationary $X$ can be easily reformulated. 
 
\begin{example}[Stationary Brown-Resnick max-stable processes] Let $Z,X,W$ be as in Example \ref{ex:BR}, $d=1$  and suppose that $X$ is stationary, which follows if $ Var(W_1(t)- W_1(s)), s,t\in \TT$ depends only on $t-s$ for all $s,t \in \TT$. 
 	 We have that $X \in \mathscr{R}_\alpha (\mcb_0, \tailmeasure_Z)$ with $\nu_Z$ having representer $Z$ and being shift-invariant.
Note that  since $\abs{Z(0)}=1$ almost surely, then $\Theta=Z$ is the local spectral process of $\nu_Z$ at $0$ and thus 
$$ 
Y^{[0]}(t)=e^{ \alpha \eta+  W_1(t)}, \ t\in \TTT, 
$$
with $\eta$ a unit exponential rv independent of $W_1$. Hence  \eqref{eqR1} reads for   $\alpha=1, x>0$  
 \bqn{ 
	\label{eqR12} \E{ \Gamma( x e^{ \eta + B^h W_1 } ) \ind{   W_1(-h)+ \eta  >- \ln x}} = 
	x \E{  \Gamma(   e^{\eta+ W_1}) \ind{ W_1(h)+ \eta  >\ln x }   }, \ \forall h\in \TT,  \forall\Gamma \in \Hh.
}  
\end{example}

\section{Applications \& open questions}
\cEG{\label{openq}
We first mention four  applications 
 considering  $X \in \RVBA$ as in \Cref{sec:rvRelements} assuming \Cref{A:5}. 
\begin{enumerate}[{AP}1)]
\item \label{app2} A  well-known application of RV is the derivation of 
the   tail behaviour of $H(X)$, for a given functional of interest  $H$. When $X$ has \cadlag\ sample paths, a canonical choice is 
$H(f)=H_K(K)=\sup_{t \in K}\norm{f},f \in \spaceD,$ with  $K$ a compact set in $\R^l$, or 
$H(f)=:H_A^*(f) =\int_A f(t) \lambda(dt), f \in \spaceD,$ with $A$ a bounded Borel set in $\R^l$ of positive Lebesgue measure. Conditions on $H$ for tractable   tail behaviour of $H(X)$ are presented  in Remark \ref{lem:simp} below;  
\item \label{app1} As already discussed in several contributions, see e.g., \cite{dombry:hal-02168872}[Prop 2.3], RV   implies the convergence of the collection of probability measures  $\mu_z(A)=\pk{z^{-1}\,\DOT\, X \in  A \mid z^{-1}\DOT X \in  B}$ for all $A\in \AA$, as $z\to\infty$. Assuming, additionally, that the Borel set $B$ belongs to $\mcb$ and is $\tailmeasure$-continuous (i.e., $\tailmeasure(\partial B)=0$) with $\tailmeasure(B)>0$, we obtain 
\bqn{ \mu_z \convvague \mu, \ z\to \IF,
}
where $\mu(\cdot )= \tailmeasure(\cdot \cap B) /\tailmeasure(   B)$.
This application is useful for the formulation of conditional limit results, as already shown in the aforementioned  contribution;
\item \label{app3} An interesting application considered for the discrete setup is developed recently in {\cite{AAP1579}} for the product of RV random matrices. The results therein can be extended to the product of random matrix functions, making use of \Cref{thmA}, \Cref{Th:RV}, and ideas given in the aforementioned contribution. Moreover, extensions to more general homogeneous functionals can be also obtained using further Remark \ref{lem:simp};
\item \label{ap4} 
One advantage of introducing the RV with respect to some boundedness is that it includes also the concept of hidden RV discussed for instance in \cite{MR2271424,MR3271332,dombry:hal-02168872}.  Indeed, in our settings hidden RV corresponds to the choice of   $\mcb_F$  with $F \subset  \spaceD$ closed being the boundedness on $\spaceD_F=\spaceD \setminus \{F\}$ defined in \Cref{sec:boundV}. Since by definition   we need $\spaceD_F=\spaceD\setminus F$ to be a measurable cone, we shall further assume that $F$ is a cone, \ie\ 
$$z\DOT F \subset F, \  \  \forall z\in \RP.$$
Define similarly $\mcb_{F'}'$ with respect to $\spaceD'_{F'}$ with $ F'\subset \spaceD'$ a closed cone. The next result is useful when considering maps of hidden regularly varying processes.
\begin{lemma}
	Let   $H: \spaceD\mapsto \spaceD' $  be  $\borelD/\borel{\spaceD'}$ measurable with $H(F)=F'$. Let  further  $\nu_z, z>0$ be $\mcb_F$-boundedly finite measures on $\borel{\spaceD_F}$ and let  $\nu$ be a $\mcb_{F'}'$  boundedly-finite measure on $\borel{\spaceD'_{F'}}$.  
	If one of the following conditions
	\begin{enumerate}[(i)]
		\item $H$ is uniformly continuous;
		\item $\spaceD$ or $F$ are compact and $H$ is continuous;
		\item $\nu(Disc(H))=0$, $H$ is a continuous and one-to-one if restricted on $F,$ which has finite number of elements  
	\end{enumerate}  
	is satisfied and $\nu_z \stackrel{v,\mcb_F}{\longrightarrow}\nu$,  then $\nu_z  \circ H^{-1}\stackrel{v,\mcb'_{F'}}{\longrightarrow}  \nu \circ H^{-1}$.
	\label{homogeneousTransf}
\end{lemma} 
\BRM  Under the conditions of \nelem{homogeneousTransf}, if 	  
$H(c\DOT f)=c \DOT H(f), \forall c> 0, f\in \spaceD$ and $\tailmeasure \circ H^{-1}$ is non-trivial, then 
$\mu \in \widetilde{\mathcal{R}_\alpha}( a_n,\mcb_F,\tailmeasure)$  implies $\mu \circ H^{-1} \in \widetilde{\mathcal{R}_\alpha}( a_n , \mcb_{F'}', \tailmeasure \circ H^{-1})$ and thus under the  settings of \cite{klem} we retrieve 
the claim of   Lemma 3.2 therein.  
\label{lem:simp}
\ERM
\end{enumerate} 
The applications and findings of \cite{PH2020} for $\TTT=\R^l, l=1$ can be  extended to our general case of stationary $X$ for all integer $l>1$, by alluding to the methodology developed therein, together with \Cref{thmA}. We do not presently repeat all calculations, but rather mention a few details and some new results on the tail behaviour of supremum of \cadlag\ processes. \\ 
The rest of this section considers  random processes $X(t), t\in \TT$ with $\TT=\Zset^l$ or $\TTT=\R^l$. In the latter case we assume that 
$X$ has \cadlag\ sample paths.  Further, $\tailmeasure_Z$ is a tail measure with representer $Z$, which has  almost surely \cadlag\ sample paths if  $\TTT=\R^l$. Suppose next that  $\E{\norm{Z(0)}^\alpha}=1$, where  $\norm{\cdot}$ is a  norm on $\R^d$.
\subsection{Stationary  case}
We now consider the case when $X$ is stationary and $X \in \widetilde{\mathcal{R}_\alpha}(a_n, \mcb_0,\nu_Z)$. Hence   
$\nu_Z$ is shift-invariant and therefore uniquely determined by  $Y=Y^{[0]}$. Note in passing  that  $\E{\norm{Z(0)}^\alpha}=1$ implies that $p_h=
\E{\norm{Z(h)}^\alpha}=1$, for all $h \in \TT$. \\
The determination of the tail behaviour of $H_K(X)$, with $H_K$ the supremum functional  in \cref{app2},  is a classical interesting problem of probability theory. As already demonstrated  in \cite{PH2020},  our results  can be applied to consider both the case $K$ does not depend on $n$ and  $K=K_n=[0,n]^l$,  when  $n$ tends to infinity. We first state a general upper bound for the growth of the supremum 
$$M_n=  \sup_{t\in [0,n]^l \cap \TT} \norm{X(t)},$$
 assuming for simplicity that  $\norm{\vk X(0)}$ is a unit $\alpha$-\mbox{Fr\'{e}chet} rv with df $e^{-x^{-\alpha}},\: x>0$.
\begin{proposition} \label{genUB} If $X \in \widetilde{\mathcal{R}_\alpha}(a_n, \mcb_0,\nu_Z)$ where $a_n=n^{1/\alpha}$ and  $\norm{\vk X(0)}$ is a unit $\alpha$-\mbox{Fr\'{e}chet} rv,   then    for all $x>0$
	\bqn{ \limsup_{n \to \IF} \pk*{M_n> a_n  ^l x} \le \theta_Y 
		x^{- \alpha},\quad \theta_Y=\E*{\frac{1}{\int_{\TT} \mathbb{I}( \norm{ Y(t)}>1 ) \lambda(dt)} } \in (0,\IF),
		\label{raki}
	}
provided that for $\TT=\R^l$ we have  
	$I_{k,n}=\E*{ 1/\int_{s \in [-k,k]^l } \ind{ \norm{\vk X(s)}>a_n \lvert X(0)>a_n  } \lambda(ds)  }$ is bounded 	for all $n$ large and some $k>0$.
\end{proposition} 
Note  that  \eqref{raki} is shown in \cite{kulik:soulier:2020}[Lem 7.5.4] for 
 $\TT=\mathbb{Z}^l, l=1$ and $I_{k,n}\le 1$ for all positive integers  $l,k,n$ if $\TT=\mathbb{Z}^l$. \\ 
 In the particular case that 
$$ \ve(Y)=\int_{\TTT} \mathbb{I}( \norm{ Y(t)}>1 ) \lambda(dt)= \IF$$ 
almost surely, then $\theta_Y=0$. Hence \Cref{genUB} implies  the following convergence in probability
\bqn{ \label{matura}
 a_n  ^{-l}M_n  \convprob 0, \quad n\to \IF.
}
In order to establish weak convergence of  $ a_n  ^{-l}M_n $ to some \mbox{Fr\'{e}chet} rv as $n\to \IF$,  we have to guarantee the positivity of  $\theta_Y$. In both the discrete setup 
(cf., \cite{BojanS,basrak:planinic:2019,HBernulli}) and the continuous case with $l=1$ dealt with in \cite{PH2020}, it is known that $\theta_Y>0$ follows from the anticlustering condition 
of \cite{Davis}, which we now present for both cases $A=\mathbb{Z}^l$ and $A=\R^l$.\\
 We say that $f$ is a scaling function, if $f:(0,\IF) \to (0,\IF)$ is non-decreasing and unbounded,  and set $\norm[\IF]{x}=\max_{1\le i \le l} \abs{x_i}, x\in \R^l$.
\begin{Condition}[$C(A)$]
	\label{anticlustering}
	 There exist scaling functions  $a$ and $r$  such that   
	\bkn{
		\label{eq:anticlustering}
		\lim_{t\to\infty} \limsup_{y\to\infty} \pr\left\{\sup_{t\leq \norm[\IF] {s}  \leq r(y), s\in A} \norm{\bsX(s)} > a(y) x \:\Big|\:  \norm{X(0)}>a(y)\right\} = 0, \quad \forall x>0. 
	} 
 \end{Condition}
  From the anticlustering condition we may derive important properties of $Y$ and $Z$, in particular that $\theta_Y>0$. 
\BEL
\label{lem:gio2} Under the assumptions of \Cref{genUB}, if $\TTT=\R^l$ and  Condition 6.1 (C($\mathbb{Z}^l)$) holds, then 
	\bqn{\label{gio}
\pk*{ \int_{\TTT} \norm{ Y(t)}^\alpha \lambda(dt) \in (0,  \IF)}=1. 
 }
\EEL 
As in \cite{PH2020}, we say that the shift-invariant tail measure $\tailmeasure_Z$ is dissipative if \eqref{gio} holds almost surely. Along the same lines of the aforementioned paper, it follows that $\tailmeasure_Z$ is dissipative if and only if $\ve(Y)$ and      
$ \int_{\R^l} \norm{ Z(t)}^\alpha \lambda(dt) 
$
 are almost surely positive and finite, implying in particular that  $\theta_Y>0$.\\
Moreover, if  $\tailmeasure_Z$  is dissipative, the PPP $N$ defined in Remark \eqref{nock},\cref{nock:2}, has the following representation 
$$  N(\cdot)= \sum_{i=1}^\IF \delta_{P_i B^{\tau_i} Q^{(i)}}(\cdot),$$
where $\sum_{i=1}^\IF \delta_{P_i, \tau_i, Q^{(i)}}(\cdot)$ is a PPP on $(0,\IF) \times \R^l \times  \spaceD$, with mean measure $\theta_Y \lambda(\cdot) \odot  v_\alpha(\cdot) \odot \mathbb{P}_Q(\cdot)$ and $Q^{(i)}$'s are independent copies of $Q$ with law 
$$ \mathbb{P}_Q(\cdot)=\theta_Y^{-1} \E{ \delta_{Y/\sup_{t\in \TTT } \norm{Y(t) } }(\cdot)/\ve(Y) }.$$
The  dissipative representation of $N$   is key  to the so-called $m$-dependent approximation (see \cite{klem} for the definition). 
Specifically, if  $\tailmeasure_Z$ is dissipative,  the max-stable stationary process  $X$ defined in \eqref{eq1} (recall $Z$ has non-negative components) has the dissipative representation 
$X(t)= \max_{i\ge 1} P_i  Q^{(i)}(t-\tau_i), t\in \TTT$, which has an $m$-approximation given by  
$$X^{(m)}(t)= \max_{i\ge 1} P_i {Q^{(i)}(t-\tau_i)}\ind{ \norm{t-\tau_i }\le m},\quad  t\in \TTT,m>0.$$
Similarly, an $m$-approximation can be derived for the $\alpha$-stable stationary $X$ with $\alpha \in (0,2)$ defined by substituting $\max$ with  $\sum$ in \eqref{eq1} and in the above dissipative representation. 
In order to avoid centering, when $\alpha \in [1,2)$, as  in \cite{PH2020},   $Z$ is further assumed to be symmetric. In both cases, $(X, X^{(m)})$ is stationary. The next results extends \cite{PH2020}[Thm 4.1, Cor 4.3, Thm 4.5, Cor 4.6] (note that $T \pk{\cdot}$ should be $\pk{\cdot}$ therein) to $l\ge 1$.  Related results are derived also in \cite{Genna04, Genna04c,MR2384479,Roy2}.
\begin{theorem}
	If $X$  as above is max-stable or $\alpha$-stable, then  $X \in \widetilde{\mathcal{R}_\alpha}(a_n, \mcb_0,\nu_Z)$ with $a_n=n^{1/\alpha}$. 
Moreover, we have 
\bqn{  a_n  ^{-l}M_n  \stackrel{d}{\to}   \eta_X^{1/\alpha} V, \quad n\to \IF,
	\label{jugendlich}
}
with $V$ an $\alpha$-\mbox{Fr\'{e}chet} rv and $  \eta_X=\theta_Y<\IF$. 
\label{thm:SP} 
\end{theorem}	
\BRM In the literature  $\theta_Y$ is commonly referred to as  the candidate extremal index, which under the assumptions of \Cref{thm:SP} is equal to the extremal index $\eta_X$ of $X$. As shown  in \cite{wangExt}, it is possible to have  $\eta_X < \theta_Y$. Note in passing that when \eqref{jugendlich} holds and $X$  as in \Cref{genUB} is stationary and regularly varying with tail process $Y$, then \eqref{raki} implies  
\bqn{  \eta_X \le \theta_Y,
}
which for $\TTT=\Zset$ is shown in \cite{kulik:soulier:2020}[Lem 7.5.4].
\ERM
\subsection{Non-stationary  case}
  The non-stationary case is significantly less tractable, when compared to the stationary one. Yet, there are a few exceptions,  for instance  $X_{f,\sigma}$ defined  in \Cref{ex:scaled}, with $X$ stationary and regularly varying. Under some growth restrictions on $f$ and $\sigma$, \Cref{genUB} and \Cref{thm:SP} can be extended for $X_{f,\sigma}$.  \\
  We shift our focus below to an interesting special case, namely $X(t)=R Z(t),t \in \TT$, with $R$ a non-negative  rv independent of $Z$, which is 
  the representer of some shift-invariant tail measure $\nu_Z$ on $\AA$. \\
 Assume next that 
$$ 
\limit{x} x^\alpha \pk{R> x } = 1,
$$
for some $\alpha>0$.  Applying Corollary
   \ref{shikon} we have that $X \in \widetilde{\mathcal{R}_\alpha}(a_n, \mcb_0,\nu_Z)$ with $a_n=n^{1/\alpha}$.  Moreover, in view of \cref{app2} or directly by \cite{AAP1579}[Lem 1.1]    	for all $n>0$
   $$ \limit{x} x^{-\alpha}\pk*{ \sup_{t \in [0,n]^l \cap \TT } \norm{X(t)} >x } = \E*{\sup_{t \in [0,n]^l  \cap \TT } \norm{Z(t)}^\alpha} \in (0, \IF).
   	$$
We consider next what happens when $n$ tends to infinity.\\
\BEL
 If $R$ possesses a probability density function $f$ such that $f(s) \le   c s^{-\alpha - 1}$ for some $c>0$ and all $s$ large, then  
\bqn{ \limsup_{n \to \IF} \pk*{ \sup_{t\in [0,n]^l \cap \TT} \norm{X(t)}> a_n  ^l x} \le   C\theta_Y 
	x^{- \alpha}
	\label{rak2}
}
is valid for all $x>0$ and some fix $C>0$.  
\label{malhot}
\EEL 
\begin{example}[Stationary $Z$] Consider  $R$  as in \Cref{malhot} assuming further that $Z$ is  stationary. Clearly, $\tailmeasure_Z$ is shift-invariant and moreover by 
	 \cite{ZKE}[Cor 1, Rem 1,ii)]  we have that  $\theta_Y=0$. 	Consequently, \eqref{rak2} implies \eqref{matura}.
\end{example}	
 \subsection{Open questions} 
RV in the discrete setup $\TT=\mathbb{Z}^l$ has played an important role in the derivation of large deviation type results, as considered in e.g.,  
\cite{MR3078288,buritica2021threshold}. Also for the discrete setup, \cite{kulik:soulier:2020,cissokho2021estimation,MR4280158} have shown that RV and shift-invariant tail measures are crucial for the estimation of various functionals of time series. \\
 The applications of RV to stationary processes are abundant, while the non-stationary case is rather intractable. Even for the simple case $X(t)= RZ(t)$ as in the previous section the asymptotic approximation of $ a_n  ^{-l}M_n $ could not be derived for general $\tailmeasure_Z$. In the recent contribution \cite{Yanis}, periodic sequences have been considered for the discrete setup. \\
With motivation from the aforementioned results and developments we formulate below four open questions. 
\begin{enumerate}[{OP}1)]
\item 
Heavy-tailed large deviation approximation of a   sequence of $\spaceD$-valued random elements $X_n,n\in \mathbb{N}$ can be introduced as in \cite{MR4036028}[Def 2.5]
also in  the general setup of this paper assuming \Cref{A:5}. Specifically, given $\gamma_n, n\ge 1$ that increases to infinity as $n$ tends to $\IF$, we require for a closed cone $F\subset \spaceD$ 
$$
 \gamma_n \pk{ X_n   \in  \cdot  }  \stackrel{v,\mcb_F}{\longrightarrow} \mu(\cdot), \quad n\to \IF,
$$
 with $\mu$ a non-trivial Borel measure.   Since $\gamma_n$ is general, $\mu$ does not need to be a tail measure.   It is of interest to consider non-compact $\TTT$, for instance $\TTT=\R^l$ and the special case where $\mu$ is obtained as a transform of the product of two tail measures as in the results derived in \cite{MR4036028}.  It remains to be investigated if such extensions yield significant applications; 
\item The  applications discussed in   \cite{kulik:soulier:2020,cissokho2021estimation,MR4280158,buritic2021variations,buritica2021threshold} for $\TTT=\mathbb{Z}^l$ can be considered also in the non-discrete setup $\TTT= \R^l$ using additionally our findings related to RV and tail measures. Still several technical conditions in the aforementioned papers  need to be translated to the non-discrete settings, which is not an easy task;  
 \item Does RV of periodic processes on $\TTT=\R^l$ offer some technical advantages in the analysis of related questions   posed in \cite{Yanis}? In particular it is of some interest to  relate (and estimate) the extremal index of periodic processes in terms of the corresponding $Y^{[h]}$'s;  
\item Several findings and applications in \cite{PH2020} are derived based on \Cref{anticlustering} (C($\R^l$)). As shown in \nelem{gio} the weaker 
\Cref{anticlustering} (C($\mathbb{Z}^l$)) can instead be imposed even when $\TTT=\R^l$. It is of interest to investigate if the weaker condition \Cref{anticlustering} (C($\mathbb{Z}^l$)) can be imposed in the applications discussed in \cite{PH2020} and also to characterise stationary processes $X$ for which both conditions are equivalent.    
\end{enumerate}
}
   	
\section{Proofs}
\label{sec:proofs}
\def\TT{\mathcal{T}}
 
\proofprop{pnu} In the proof below we use several times the Fubini-Tonelli theorem which is applicable since $\nu$ is $\sigma$-finite. Since $\norm[t]{\cdot},t\in \TO$ is measurable,  $1$-homogeneous and the outer multiplication $(z,f) \mapsto z \DOT f$ is jointly measurable, then for all maps $\Gamma \in \Hh$  
and  all $h , t\in \TO$ such that $p_hp_t>0$ and all $x >0$, by the definition of $Y^{[h]},Y^{[t]}$  
and  $-\alpha$-homogeneity of $\nu$ 
\bqny{ 	p_h \E*{ \Gamma( x  \DOT Y^{[h]}  ) \ind{ x\norm[t]{ Y ^{[h]}}> 1}} 	&=& 
	\int_{\spaceD}  \Gamma(x\DOT f) \ind{ x\norm[t]{ f}> 1, \norm[h]{ f}> 1 } \nu(df) \\
	&=& 	x^{\alpha}	\int_{\spaceD}  \Gamma(f) \ind{\norm[t]{ f}> 1, \norm[h]{ f}> x } \nu(df) \\
	&=& 	x^{\alpha}	\int_{\spaceD}  \Gamma(f) \ind{\norm[h]{ f}> x, \norm[t]{ f}> 1  } \nu(df) \\
	&=& x^{\alpha} p_t \E*{ \Gamma(   Y^{[t]}  ) \ind{ \norm[h]{ Y ^{[t]}}> x}}.
} 
If  $p_h=0$ and $p_t>0$,  as above taking $\Gamma$  bounded by some constant $C>0$ we obtain 
\bqny{  x^{\alpha} p_t \E*{ \Gamma(   Y^{[t]}  ) \ind{ \norm[h]{ Y ^{[t]}}>   x}}
	&=& 
	\int_{\spaceD}  \Gamma(x\DOT f) \ind{ x\norm[t]{ f}> 1, \norm[h]{ f}\ge  1 } \nu(df)
	\\
& 	\le &	 C \int_{\spaceD}   \ind{\norm[h]{ f}\ge  1 } \nu(df)
	\\
& 	= &	 C \int_{\spaceD}   \ind{\norm[h]{ f}>  1 } \nu(df)
	\\
&= &	 Cp_h\\
&= &0
}
for all $x\in (0,\IF)$, where the third last equality follows from  \eqref{thus}. Since $\norm[h]{\cdot}$ is non-negative we have thus  $\pk{\norm[h]{ Y ^{[t]}}=0}=1$ and further  \eqref{tiltY} holds.
\\
 A direct implication of \eqref{tiltY} is that  $R=\norm[h]{\vk Y ^{[h]}}$ is an $\alpha$-Pareto rv. 
\mb{In particular,}  for all $x\in ( 1,\IF), h\in \TO $ using \eqref{tiltY} and \mb{that} $\pk{ \norm[h]{   Y^{[h]}} >1 >  1/x}=1$  we obtain
\bqny{ \lefteqn{ \pk{ [\norm[h]{  Y^{[h]} } ]^{-1}\DOT \vk Y^{[h]} \in A, \norm[h]{   Y^{[h]}} > x }}\\
	 &=& 
	\E*{ \ind{[\norm[h]{  Y^{[h]}}]^{-1}\DOT  \vk Y^{[h]} \in A} \ind{\norm[h]{   Y^{[h]}} > x } } \\
	&=& x^{-\alpha} \E*{ \ind{   [x \norm[h]{  Y^{[h]}} ]^{-1}\DOT  (x \DOT   Y^{[h]})  \in A} \ind{\norm[h]{   Y^{[h]}} > 1/x } } \\
	&=&  x^{-\alpha} \E*{\ind{ [\norm[h]{     Y^{[h]}}]^{-1} \DOT Y^{[h]} \in A} }, \  \ \forall A \in \AA
}
implying that    $ R$ is independent of $\Theta^{[h]}=  \norm[h]{   Y^{[h]}}^{-1} \DOT  Y^{[h]}$. \eEE{Hence $\pk{\norm[h]{ \Theta ^{[t]}}=0}=1$ for $h,t$ such that $p_h=0, p_t>0$ follows from  $\pk{\norm[h]{ Y ^{[t]}}=0}=1$ shown above and the $1$-homogeneity of $\norm[t]{\cdot}$'s.}\\
By the definition  for all $h\in \TO$ such that $p_h>0$ we have that $\pk{\norm[h]{   \bTheta^{[h]}}=1}=1$   and thus 
\eqref{kksuter} follows.  Further for all $h,t\in \TO$ such that $p_hp_t>0$ and all $\Gamma \in \Hh_0$ by \eqref{tiltY} (recall $\nualpha(dr)= \dualpha$)
\bqny{p_h\E{ \norm[t]{  \tPr}^\alpha \Gamma  ( \bTheta^{[h] })    } 
	&=& p_h\E*{ \frac{ \norm[t]{   Y^{[h]}}^\alpha }  
		{ \norm[h]{   Y^{[h]} } ^\alpha} \Gamma  (  Y^{[h]})} \notag \\
	&=&  \int_{\spaceD} \frac{ \norm[t]{   y}^\alpha }  
	{ \norm[h]{  y} ^\alpha} \ind{ \norm[t]{y}>0} \Gamma (  y) \ind{\norm[h]{   y}>1 } \nu(d  y)\notag \\
	&=& \int_0^\IF   r^\alpha   \int_{\spaceD}  \frac{ 1 }   
{ \norm[h]{r\DOT  y} ^\alpha}   \Gamma ( r y) \ind{r\norm[h]{   y}>r } \ind{r \norm[t]{   y}> 1 } \nu(d  y) \nualpha(dr) \notag \\
	&=& \int_0^\IF   r^\alpha   \int_{\spaceD}  \frac{ 1 }   
{ \norm[h]{r\DOT  y} ^\alpha}  \Gamma ( r y) \ind{r\norm[h]{   y}>r } \ind{r \norm[t]{   y}> 1 } \nu(d  y) \nualpha(dr) \notag \\
	&=& \int_0^\IF   r^{2\alpha }  \int_{\spaceD}  \frac{ 1 }   
{ \norm[h]{   y} ^\alpha}  \Gamma (  y) \ind{\norm[h]{   y}>r } \ind{ \norm[t]{   y}> 1 } \nu(d  y) \nualpha(dr) \notag \\
	&=&   \int_{\spaceD}  \frac{ 1 }   
{ \norm[h]{   y} ^\alpha}   \Gamma (  y) \ind{ \norm[t]{   y}> 1 } \int_0^\IF  \alpha  r^{\alpha -1} \ind{\norm[h]{   y}>r }  dr   \nu(d  y)\notag \\
	&=& p_t  \int_{\spaceD} \frac{1}{p_t}  \ind{ \norm[h]{ y}>0}  \Gamma (  y) \ind{\norm[t]{    y}>1 } \nu(d   y) \notag \\
	&=& p_t\E*{ \ind{ \norm[h]{   \bTheta^{[t]}}  > 0} \Gamma(   \bTheta^{[t]} )},
}
where we used the $-\alpha$-homogeneity of $\tailmeasure$ in the derivation of last forth equality above,  
hence \eqref{tsf} follows. 
Next  	$ Y^{[h]}, h\in \TO:\eEE{p_h>0} $  uniquely define $\nu$, \ie\ we have that 
	$$ \nu( \{ f \in A: \norm[h]{  f} > 1 \})= \tailmeasure(A \cap A_{1h}),\quad  A \in \AA, h\in \TO:\eEE{p_h>0} $$
	 determine $\nu$, which follows from \nelem{triv}.
\QED

\def\THZ{ \bTheta^{[h]} _{ Z}}
\prooflem{Th:Rep}  
{First note that in view of Remark \ref{remNA}, $Y^{[h]}, \Theta^{[h]}, h\in \TO$ and $N$ can be defined in the same probability space $\Pspace$ and all these random elements are independent. \eEE{By assumptions, $q_t$'s and $p_t$'s are such that $\E{p_N}< \IF$.} Note that by   \neprop{pnu} also the $\alpha$-Pareto rv $R= \norm[t_0]{  Y^{[t_0]} } $ is defined on this probability space. Therefore also $Z=Z_N$ is defined on $\Pspace$ (once we show that it is well-defined).} 
Set 
$\SSY(\Theta^{[h]})=1$ if $p_h=0$.  Since $\pk{\norm[h]{   \Theta^{[h]} }=1}=1$ and $q_h>0$ for all $h\in \TO: p_h>0$, then  $\pk{ \SSY(\Theta^{[h]}) >0}=1$ for all $h\in \TO$. 
By the independence of $N$ and $\Theta^{[h]}$'s  
$$ \pk{\SSY(\Theta^{[N]})>0} = \intTO \pk{\SSY(\Theta^{[\mb{t}]})>0} q_t\lambda(dt)=\intTO q_t\lambda(dt)=1, $$ 
where $\lambda=\LK$ is the counting measure on $\TO$. It follows from \eqref{tsf}, that for all $h$ such that $p_h>0$ 
$$ 
\E{\SSY(\Theta^{[h]})}=  \intTO  \E{ \norm[t]{   \Theta^{[h]} }^\alpha }q_t\lambda(dt) \le
\frac{1}{p_h} \intTO  p_t  q_t\lambda(dt) = \E{p_N}/ p_h< \IF
$$
and thus   $\SSY( \Theta^{[N]}) $ is  almost surely positive and finite. Hence  by the cone measurability assumption and the independence of $N$ with  $Y^{[h]}$'s, the random element $Z=Z_N$ is well-defined and by \eqref{kksuter} we have that $\pk{ Z_{\TO}=  0_{\TO}  }=0$ follows. 

 Next  for all $h:p_h>0$,  
let $ \THZ$ be random elements with probability law  given by  
$$  
\mu_h(A)=\E{ \norm[h]{  Z }^\alpha \ind{[\norm[h]{Z}]^{-1} \DOT Z \in A }},\, A \in \AA.$$
{By the joint measurability of the  pairing  $(z, f) \mapsto z \DOT f$ we have that $\mu_h$ is a well-defined probability measure on $\AA$ and further} 
     $\nu_Z$  specified in \eqref{repNU} is also well-defined. Assume next that  $Z,\Theta_Z^{[h]}, \Theta^{[h]}, \forall h\in \TO:p_h>0$ are all defined in the same probability space $\Pspace$. 
Using \eqref{spd} we have (note that  $ \norm[h]{  \THZ }=1,\forall h\in \TO$ almost surely)
$$ \pk{ \THZ \in A}  = \E*{ \ind{ [\norm[h]{  \vk \THZ } ]^{-1} \DOT \THZ \in A } } 
=	\E*{ \ind{ \vk \Theta^{[h]} \in \norm[h]{  \vk \Theta^{[h]} } \DOT  A } } 
= \pk{\vk \Theta^{[h]} \in A}, \ \forall A \in \AA.
$$
Since $\nu$ and $\nu_Z$ have local spectral tail process $\Theta^{[h]}$ and $\Theta^{[h]}_Z, h\in \TO$, respectively, by 
\neprop{pnu}   $\nu=\nu_Z$. \eEE{The latter and \Cref{L:P3} imply  
\bqny{ 0=\nu\Bigl( \Bigl\{ f\in \spaceD: \sup_{t\in \TO} \norm[t]{f} =0 \Bigl\} \Bigr) = \int_0^\IF   
	 \pk*{r\sup_{t\in \TO} \norm[t]{Z}=0} \nualpha(dr)
}
and thus $\pk*{\sup_{t\in \TO} \norm[t]{Z}=0}=0$. Hence using further \eqref{spd} we establish  \eqref{harixh}.} Since $\nu=\nu_Z$  yields
\bqn{\label{twi}
 \tailmeasure \left(\left\{\bsy\in\spaceD:  \sup_{t\in K} \norm[t]{f} >	1\right\}\right) =
\int_0^\IF  \pk*{  r \sup_{t\in K} \norm[t]{Z} > 1}  \nualpha(dr) = 
	\E*{ \sup_{ t\in K  } \norm[t]{\bsZ}^\alpha}   
} 
for all $K \subset \TO$, then \mb{\eqref{borsh} is equivalent to \ref{L:P4} establishing the proof.}
\QED

\prooflem{lem:repA} Let the $\AA/\borelR$ measurable map $ H: D\mapsto \R$ be  {such that for some $\ve_H>0$} and 
all $f\in \spaceD $ we have 
 $H(f)=0$ if $f^*_{K_0} \le \ve_H$.  Hence for all  sets $K$  such that $K_0\subset K\subset \TTT$,  
since   $q_t$'s are positive and { the maps $\norm[t]{\cdot}: \spaceD\mapsto   [0,\IF], \forall 
	t\in \TO$ are 1-homogeneous} 
\bkn{\label{heac} 
	H(f)= H(f) \ind{  f^*_{K} > \ve }=H(f) \ind{\exc_{K}(\ve^{-1}\DOT f) >0}, \quad  \forall \ve\in (0, \ve_H],\forall f\in \spaceD .
} 
If $\lambda=\LK$ is the counting measure on $\TO$, then  $\exc_{K}(\cdot) \in \Hh$. 
Since further by assumption $I=\int_{K} p_t q_t\lambda(dt) < \IF$,  
 by \Cref{L:P1} and the $\sigma$-finiteness of $\tailmeasure$, applying  the  Fubini-Tonelli theorem we obtain 
\bkny{  
	\int_{\spaceD} \exc_{K}(\ve^{-1} \DOT f) \tailmeasure(df) &=&
	 \ve ^{-\alpha} \intTKP \int_{\spaceD} \ind{ \norm[t]{f }> 1 } \tailmeasure(df) q_t\lambda(dt)\\
	&=& \ve ^{-\alpha}\intTKP  p_t q_t\lambda(dt) < \IF
}
and hence  $\tailmeasure(\{f\in \spaceD: \exc_{K}(\ve^{-1} \DOT f) =\IF \})=0$.
By \eqref{mada},\eqref{heac} and  the Fubini-Tonelli theorem  
\bkny{ \tailmeasure[H]&=& \int_{\spaceD} H(f) \tailmeasure(df)=\int_{\spaceD} H(f)  \ind{\exc_{K}(\ve^{-1}  \DOT  f) >0}\tailmeasure(df) \\
&=& \int_{\spaceD} H(f)\ind{\exc_{K}(\ve^{-1} \DOT f) >0} \frac{  \exc_{K}(\ve^{-1} \DOT f) } {\exc_{K}(\ve^{-1} \DOT f)} \tailmeasure(df)\\
&=& \intTKP  \int_{\spaceD} H(f)\ind{\exc_{K}(\ve^{-1} \DOT f) >0} \frac{  \ind{\norm[r]{ f }> \ve} } {\exc_{K}(\ve^{-1} \DOT f)} \tailmeasure(df)
	q_r \lambda(dr)\\
	&=& \ve^{-\alpha}
	\intTKP \int_{\spaceD} \frac{  H(\ve \DOT  f) \ind{ \norm[r]{ y }  >1 }   }  { \exc_{K}( f)} 
	\tailmeasure(df) q_r \lambda(dr)\\
	&=& \ve^{-\alpha} 	\intTKP  \E*{  \frac{ H(\ve \DOT  Y^{[r]})  
		} {\exc_{K}(  Y^{[r]})}} p_r q_r \lambda(dr)
}
establishing the claim. 
\QED 

\def\SKQ{\mathcal{S}_K^q}

\prooflem{capit2}  \Cref{item:toniRA}:  Let $\nu$ \mb{be defined through a} stochastic representer $Z=Z_N$ as  in  \nelem{Th:Rep}. Since $\nu$ satisfies \Cref{L:P1}-\Cref{L:P3}, then  by 
\nelem{triv} $\tailmeasure$ is unique, hence the claim follows. \\
\Cref{item:toniRC}:  The claim follows by showing that  \eqref{twi} holds for all compact $K\subset  \TTT$, which is implied by  Remark \ref{mbab}. 
\cEV{When \Cref{A:2} holds we can use additionally   
Remark \ref{mbab},\Cref{mbab:0}.}\\
\Cref{plusT}: Let   $\Theta^{[h]},h\in \TO$ be a family of spectral tail processes   satisfying \eqref{kksuter} and \eqref{tsf} and set $Y^{[h]}= R\DOT \Theta^{[h]}, h \in \TO$.
For all $h,t\in \TO$ such that $p_hp_t>0$ we have by (\ref{tsf}) that  $\norm[t]{\vk \Theta^{[h]}} $ is positive with non-zero probability and almost surely finite. 
Given $x>0$ and $\Gamma \in \Hh$  by the $1$-homogeneity of 
$\norm[t]{\cdot},t\in \TO$ and the independence of the $\alpha$-Pareto rv $R$ with $\Theta^{[h]}$'s 
(set $B_h= \norm[h]{\vk \Theta^{[t]} }$  and recall that $\pk{B_t=1}=1$  and $\nualpha(dr)= \dualpha$)
\bqny{ 
	\lefteqn{	x^\alpha p_t  \E*{ \Gamma (  Y^{[t]}) \ind{ \norm[h]{ Y^{[t]} }> x}} }\\
	&=& 
	x^{\alpha} p_t
	\int_0^\IF 	\E{ \Gamma (r \DOT  \Theta^{[t]}) \ind{ r B_h > x,  r> 1, 0 < B_h < \IF} }   \nualpha(dr)\\
	&=& 			 
	p_t 	\int_0^\IF 	\E*{ B^\alpha_h  \Gamma  \Bigl((r x)\DOT   \frac{ \Theta^{[t]} }{ B_h} \Bigr) \ind{ r > 1   ,  r \frac{B_t}{B_h}> 1/x, 0 < \frac{B_t}{B_h} < \IF }} \nualpha(dr)\\
	&=& 			 
	p_h 	\int_0^\IF 	\E*{   \Gamma ((r x) \DOT  \Theta^{[h]}  ) \ind{ r > 1   ,  r\norm[t]{ \Theta^{[h]}}> 1/x, 0 < \norm[t]{\vk \Theta^{[h]} } < \IF }} 
	\nualpha(dr)\\  
	&=& p_h	 \E*{ \Gamma ( x \DOT Y^{[h]}  ) \ind{ \norm[t]{ Y ^{[h]}}> 1/x}},
}	
where we used \eqref{tsf}  in the second last  line above and $\pk{ B_h \in [0,\IF)}$ which follows by    Remark \ref{shurdh}, hence the proof is complete. \QED

\prooflem{lem:kulikB22} The claim follows with the same arguments as given in the proof of \cite{kulik:soulier:2020}[Thm B.2.2 (b)]. \QED

\prooflem{lem:engj} For all $z\in \RP,$ if $\Gamma: \spaceD \mapsto \RPP$  is a bounded continuous map and $supp(H) \in \mcb$, then  by assumption \Cref{item:boundProd} and the continuity of the  pairing  $(z,f) \mapsto  z\DOT\, f$, also $ \Gamma_z(f)=\Gamma(z\ \DOT\  f), f\in \spaceD$ is a bounded continuous map supported on $\mcb$ 
for all $z\in \RP$. Consequently,  the assumption $g \in \RA$ implies     
$$ \mu[\Gamma]= \limit{x}\mu_{xz}[\Gamma(xz \ \DOT\  )]= \limit{x} \frac{g(xz )}{g(x)} g(x)\mu[\Gamma(xz \ \DOT\  )]=z^\alpha \mu[\Gamma_z]=\tailmeasure_z[\Gamma], \ \forall z \in \RP.$$
Since $z$ can be chosen arbitrary 
$$ \mu[\Gamma_{s}]= \tailmeasure_{z/s}[\Gamma_s]= s^{-\alpha} z^\alpha \nu[\Gamma_z]= \mb{s^{-\alpha} \tailmeasure[\Gamma]},
\ \forall s \in \RP,
$$
hence the claim \Cref{item:AL1} follows from Remark \ref{lepuruxhi}.\\
By assumption
$H_{t_0}(f)= \norm[t_0]{ f }$ is a $\borelD/\borel{\R}$ measurable function.
The assumption that  $H_{t_0}^{-1}(B) \in \mcb$
for all $B\in \borel{\R}$ with $B \in \mcb_0(\R)$  and $\tailmeasure(Disc(H_{t_0}))=0$ imply in view of \netheo{mapping}  \mb{that}
$$\limit{z} g(z) \pk{ \norm[t_0]{X }\in  z \DOT A}=\tailmeasure(\{f\in \spaceD: \norm[t_0]{f } \in A  \})=: \nu_*(A), \quad \forall A \in \mcb_0(\R) \cap 
\borel{\R}.$$
Since for $A=(x,\IF)$ we have $\nu_*(A)=p_{t_0}(x)\in (0,\IF)$, the measure $\nu_* $ is non-zero and hence the assumption that $g$ is Lebesgue measurable implies that  $g\in \RA$ for some $\alpha>0$ using for instance
\cite{kulik:soulier:2020}[Thm 1.1]. Hence by statement  \Cref{item:AL1} we have that 
$\tailmeasure$ is $-\alpha$-homogeneous and \mb{thus} statement \Cref{item:AL2} \mb{holds}. \\
We show next \Cref{item:AL3}    \mb{along the lines of} \cite{kulik:soulier:2020}[Thm B.2.2]. Since $\nu$ is non-trivial and $\mcb$-boundedly finite, we can find an open set 
$A$ such that  $A \in \mcb$ and $\nu(A) \in (0, \IF)$.
  \mb{Further, by our } assumption  $z\DOT A \subset A$ for all $z\ge 1$ (thus $A$ is a semi-cone). 
  As shown in  \cite{kulik:soulier:2020}[p.\ 521]
  \begin{equation}\label{seep}
  	t\DOT A \subset s \DOT A, \  \  \forall  t\ge s>0.
  \end{equation}
 Consequently,  $z \mapsto \tailmeasure(z \, \DOT \, A) $ is decreasing and by \Cref{item:boundProd} also finite for all $z\in \RP$.  Further by 
 \Cref{item:leben}   
$$t\,  \DOT\,  \overline{A} \subset s\DOT A,\ \  \forall t>s>0$$
 implies that $\tailmeasure(\partial( z\DOT A))=0$ for almost all $z\in \RP$. Hence  \mb{by} assumption we can find some $k>0$ such that $\tailmeasure(\partial(A_k))=0,  A_k=k\DOT A$ and further $\nu(A_k)\in (0,\IF)$. \mb{By}  the continuity of the  pairing  we have that $z\DOT A$ is  open for all $z\in \RP$, \mb{and} then 
$\mu \in \RVP $ implies for almost all $s\in \RP$
\bkn{ \label{diell}
	\mb{\limit{z} \frac{ g(z/s)}{g(z)}= \limit{z} \frac{ g(z/s)}{g(z)} \frac{ \mu((z/s) \DOT (  ks \DOT A))}{\mu(z \DOT ( k\DOT  A))}} = \frac{ \tailmeasure(s \DOT A_k)}{\tailmeasure(A_k)} < \IF,
}
where the last inequality follows since  
 $s\DOT A_k =(k s)\DOT A \in \mcb$ and $\tailmeasure$ is $\mcb$-boundedly finite. Note that since $g$ is non-negative   and 
 $$\limit{z} g(z)  \mu( (kz) \DOT A)= \nu(A_k) \in (0,\IF)$$
  we have that $\mu( (kz) \DOT A)$ is positive and finite for all $z$ large, hence 
 $\mu( (kz) \DOT A)/\mu( (kz) \DOT A)=1$ for all $z$ large justifying the second expression in \eqref{diell}.\\
 Next, by the countable additivity of $\tailmeasure$,  (\ref{seep}) and the assumption  
 $\cap_{z\ge 1}  (z \DOT A)$ is empty we have 
 $$ \limit{z} \tailmeasure(z \DOT A_k)=\tailmeasure(\emptyset)=0$$
  and \mb{by} \eqref{seep} we have that the limit in \eqref{diell} cannot be constant. \mb{Then,} since  $g$ is Lebesgue measurable,  by \cite{kulik:soulier:2020}[Thm 1.1.2] $g\in \RA$ and moreover necessarily  $\alpha>0$, 
 hence  statement  \Cref{item:AL3} follows. 
  \QED

\prooftheo{prophope}  
 Let  $H:\spaceD \mapsto \R, supp(H)\in \mcb$ be a bounded   continuous map. The assumption on $\mcb$  implies that there exists $\ve>0$ and some \eEE{$K \subset\TO $} such that $H(\bsy)=0$, if $\bsy^{*}_{K}=
 \sup_{ t \in K\cap \TO } \norm[t]{\bsy}\leq c\epsilon$ for some fixed given $c>1$. 
Hence we have  
\bkn{ \label{hEps} 
	H(f) = H(f)\ind{\exc_{K}((c\epsilon)^{-1} \DOT f)>0} = \eEE{H(f)\ind{ \sup_{ t \in K\cap \TO } \norm[t]{\bsy}> c\epsilon }}, \quad \forall f\in \spaceD.
} 
Recall that  $\exc_{K}$ is defined by 
$$
\exc_{K}(f ) = \int_{K} \ind{ \norm[t]{f}> 1} q_t\lambda(dt), \quad   f: \spaceD \mapsto \R^d, 
$$
with   $q_t$'s positive constants. Next, for all $\eta,z$ positive, by the Fubini-Tonelli theorem,
\bkny{ \E{ H(z^{-1} \DOT \bsX)}	  
	& \eEE{\ge} & \E{ H(z^{-1} \DOT \bsX) 
			\ind{\exc_{K}(\epsx \DOT  \bsX)>\eta}}  \\
	& =&\E*{ H(z^{-1} \DOT \bsX) 
		\ind{\exc_{K}(\epsx  \DOT \bsX)>\eta} \frac{\exc_{K}(\epsx \DOT  \bsX)   }{\exc_{K}(\epsx  \DOT  \bsX)}}\\  
	& =& \intTKP  \E*{  \frac{H(z^{-1} \DOT\bsX)\ind{\exc_{K}(\epsx \DOT \bsX)>\eta}\ind{\norm[t]{\bsX }>\epsilon z} 
		} {\exc_{K}(\epsx \DOT \bsX)} } 
	q_t\lambda(dt) .
}
Note that for the derivation of the second equality above we have used that $\exc_{K}((\epsilon z)^{-1} \DOT \bsX)$ is finite almost surely, which is consequence of the choice of $q_t$'s since  
$$ \E{ \exc_{K} (\epsx \DOT \bsX) } = \E*{ \intTKP \ind{ \norm[t]{\epsx \DOT\bsX  }> 1  } q_t\lambda(dt)}
\le   \int_{\TO}  \max(1, p_t)    q_t\lambda(dt) < \IF  .
$$
The assumption of the continuity of the   pairing  $(z,f) \mapsto z\, \DOT\, f$ implies that $ H_\ve(f)=H(\ve \DOT f): \spaceD \mapsto \RPP$ 
is also a bounded continuous map. Moreover, by \Cref{item:boundProd}   $H_\ve$ satisfies $supp(H_\ve) \in \mcb$.  
Hence, by the  RV  of $\norm[t_0]{\bsX}$, condition \eqref{klubge},  the continuity of  $H_\ve$ and  the fact that $\exc_{K}(f), f\in \spaceD $ is almost surely continuous with respect to the law of $Y^{[h]}$ (hence \mb{the} continuous mapping theorem can be applied) and 
the dominated convergence theorem, for almost all $\eta>0$ we  obtain
\bkny{
	\lefteqn{   \lim_{z\to\infty} \frac{ \E*{ H(z^{-1} \DOT  \bsX) \ind{\exc_{K}(\epsx  \DOT \bsX)>\eta}}}{\pk{ \norm[t_0]{\bsX }>z}}}\\
	&  =   & \lim_{\ve z\to\infty}   \intTKP  \E*{ \frac{H_\ve (  (\ve z)^{-1} \DOT\bsX)\ind{\norm[t]{\bsX }>\epsilon z} 
			\ind{\exc_{K}(\epsx  \DOT \bsX)>\eta}} {\pk{ \norm[t]{\bsX}>\ve z}\exc_{K}(\epsx \DOT \bsX)} }\\
	& & \times  \frac{\pk{\norm[t]{\bsX }>\ve z} }{\pk{\norm[t_0]{\bsX }>\ve z}} \frac{\pk{\norm[t_0]{\bsX }>\ve z} }{\pk{\norm[t_0]{\bsX}> z}}	q_t\lambda(dt) \\
	& =&
	\frac{1}{\epsilon^{\alpha} p_{t_0}} \intTKP \E*{  \frac{H(\epsilon \DOT \bsY^{[t]})
			\ind{\exc_{K}(\bsY^{[t]})>\eta}} {\exc_{K}(\bsY^{[t]})} } p_t q_t\lambda(dt) . 
}
The monotone convergence theorem  leads to (recall \eqref{hEps}) 
\bkn{   \label{eq:lim-mu}
	\lim_{\eta\to 0}    \epsilon^{-\alpha}  \intTKP \E*{  \frac{H(\epsilon \DOT\bsY^{[t]})
			\ind{\exc_{K}(\bsY^{[t]})>\eta}} {\exc_{K}(\bsY^{[t]})} } p_t q_t\lambda(dt) 
	= 
	\epsilon^{-\alpha}  \intTKP \E*{  \frac{H(\epsilon \DOT \bsY^{[t]}) \COM{ \ind{\exc_{K}(\bsY^{[t]})>0, \exc_{K}(\bsY^{[t]}/2)>0  }} } {\exc_{K}(\bsY^{[t]})} } p_t q_t\lambda(dt)=\nu[H].
}
	By the above and \eqref{eq:TT}  for some $C^*>0$
\bkny{
	\lefteqn{  
		\lim_{\eta\downarrow 0} \limsup_{z\to\infty} 	
		\frac{	\E*{H(z^{-1} \DOT  \bsX)\ind{\exc_{K} ((c\epsilon)^{-1} z^{-1} \DOT \bsX)>0, \exc_{K}(\epsx \DOT \bsX) \leq \eta}}}{\pk{ \norm[t_0]{X}>z}}}	\notag \\
	& \leq   &\mb{ C^* 	\lim_{\eta\downarrow 0} \limsup_{z\to\infty} 	
	 \pk{ \bsX_{K}^* >c\epsilon z,\exc_{K}(\epsx \DOT \bsX) \leq \eta}  }\\
	& = & 0.
} 
Consequently,	 \eqref{eq:lim-mu} yields	
\begin{align*}
	\lim_{z\to\infty} \frac{p_{t_0}\E*{H(z^{-1} \DOT\bsX)}}{\pk{ \norm[t_0]{\bsX }>z}} = \tailmeasure[H].
\end{align*}
\eEE{Since $H$ is bounded, then for some constant $\tilde C>0$ by \eqref{wtnH} and \eqref{hEps}
$$ \lim_{z\to\infty} \frac{p_{t_0}\E*{H(z^{-1} \DOT\bsX)}}{\pk{ \norm[t_0]{\bsX }>z}} \le  \tilde C 
\lim_{z\to\infty} \frac{p_{t_0}\pk{ \sup_{ t \in K\cap \TO } \norm[t]{X}> c\epsilon   }}{\pk{ \norm[t_0]{\bsX }>z}}< \IF
$$
implying  $\tailmeasure[H]< \IF$. 
}
By the assumption  $Y^{[h]}, h \in \TO$ is  a family of tail processes, hence in view of  \nelem{capit2},\Cref{item:toniRA} it defines a unique tail measure $\nu_*$. From \nelem{lem:repA} we have that $\nu_*[H]= \nu[H]$.  In view of Remark \ref{lepuruxhi}, $\nu_ *[H]$ uniquely defines $\nu_*$ for $H$ as chosen above. Consequently,  $\nu=\nu_*$ follows. 
In the above calculations we  choose $q_t$'s positive such that $\sum_{t\in \TO}   \max(1,p_t)q_t <\IF$ and take $K$ such that $t_0\in K$, which is possible. Since 
$Y^{[h]},h\in \TO$ is a family of tail processes and $\nu[H]< \IF$, then   Remark \ref{mbab},\Cref{mbab:0} implies 
that  
 $\nu$ is a $\mcb$-boundedly  finite  non-trivial Borel  tail measure and 
thus the  claim follows by the definition of RV.  
\QED 
   
\prooflem{rusl}
We can extend $\tailmeasure$ to $\AA$ as in \eqref{extendNU}. 
Note that 
(set $f^*_\TO=\max_{t_i\in \TO} \norm[t_i]{f})$  
\bkny{ 
	D  	&=&  \bigcup_{ t_i  \in \TO} \bigcup_{k  \in \Nset} (\{ f\in \spaceD : f_{\TO}^*= 0\} \cup\{f\in \spaceD : \norm[t_i]{ f}> 1 /k\}   )=:  \bigcup_{ t_i  \in \TO} \bigcup_{k \in \Nset} (A_0 \cup A_{ki}) .
}
The joint measurability of the the outer multiplication and the measurability of $\norm[t]{\cdot}$'s yield
$$A_0 \in \AA, \ \  
\ A_{ki} \in \AA,\ \  \forall k\in \Nset,\forall t_i\in \TO.$$ 
Since $\tailmeasure$ is non-zero it follows that \mb{it is impossible that}
$p_{t}(x)=\tailmeasure(\{f\in \spaceD: \mb{\norm[t]{f}}>x\})=0$ for all  $t\in \TO$ and all $x$ in a dense set of $\RPP$.  Hence the assumption that $\tailmeasure$ is non-trivial  implies that for some $t_0\in \TTT,x>0$  
\begin{equation}\label{e:pt0}
	p_{t_0}(x)=\tailmeasure(\{f\in \spaceD: \mb{\norm[t_0]{f}}>x\}) \in (0,\IF).
\end{equation}
Since $\tailmeasure$ is $\mcb_0$-boundedly finite and we have  set $\tailmeasure(\{0\})=0$, it follows utilising further \eqref{heart} that  $\tailmeasure$ is $\sigma$-finite. 
In view of \netheo{sam},\Cref{tev:7} for a countable dense set $\TO$ 
$$\tailmeasure(Disc(\projt))=0, \  \forall t\in \TO.$$
\cEE{Set $\nu_z(A)= g(z) \pk{ z^{-1} \DOT X \in A}, A \in \AA$. By  \netheo{mapping} 
 $\nu_z \circ \projt^{-1}\convvagueO\tailmeasure \circ \projt^{-1}$, whenever $p_{t}>0$ and in particular by \eqref{e:pt0} this holds for $t=t_0$}. This then implies  (use for instance  \cite{kulik:soulier:2020}[Thm B.2.2]) that 
$g \in \RA$ for some $\alpha>0$ and  hence $\tailmeasure$ is $-\alpha$-homogeneous follows from 
\nelem{lem:engj}, statement \Cref{item:AL1}. By the above we 
 can take 
$$g(t)= p_{t_0}/ \pk{\mb{\norm[t_0]{X}}> t}.$$  
Since   $\tailmeasure(Disc(\projt))=0,  \forall t\in \TO$
we have \eqref{klubge} holds and thus $\tailmeasure$ satisfies \Cref{L:P2}. We conclude that $\tailmeasure$ is a tail measure on $\spaceD$ with index $\alpha$. \vEE{In view of \nelem{Th:Rep} we have that $\tailmeasure=\tailmeasure_Z$ and $\pk{\sup_{t \in \TO} \norm[t]{Z}>0}=1$. 
Since  $\norm{x}=0$ if and only if $x=0$ and $Z$ has almost surely \cadlag\ sample paths, 
then  $\pk{Z =0}=0$.} \cEV{The last two claims follow from Remark \ref{mbab},\Cref{mbab:0}}.  
\QED

	\prooftheo{thmA} Since by the choice of $U_n$ we have $\nu(disc(\proj_{U_n}))=0$, the first implication is direct consequence of the continuous mapping theorem 
	(utilising \netheo{sam},\Cref{tev:7}-\Cref{tev:9}) and the characterisation of $\mcb_0$ and $\mcb_0(\spaceD_{ U_n})$. 
	 { In particular $\nu^{(n)}= \nu \circ \mathfrak{p}_{U_n}^{-1}$ and it follows easily that the local tail processes of $\nu^{(n)}$ denoted by $Y_n^{[h]},h\in 
		\TO_n=\TO \cap U_n$ satisfy 
		$$Y^{[h]}_n = Y^{[h]}_{U_n}$$
	almost surely, where $Y^{[h]}$'s are the local tail processes of $\nu$ and $Y^{[h]}_{U_n}$ is their restriction on $U_n$.
} 

We show next the converse assuming 
	$X_{U_n} \in \mathscr{R}_\alpha(\mcb_0(\spaceD_{ U_n}),\nu^{(n)})$ for all $n\in \Nset$.\\
	\underline{Step I (existence of $\nu$)}:\\	 The sets $U_n$ are increasing and    $\bigcup_{n=1}^\IF  U_n =
	\R^l$. Each measure $\nu^{(n)},  n\in \Nset$ is \eEE{$K(\TO_n)$-bounded (or compactly-bounded)} with $\TO_n=\TO \cap U_n$ and has a unique family of corresponding tail processes $Y^{[h]}_{n}, h\in \TO$. 
	Since all spaces $\spaceD_{U_n}, \spaceD$ are Polish we can consider all local tail processes to be defined on the same non-atomic complete probability space, \cite{MR100291}[Lem p.\ 1276].\\ 
	Applying the continuous mapping theorem (we utilise \netheo{sam},\Cref{tev:7}-\Cref{tev:9}) to the projection of $U_{n+1}$ to $U_n$ denoted by 
	$\mathfrak{p}_{U_{n+1},U_n}$ shows that $ \nu^{(n)} = \nu^{(n+1)}\circ \mathfrak{p}_{U_{n+1},U_n}^{-1}$. It follows that 
	 the restriction  of tail processes $Y^{[h]}_{n+1}$ of $\nu^{(n+1)}$ on $U_n$ denoted by $Y^{[h]}_{U_n},  h\in U_{n} $
	are also tail processes. By the uniqueness of the family of the tail processes it follows that $\nu^{(n)}$ has local tail processes $Y^{[h]}_{U_n}, h\in \TO$, i.e., 
	$Y^{[h]}_n = Y^{[h]}_{U_n}$
	almost surely.  We can extend all $Y^{[h]}_n$'s to be \cadlag\ processes on $\spaceD$. Applying  \netheo{sam},\Cref{tev:4} or \netheo{sam},\Cref{tev:10} we obtain that 
	$Y^{[h]}_n$ converges weakly  on $\spaceD$ as $n\to \IF$ to a $\spaceD$-valued random element 
	$Y^{[h]}$. It follows easily that $Y^{[h]}$ restricted on $U_n$ coincides almost surely with $Y^{[h]}_n$ and moreover $Y^{[h]}, h\in \TO$ is a family of tail processes.  Let $\nu$ denote the corresponding tail measure defined by \eqref{eq:goodschl}. By the definition of local tail processes and the above we  have that 
	\begin{equation}\label{xhinxh}
		\nu_{\lvert U_n}=\nu^{(n)}=\nu_{Z_{n}},\quad n\in \Nset,
	\end{equation}
where $Z$ is a representer of $\nu$ constructed from $Y^{[h]}$'s  {and $Z_{n}$ is given by $\eqref{zn}$.}
Hence  we have 
\begin{equation}\label{ja}
	\E{(\norm[0]{Z_{n}})^\alpha}= \nu^{(n)}(\{f\in \spaceD_n: \norm[0]{f}> 1\} ) 
	\to  	p_0=\E{(\norm[0]{Z})^\alpha}< \IF, \quad n\to \IF. 
\end{equation}
 It follows that $Y^{[h]}$ satisfies \eqref{ifo} (since that holds for $Y^{[h]}_n$'s) and hence $\nu$ is $\mcb_0$-boundedly finite. \\
	\underline{Step II (RV  of $X$)}:\\ 
	By (\ref{heart}) and the definition of $\norm[t]{\cdot}$, see \Cref{A:2} and \eqref{heart} the boundedness $\mcb_{0}$  on $\spaceD_0$   can be  generated  (see also \cite{kulik:soulier:2020}[Example B.1.7]) by the open sets 
	$$O^\IF_{k}=\Bigl \{ f\in \spaceD: \sup_{t\in [-k,k]^l\cap \TO} \mb{\norm[t]{f}}> 1/k  \Bigr\}, \ \ k\in \Nset .$$
	Since $\nu$ is $-\alpha$-homogeneous by Remark \ref{simonDiv}   $\nu(\partial O^\IF_k)=0$ for all $k \in \Nset$.\\
By \eqref{ja} we can assume without loss of generality that  $p_{0}=1$. 
In view of 
Remark \ref{nock},\Cref{nock:1} and recalling that $\nu$ is $\mcb_0$-boundedly finite,  the claim follows if we show the following weak convergence: 
\begin{equation}\label{posht}
	\mu_{k,z}(\cdot )= \pk{ z^{-1}  \DOT X \in \cdot  \cap {O^\IF_k} }/\pk{ \mb{\norm[0]{X}} > z} \convweak \nu(\cdot \cap {O^\IF_k})=:\nu_k(\cdot), \  z\to \IF
\end{equation}	
  for all positive integers $k$. Note that $\nu_k$ is a finite measure and set 
 	$$B_k^n=\Bigl\{ f\in \spaceD_n: \sup_{t\in [-k,k]^l\cap \TO} \mb{\norm[t]{f}}> 1/k  \Bigr\}. $$
Next, fix an integer $k>0$. In the light of  \netheo{sam},\Cref{tev:10} the stated  weak convergence is equivalent to 
$$\mu_{k,z}^*(\cdot)= \pk{ z^{-1}  \DOT X_{U_n} \in \cdot  \cap {B^n_k} }/\pk{ \mb{\norm[0]{X}} > z} \convweak \nu_k^*(\cdot), \   z\to \IF, $$
where 
$$\nu_k^*(A)= \nu_k(\{f\in D: f_{U_n} \in A\}), \  \  A \in \borel{D_{U_n}}$$
 for all $n$ large.
	The properly localised boundedness $\mcb_0(U_n)$   can be generated  by the open sets on $\spaceD_n$ given by 
	$$O_k^n=\Bigl\{ f\in \spaceD_n: \sup_{t\in U_n\cap \TO} \mb{\norm[t]{f}}> 1/k  \Bigr\}, \ \ n\in \Nset. $$ 
In particular, for all $n$ large, $B_k^n \subset O_k^n$, implying   $B_k^n \in \mcb_0(U_n)$. Moreover, $\nu^{(n)}(\partial B_k^n )=0$, since $\nu^{(n)}$ is $-\alpha$-homogeneous and so we can apply Remark \ref{simonDiv}.  By the assumption 	$X_{U_n} \in \mathscr{R}_\alpha(\mcb_0(U_n),\nu^{(n)})$  
Remark \ref{nock},\Cref{nock:1} implies 
 the weak convergence 
 $$\mu_{k,z}^*(\cdot) \convweak  \nu^{(n)}(\cdot \cap B_k^n )=\nu_k^*(\cdot), \ z\to \IF,$$
 where the last equality above follows from  \eqref{xhinxh} establishing the proof.  \QED

\proofkorr{shikon}  It follows as in the proof of \cite{RDTM:2008}[Thm 3.3]  that  $X$ has almost surely sample paths in 
$D(\R^l,\R^d)$. By \cite{RDTM:2008}[Lem 3.1, Thm 3.3] we have that $X_{U_n}$ is regularly varying with $U_n$ as in \netheo{thmA}. Moreover, 
$X_{U_n}$ has de Haan representation \eqref{eq1} with $Z^{i}$'s independent copies of $Z_n$ determined in \eqref{zn}. The claim follows from the converse in the aforementioned theorem and Remark \ref{enu},\Cref{enu:2}. Note that the case $l=1$ is already shown in  \cite{PH2020}[Thm 4.1] using a direct proof. \QED

\prooftheo{Th:RV}
 \Cref{item:RVD} $\implies$ \Cref{item:tightness}: Since $\tailmeasure$ is $\TTTK$-bounded (and also $\mcb_0$-boundedly finite) it has a $\spaceD$-valued representer $Z$ that satisfies \eqref{borsh}, which in view of \Cref{Th:Rep} is equivalent with \Cref{L:P4}. As mentioned for instance in 
 \cite{MR0402847}[p.\ 205] the set of stochastic continuity points of $Z$  denoted by $Z_{\mathbb{P}}$, i.e., 
 $t\in Z_{\mathbb{P}}$ such that  $\pk{Z(t)\not=Z(t-) }=0$ is the same as the set of points $\{ t\in \TTT: \mathbb{P} \circ Z^{-1} (\{f\in \spaceD: \proj_t  \text{ is continuous at } f \} )=0$, \ie\  
$ \proj_t: \spaceD \mapsto 
 \R^d$  is continuous almost everywhere $\mathbb{P} \circ Z^{-1}$. Hence for all $t\in Z_{\mathbb{P}}$  we have 
 $$\tailmeasure(\{f\in \spaceD: f(t)\not= f(t-) )=\int_0^\IF \pk{ Z (t)\not= Z(t-) }   
  \nualpha(dr)= 0 $$
and thus  $\proj_t, t\in Z_{\mathbb{P}}$ is $\tailmeasure$-continuous almost everywhere.\\ 
Let $\TO \subset Z_{\mathbb{P}}$ be a dense set of \mb{$\TTT$} and let $a< b, a,b\in Z_{\mathbb{P}}$ be given and set $K=[a,b]$.   The existence of $T_{a,b} \subset K$ which is up to a countable set equal $K$ such that  
 \eqref{eq:fidi-rv} holds  for all $t_1,\dots,t_k\in  \tE{T_{a,b}} , k\ge 1$  follows \mb{by} arguments mentioned  in \cite{PH2020} where the stationarity has not been used and the proof relies on \cite{hult:lindskog:2005}[Thm 10, (ii)$\implies$(i)]. 
 \eEE{ In the rest of the proof, by the equivalence of  the norms on $\R^d$ we shall suppose without loss of generality that $\norm{\cdot}$ equals the norm $\norm[*]{\cdot}$ on $\R^d$ used also for the definitions of $w'$ and $w''$ below.}\\ 
  Taking $T_0$ to be the union of $T_{a,b}$'s, then \eqref{eq:fidi-rv} holds  for all $t_1,\dots,t_k\in  \tE{T_{0}} , k\ge 1$, with $T_0 \subset \TTT$ such that $\TTT \setminus T_0$ is countable. Moreover, from \cite{hult:lindskog:2005}[Eq.~(7),(8),(9)] and \cite{billingsley:1999}[Eq.~(12.32)] we obtain   \eqref{eq:tightness}. \\
  \Cref{item:tightness}$\implies$\Cref{item:conditional}: Condition \eqref{eq:tightness2A} follows immediately from \eqref{eq:tightness}. 
  For all $h\in T_0, \ve>0, z>0$   
$$\pk{ w'(\bsX,K,\delta) >  \ve z   , \mb{\norm[h]{\bsX}}>z } \leq \pk{ 
	w'(\bsX,K,\delta) > \ve z  }$$
and thus if $p_h>0$,  in view of \eqref{eq:tightness} and \eqref{klubge} we obtain 
\bkn{      \label{eq:tightness-conditional}
	\lim_{\delta\downarrow 0} \limsup_{x\to\infty} \pk{ w'(\bsX,K,\delta) >\ve  z   \mid \mb{\norm[h]{\bsX}}>z } = 0.
}  
Using \cite{PH2020}[Thm B.1], \eqref{eq:tail-process-time-series-tailmeasure}  and \eqref{eq:tightness-conditional} 
{we obtain the weak convergence in $(\spaceD, \Dtau)$ of 
	$z^{-1}\DOT \bsX$  conditionally on $\mb{\norm[h]{\bsX}}>z$} to $\bsY^{[h]}$ and further the limiting process $Y^{[h]}$ has almost surely paths in $\spaceD$. Since $\spaceD$ is Polish, then by  \cite{MR100291}[Lem p.\ 1276]  $Y^{[h]}$ can be realised as a random \mb{element}  on 
the  non-atomic probability space $\Pspace$.\\
In order to establish the claim we need to show that $Y^{[h]}, h\in \TO$  are local tail processes of a  tail measure $\nu$ on $\AA$.  By \neprop{pnu} we have that $Y^{[h]}(t),t\in \TO$ are local tail processes when we take $\TTT=\TO$, \ie\  
\eqref{kksuter} and \eqref{tiltY} hold.  Since by \netheo{sam},\Cref{tev:2} $\sigma(\projt, t\in \TO)=\borelD$  
 it follows that for all $h,t\in \TO$ such that $p_hp_t>0$  \eqref{tiltY} holds also for $Y^{[h]}(t),t\in \TTT$ and thus $Y^{[h]}, h\in \TO$ are local tail processes as $\spaceD$-valued random elements. Hence by \nelem{capit2},\Cref{item:toniRA}, there exists a unique tail measure $\nu$ corresponding to these tail processes.\\
\Cref{item:conditional}   $\implies$\Cref{item:RVD}:
{
	 Let $\TO \subset T_0$ be a dense subset of $\TTT$} and let $w',w''$ be the two moduli  on $\spaceD$ defined in (\ref{wp})),(\ref{wpp}) below. 
For all $z$ sufficiently large and $\delta, \ve $ positive 
\bkny{ 
 	\pk{ X^*_{K}> \ve z}  & \le &  	 \pk{ w'(\bsX,K,\delta)>  \ve z/2}
 	+ \pk{w'(\bsX,K,\delta)\le   \ve z/2, X^*_{K}> \ve z}\\
& \le &  	 \pk{ w'(\bsX,K,\delta)>  \ve z/2}
+ \pk{ \max_{0 \leq i \leq  m} \mb{\norm{\bsX(t_i)}}>  \ve z/2}
}
 	for all  $K=[-n,n] , n\in \Nset$ and every sequence
 	$(t_0,\dots,t_m){\in \TO^m}$ such that $-n=t_0<\cdots<t_l=n$ and $t_{i}-t_{i-1}\leq\delta$ for $i\le m$ since 
 $w''$ is dominated by $w'$, which is shown in \cite[Eq.~(12.28)]{billingsley:1999}.
 By the assumptions \eqref{eq:tightness-conditional} holds, which together with \eqref{eq:tightness2A} implies \eqref{eq:tightness}. Consequently,   the arbitrary choice of $\delta>0$ yields  (\ref{wtnH}). \eEE{In particular, (\ref{wtnH}) implies that 
\begin{equation}\label{CC}
p_t  \le C p_{t_0}< \IF, \quad \forall t\in K
\end{equation}
for some constant $C>0$}. 
  
A crucial implication  of \neprop{pnu} is that $Y^{[h]}$ has the same law as 
$R \DOT \Theta^{[h]}$ with $R$ an $-\alpha$-Pareto rv independent of the local spectral tail process $\Theta^{[h]}$. 
This implies that $\exc_{K}(f), f\in \spaceD$ is almost surely continuous with respect to the law of 
$Y^{[h]}$ (see also \cite{kulik:soulier:2020}[Rem 6.1.6]). Recall that   
$$
\exc_{K}(f ) = \int_{K} \ind{ \mb{\norm[t]{f}}> 1} q_t\lambda(dt), \quad   f: \TTT \mapsto \R^d,
$$
with $\lambda=\LK$ counting measure on $\TO$ 	and  we take $q_t>0, \, t\in \TO$ satisfying $\sum_{t\in K}\eEE{ \max(1,p_t)} q_t< \IF$. 
\eEE{In view of \eqref{CC} the last condition is satisfied if we show that positive $q_t$'s can be chosen such that $\sum_{t\in \TO} q_t < \IF$.}

The proof of the claim follows by showing that condition \eqref{eq:TT} in \nelem{prophope} holds for $c=2$ and appropriate $q_t$'s constructed below.\\
 { 
Consider therefore the following construction of a density $q$, which is needed for the proof below. Consider a compact set  $K \subset [0,1]$.   For fixed $k\in\mathbb{N}$, we pick a single arbitrary and distinct point from each of the sets
\begin{align*}
s_j^{(k)}\in \overline B(m 2^{-k}, 2^{-k})\cap \TO, \quad m\in 0,\dots, 2^k,
\end{align*}
where $\overline B(a, r)$ denotes the closed ball with center in $a$ and radius $r$.
Notice that for any $k$, the above system of balls covers $\TO\cap[0,1]$. 
Assign the same point density 
\begin{align*}
q_s=2^{-2k -2}
\end{align*} 
to each of the $2^k$ distinct points $s_j^{(k)}$. The sum of all these masses is equal to 
\begin{align*}
\sum_{k=0}^\infty\sum_{j=0}^{2^k}2^{-2k-2}=\sum_{k=0}^\infty2^{-k-2}=\frac{1}{2}.
\end{align*}
Spreading out the remaining mass $1/2$ among the non-chosen points in $K\cap\TO$, we obtain
$$\intTKP 	q_t\lambda(dt)=1.$$
Now consider an interval of length $v-u$ with $u,v\in K\cap \TO$. Let $n_1$ be the smallest natural number such that $2^{-n_1}<v-u$. In particular it follows that $(v-u)/2\le 2^{-n_1}$. Hence, there exists a ball $B_1=\overline B(m_1 2^{-n_1-2},2^{-n_1-2})\subset (u,v)$. In particular, there is $s_1\in B_1\cap\TO$ having mass
\begin{align*}
q_{s_1}=2^{-2(n_1+2)-2}=2^{-6}(2^{-n_1})^2\ge2^{-6}\left(\frac{v-u}{2}\right)^2.
\end{align*}
For a general $K$ not included in  $[0,1]$, a constant $C(K)$ can similarly be found such that $s_1\in [u,v)$ and 
$$q_{s_1}\ge C(K)(v-u)^2.$$
If $f\in\spaceD$ is such that $\tE{w'}(f,K,\eta) \leq \epsilon/2$, and
  $\exc_{K}((2\epsilon)^{-1}\, \DOT \, f) > 0$, then there exists $t\in [a,b)$ such that $f(t)>2\epsilon$ and
  an interval $[u,v)$ such that 
  $$t\in[u,v), \ v-u \geq \eta, \ \ \sup_{[u\leq s,s'<v)} \norm[*]{f(s)-f(s')}\leq \epsilon.$$
   Consequently, $f(s)>\epsilon$ for all
  $s\in[u,v)$ and  
\bqny{  \exc_{K}(\epsilon^{-1}f)&=& \intTKP \ind{\norm[t]{f}>\epsilon} q_t\lambda(dt)= \intTKP \ind{\norm[*]{f(t)}>\epsilon} q_t\lambda(dt)\\
 &\ge &\int_{ [u,v)}q_t\lambda(dt) \ge  C(K)(v-u)^2
 \ge   C(K)\eta^2.
}
\tE{Since $X$ has almost surely locally bounded sample paths,  the above yields for some constant $ \tilde C>0$} 
\bkny{
    \pk{ \bsX_{K}^* >2\epsilon  	 z,\exc_{K}(\epsx \DOT\bsX) \leq \eta} 
    & \leq &   \tilde C  \pk{\tE{ w'}(\bsX,K,\sqrt{\eta/C(K)})>\epsilon z},
} 
hence \eqref{eq:tightness} implies condition \eqref{eq:TT} for $c=2$ and thus the claim follows from \nelem{prophope}.\\ 
}
   \Cref{item:conditional} $\implies $\Cref{item:proj}: As shown above when    \Cref{item:conditional} holds,  the tail measure $\nu$ with index $\alpha$     is $\mcb_0$-boundedly finite. By \nelem{Th:Rep}  and Remark \ref{mbab},\Cref{mbab:0} $\nu=\nu_Z$ satisfying further \Cref{L:P4} which is equivalent with \eqref{borsh} as mentioned above. In particular $\nu_{k}$ with   representer $Z_{K_k}$ is a $\mcb_0(\spaceD_{K_k})$-boundedly finite tail measure with index $\alpha>0$ on 
   $\borel{\spaceD_{K_k}}$ for all compact intervals $K_k \subset \R$ containing some open interval that includes 0. \\   
 Let 
 $s_k<t_k, k\in \Nset$  in  $T_0$ satisfying   
 $$-\limit{k}s_k=\limit{k}t_k=\IF.$$
  Suppose for simplicity that $-s_k, t_k,k\ge 1$ are strictly increasing positive sequences and chose them to belong to $Z_{\mathbb{P}}$. This is possible since $T_0$ and $Z_{\mathbb{P}}$ are \mb{equal} up to a countable set.  In view of  \eqref{eq:tightness} we have that for all $\ve>0$
\bkny{   
	\lim_{\tE{\eta}\downarrow 0} \limsup_{z\to\infty} \frac{ \pk{  w''(\bsX, K_k,\tE{\eta}) > \ve  z  }}{\pk{ \mb{\norm[t_0]{\bsX}}>z}} = 0
}
and thus using further \cite{billingsley:1999}[Eq.~(12.31)] and the definition of $\norm[t_0]{\cdot}$ as well as the equivalence of the norms on $\R^d$
\bkny{   
	\lim_{\tE{\eta}\downarrow 0} \limsup_{z\to\infty} \frac{ \pk{  \norm{X(s_k+\eta)- X(s_k) } > \ve  z  }}{\pk{ \mb{\norm[t_0]{\bsX}}>z}} = 
		\lim_{\tE{\eta}\downarrow 0} \limsup_{z\to\infty} \frac{ \pk{  \norm{X(t_k-)- X(t_k-\eta) } > \ve  z  }}{\pk{ \mb{\norm[*]{\bsX(t_0)}}>z}} =0. 
}
Since by \cite{billingsley:1999}[p.~132] almost surely  
$$w(X,[s_k, s_k+\eta),\eta) \le 2  [  w''(\bsX, K_k,\tE{\eta}) + \norm[*]{X(s_k+\eta) - X(s_k) }],$$
$$w(X,[t_k-\eta, t_k),\eta) \le 2  [  w''(\bsX, K_k,\tE{\eta}) + \norm[*]{X(t_k-) - X(t_k-\eta) }],$$
then it follows as in the proof     \Cref{item:conditional}  $\implies$  \Cref{item:RVD} (\mb{along the lines of} \cite{hult:lindskog:2005}[Thm 10,(i)]) that $X_{K_k} \in \mathcal{R}_\alpha( \mcb_0(\spaceD_{K_k}), \nu_{k}).$ \\
 \Cref{item:proj}   $\implies$\Cref{item:tightness}:  The proof follows from   \cite{hult:lindskog:2005}[Thm 10,(ii)] and \cite{billingsley:1999}[Eq.~(12.32)].
 \QED

\prooflem{chans}  Let $U_n,n\in \Nset, \nu^{(n)}$ be as in \netheo{thmA}. By  \netheo{mapping}	$X \in \widetilde{\mathcal{R}}_\alpha(a_n, \mcb_0,\nu)$ 
implies 	$X_{U_n} \in \widetilde{\mathcal{R}}_\alpha(a_n, \mcb_0,\nu^{(n)})$ for all $n\in \Nset$. The Polish space $\spaceD_{U_n}$ is a star-shaped metric space and thus  \cite{hult:lindskog:2006}[Thm 3.1] implies  $X_{U_n} \in  {\mathcal{R}}_\alpha(a_n, \mcb_0,\nu^{(n)})$, hence the claim follows from \netheo{thmA}.
\QED

\cEG{
\def\TT{\TTT}
\proofprop{genUB}  
Note first that by the assumption on $\norm{\vk X(0)}$ for any $c>0$ we have 
$$\lim_{n\to \IF} n^c \pk{\norm{\vk X(0)}>(a_n x)^c } = x^{-c/\alpha}$$
 for all $x>0$.
We consider for simplicity only the case $\TT=\R^l$. By the stationarity of $\vk X$, using \cite{GeoSS}[Thm 2.1]  we obtain 
\bqny{ 
	\lefteqn{ \limit{m}   \limit{n} n m^{-l} \pk*{\sup_{t\in [0,m]^l  }  \norm{\vk X(t)} > a_n }}
	\\
	&=& \lim_{m\to \IF}  m^{-l}   \lim_{n\to \IF} n \pk{\norm{\vk X(0)}>a_n } \int_{t \in [0,m]^l }  \E*{  
		\frac{1}{\int_{s \in [0,m]^l } \ind{ \norm{\vk X(s)}>a_n  } \lambda(ds)}	\Bigl \lvert \norm{\vk X(t)}>a_n} \lambda(dt)\\
	&=& \lim_{m\to \IF} m^{-l}   \lim_{n\to \IF}  \int_{t \in [0,m]^l }  
	\E*{   
		\frac{1}{\int_{s \in [0,m]^l } \ind{ \norm{B^t \vk X(s)}>a_n  } \lambda(ds)}
		\Bigl \lvert \norm{\vk X(0)}>a_n } \lambda(dt) \\
	&\le & \lim_{k\to \IF} \lim_{m\to \IF} m^{-l}   \limsup_{ n\to \IF}  \int_{t \in [k,m]^l }  
	\E*{ \frac{1 }{\int_{s \in [-k,k]^l } \ind{ \norm{\vk X(s)}>a_n  }\lambda(dt)}  				\Bigl \lvert \norm{\vk X(0)}>a_n } \lambda(dt) \\ 
	&= & \lim_{k\to \IF} \lim_{m\to \IF} m^{-l}     \int_{t \in [k,m]^l }  
\E*{ \frac{1 }{\int_{s \in [-k,k]^l } \ind{ \norm{\vk Y(s)}>1  }\lambda(dt)}  } \lambda(dt) \\ 
	&=& \lim_{k\to \IF}  \E*{  \frac{1}{ \int_{s \in [-k,k]^l  } \ind{ \norm{ Y(s)}>1 } \lambda(ds)} } \\
	&=&  \E*{  \frac{1}{ \int_{t \in \R^l  } \ind{ \norm{ Y(t)}>1 } \lambda(dt)} }, 
}
where the third  last line follows from the weak convergence of $X/a_n$ conditioned on $\norm{\vk X(0)}>a_n$ to the tail process $Y$ and  continuous mapping theorem and the second last line is consequence of dominated convergence theorem. Using again the stationarity of $X$ and the above bound gives 
(write $[n/m]$ for the largest integer smaller than $n/m$)
\bqny{ 
	\limsup_{n\to \IF} \pk{ M_{n}> a_n^l x}& \le&  \limsup_{m\to \IF} \limsup_{n\to \IF} ([n/m]+1)^l
	\pk{ M_{m}>  a_n ^l x} \\
	&= & 
	\lim_{m\to \IF} \frac 1 {m ^{l}} \lim_{n\to \IF}  n^l  
	\pk{ M_{m}>  a_n^l x}\\
	&= & 
	\E*{  \frac{1}{ \int_{t \in \R^l  } \ind{ \norm{ Y(t)}>1 } \lambda(dt)} } x^{-\alpha}.
} 
The finiteness of the expectation above follows from \eqref{ifo}, hence the proof is complete.  
\QED
}

\prooflem{homogeneousTransf}  
The claim under the first two conditions follows by \netheo{mapping}  and identical arguments as in\cite{MR3271332}[Cor 2.1, 2.2]. The last claim is again consequence of \netheo{mapping}    if we show that 
\begin{equation}\label{nGST}
	H^{-1}(B) \in \mcb, \ \  \forall B \in \mcb' \cap \borel{D'} \ \text{  and  } \tailmeasure(Disc(H))=0.
\end{equation}
If $B'\in \mcb'$, then by definition 
$\Dtau'(f',F')>\ve$ for all $f' \in B'$ and some $\ve>0$. By continuity of $H$  for all $f\in F$ 
we can find $\delta_f>0$ 
such that for all $x\in \spaceD$  satisfying 
$$\Dtau(f,x) \le \delta_f$$  we have $\Dtau'(H(f),H(x))\le \ve$. Since $F$ has finite number of elements, then $\delta=\min_{f\in F}\delta_f>0$. Since $H(F)=F'$ and $\Dtau'(f',F')>\ve$ for all $f' \in B'$, then  $\Dtau(F,H^{-1}(B'))> \delta$ and thus 
\eqref{nGST} follows establishing the claim. Note that when $F$ has one element this is already proven in \cite{hult:lindskog:2006}[(ii), p.\ 125]. \QED

\prooflem{lem:gio2} Since $\pk{Z\not=0} = \pk{\norm{Y(0)}>1} =1$, then the integrals in \eqref{gio} are almost surely positive. 
Clearly, RV  of $X$ implies the RV of $X(t),t\in \mathbb{Z}^l$ with limit measure which is shift-invariant and has tail process 
$Y(t), t\in \mathbb{Z}^l$. As shown in \cite{PH2020}[Lem 3.5] 
we obtain 
\bqny{
	\norm{Y(t)}\to 0, \quad \norm[\IF ]{t}\to \IF, \quad t\in \mathbb{Z}^l
} 
almost surely, which in view of \cite{klem}[Prop 2.18] is equivalent
with 
\bqn{ \label{melone} 
	\pk*{ \sum_{t\in \R^l} \norm{ Y(t)}^\alpha < \IF}=1.
}  We have that  $Y^*(t)=\norm{ Y(t)}, t\in \mathbb{Z}^l$ is a tail process with representation $R \Theta^*(t)=R \norm{ \Theta^{[0]}(t)}, t\in \mathbb{Z}^l$, where $R$ is $\alpha$-Pareto independent of $\Theta^*$, which is a spectral tail process. Hence   in view of \cite{ZKE}[(4.6)]  we have that \eqref{melone}  is equivalent with $\pk{\int_{\R^l} \norm{ Y(t)}^\alpha \lambda(dt) <  \IF}=1$ establishing the claim.  
\QED 

\cEG{
\prooftheo{thm:SP} The RV of $X$ being max-stable has been shown in \Cref{shikon}. If $X$ is $\alpha$-stable RV can be established as in \cite{PH2020}. Alternatively, since for this case \Cref{enu},\ref{enu:2} holds and RV of $X$ for compact $\TTT$ has been established in \cite{RDTM:2008}, the RV of $X$ follows from \netheo{thmA}. If $\theta_Y=0$, then \eqref{jugendlich} follows from \eqref{raki} and when  $\theta_Y>0$  the corresponding proofs in \cite{PH2020} can be modified to cover the case $l>1$. \QED
}
\def\TT{\TTT}

\cEG{\prooflem{malhot} For all $x>0,n \in \mathbb{N}$  and all  $ y>0$ large  
\bqny{  \pk{ a_n  ^{-l}M_n  > x} &=& 
	\pk{ a_n  ^{-l}M_n   > x , R\le y}+ \pk{ a_n  ^{-l}M_n   > x , R> y}
	\\
	&\le & \pk*{\sup_{t \in [0,n]^l\cap \TT} \norm{Z(t)}  \ge  a_n^l x/y}+\pk{ M_n  > a_n^{l} x , R> y} \\
	&\le & \frac{y^\alpha }{n^l x^\alpha}\E*{\sup_{t \in [0,n]^l\cap \TT } \norm{Z(t)}^\alpha}+ c\int_y^\IF s^{-\alpha - 1}\pk*{ \sup_{t \in [0,n]^l \cap \TT} \norm{Z(t)}> a_n^{l}x/s } ds,    
}
where we used the Markov inequality and the assumption $f(s) \le c s^{-\alpha-1} $ for all $s$ large  to derive the last line above. The shift-invariance of  $\tailmeasure_Z$ implies 
\eqref{tcfN}. Hence as in \cite{debicki2017approximation, PH2020} it follows that 
$$  \limit{n} n^{-l} \E*{\sup_{t \in [0,n]^l\cap \TT} \norm{Z(t)}^\alpha}=\E*{\sup_{t \in \TT} \norm{Y(t)}^\alpha / \int_{\TT} \norm{ Y(t)}^\alpha \lambda(dt)   }= \theta_Y.$$
Given $\ve>0$ for all large $y$ we have that 
\bqny{ 
	\int_y^\IF  s^{-\alpha - 1}\pk*{ \sup_{t \in [0,n]^l\cap \TT} \norm{Z(t)}> a_n^{l} x/s } ds &<& 
	 (1+\ve) \int_y^\IF  s^{-\alpha - 1} e^{-s^{-\alpha}} \pk*{ \sup_{t \in [0,n]^l\cap \TT} \norm{Z(t)}> a_n^{l}x/s } ds\\
	 & = & (1+\ve) \alpha^{-1}\pk{M_n^* > a_n^lx},
}
where $ M_n^*$ is defined as $ M_n$ taking  $R=\Gamma_1^{-1/\alpha}$ with $\Gamma_1$ a unit exponential rv. 
If $\widetilde{X}$ is a max-stable process with 
representation \eqref{eq1} where $Z^{(i)}$'s are independent copies of $\norm{Z}$, then  we have 
$$ M_n^*  \le \sup_{t \in [0,n]^l\cap \TT} \widetilde{X}(t) ,\quad n \in \mathbb{N}$$
almost surely. Since when $\nu_Z$ is shift-invariant \eqref{tcfN} holds and thus as mentioned before by \Cref{shikon} $\widetilde{X}$ is   stationary and therefore regularly varying as shown in \Cref{thm:SP}. Hence applying  \Cref{genUB} yields \eqref{rak2} establishing the claim. 
\QED
}
     
\def\SS{\TT}
\def\IqZ{\mathcal{S}^q(Z)}

\begin{appendices}
\section{Space $D(\R^l, \R^d)$ \& the mapping theorem}
\label{sec:spaceD}
The space of generalised \cadlag\ functions $f:\R^l\mapsto \R^d$ denoted by $\spaceD=D(\R^l, \R^d)$  is the most \mb{commonly used when defining} random processes. If $U$ is a hypercube of $\R^l$ we define similarly $\spaceD_U=D(U,\R^d)$ which is Polish (see \eg\ \cite{MR2802050}[Lem 2.4]) and will be equipped with the $J_1$-topology, the  corresponding metric and its Borel $\sigma$-field. 
The case $l=1$ is \mb{the most} extensively studied in numerous  contributions as \mb{highlighted} in \Cref{sec:prim}. There are only a few articles  dealing with properties \mb{of $\spaceD$} when $l>1$ focusing mainly on the space of \cadlag\ functions $D([0,\IF)^l, \R^d)$, see \cite{MR625374,PhDTH}.\\
  The definition of $\spaceD$  for $l\ge 1$ needs some extra notation and therefore we directly refer to \cite{PhDTH} omitting \mb{the} details. \\
Let $\TO$ be a dense set of $\R^l$. Given a hypercube $[a,b] \subset \R^l$ set $K=[a,b]  $ and  write $P_k(K,\eta),\eta>0$ for a partition  of $K= \bigcup_{i=1}^k K_i$ by disjoint hypercubes   $K_i=[a_i,b_i]  \subset {\TTT }, i\le k$ with {smallest length} of $[a_i,b_i]$'s  exceeding  $\eta$   and  let $P(K,\eta)$ be the set of all \tE{such partitions}.\\
 \cEE{We define for a given norm $\norm[*]{\cdot}$ on $\R^d$ and $\eta>0, f\in \spaceD$} 
\bkn{	\label{wp} 
		\tE{	w(f,K,\eta)} & =& \sup_{  s,  t \in K \cap \TO}  \norm[*]{f(t)-f(s)}, \\
		 	w'(f,K,\eta) &=&  \tE{\inf_{  
			P_k(K,\eta) \in 	P(K,\eta)} \max_{1 \leq i \leq k} \max_{ s,t \in K_i \cap \TO} \norm[*]{f(t)-f(s)} },\\
				w''(f,K,\delta)  &=&  \sup_{  \tE{  s,t,u   \in K \cap \TO } \atop   s \le t \le u \le s+\delta 
	} \min( \norm[*]{f(t)-f(s)},  \norm[*]{f(u)-f(t)}) .  \label{wpp}
}

Let $\tau$ be time changes 
$\R^l \mapsto  \R^l$,  \ie\ its  components denoted by $\tau_i: \R \mapsto \R  , i\le l$ 
are strictly increasing, continuous, $\tau_i(-\IF)=- \tau_i(\IF)=-\IF$  and     such that their slope norm 
$$\abs{\normE{\tau_i}}=\sup_{s\not= t, s,t\in \R} \abs{ \ln (( \tau_i(t)- \tau_i(s))/(t-s))} $$ is finite.   Write $\Lambda$ for the set of all $\tau$'s. As in \cite{PhDTH} we define the metric $\Dtau$ for all $f,g\in \spaceD$ by 
\begin{equation}\label{Ddt}
	\Dtau(f,g)= \sum_{ j=1}^\IF  2^{-j} \min(1, d_{N(j)}(f,g)), \quad f,g\in \spaceD,
\end{equation} 
where $N(j), j\le \Nset^l$ is an enumeration \mb{of} $\Nset^l$ and $d_{N(j)}(f,g)$ is as in \cite{PhDTH}[Eq. (2.26)], i.e., 
$$  d_{N}(f,g)= \inf_{\tau \in \Lambda} \mb{\left( \sum_{i=1}^l \abs{\norm {\tau_i}}+ \max_{t\in \R^l \cap \TO} \norm[*]{(k_N(\tau(t)) 
	 \cdot f(\tau(t)) - 
k_N(t)\cdot g(t)}\right)}, \quad N =(N_1 \ldot N_l)\in \Nset^l, $$
where $\cdot$ is the Hadamard product (the usual component-wise product) with  $k_N: \R^l \mapsto \R^d$ having components 
$$
k_{N_i}(t)=1, \ t\in [-N_i,N_i], \ k_{N_i}(t)=0, \ t\in [-N_i-1,N_i+1]^c,
$$ 
 and for  other $t \in \R^l$ the components  $k_{N_i}(t), i\le d$ are defined by linear interpolation.  Here the hypercube $[-N,N]$ is defined as usual for $N \in \Nset^l$.\\
Since for all $N(j)\in \Nset^l$ such that $[-N(j),N(j)] \subset [-N,N]=[-k, k]^d, k\in \Nset$ we have $d_{N(j)}(f,g) \le d_N(f,g)$    by the definition of $d_N$ it is clear that 
$$ 
d_{N(j)}(f,g) \le \sup_{t\in \TO \cap [-k,k]^d}  \norm[*]{f- g}
$$
we conclude that for all $\eta>0$ \eEE{we can find $k>0$ independent of $f$ and $g$ such that}
\akn{ \Dtau(f,g) \le  \sup_{t\in \TO \cap [-k,k]^d}  \norm[*]{f- g} + \eta  
\label{prs}
}
holds for all $k$ sufficiently large. 
Let $J_1$ be the Skorohod topology on $\spaceD$, \ie\ the smallest topology on $\spaceD$ such that  $\limit{n}f_n = f$ \mb{holds} if and only if there exists $\tau_n \in \Lambda$ such that: 
\begin{enumerate}[{${J_1}$}a)]
	\item \label{Si} $\limit{n} \sup_{s\in \R}\abs{\tau_{ni}(s)- s}=0, \ \forall 1\le i \le l $;
\item  \label{Sii} $\limit{n}\sup_{t\in [-N,N]}\norm[*]{ f_n(\tau_n(t)) - f(t)}=0, \ \forall  N \in \Nset^l.  $
\end{enumerate}  
Let  $X_{\mathbb{P}}$ denote the set of stochastic continuity points $u\in \TTT$ of the $\spaceD$-valued random element $X$ defined on $\Pspace$, \ie\ $X_{\mathbb{P}}$  consists of  all  $u\in \TTT$ such that the image measure  $\mathbb{P} \circ X^{-1}$ assigns mass 0 to the event 
$\{ f\in \spaceD: \text{ $f$ is discontinuous at $u$}  \}$.  Write similarly $T_\nu$ for the sets of continuity points of a measure $\tailmeasure$ on $\AA$. 
The next result is \mb{largely} a collection of several known results in the literature.  
\begin{theorem}
\begin{enumerate}[(i)]
		\item \label{tev:0} $f\in \spaceD$ if and only if 
		$$\lim_{\delta \to 0} w'(f,K,\delta)=0,\quad \  \sup_{t\in \TO \cap K} \mb{\norm[*]{f(t)}}< \IF, \ \ \forall 
		K\in \TTTK;
		$$
	\item \label{tev:1} $(\spaceD,\Dtau)$ is a Polish   space and $\Dtau$ generates the $J_1$ topology;
\item  \label{tev:2} 	 $\sigma(\proj_t, t\in \TO)=
\borel{\spaceD};$
\item  \label{tev:3} The  pairing  $(z, f) \mapsto z f$ which is \mb{a} group action of $\RP$ on $\spaceD$ is continuous in the product topology on $\RP\times \spaceD$;
\item \label{tev:Mreis} \eEE{$A \in \mcb_0$ if and only if \eqref{heart} holds for some  $\ve_A>0$ and some hypercube $K_A \subset \R^l$;} 
\item  \label{tev:6} For all  $f\in \spaceD$ such that $\norm{f(0)}\ge 1$ 
we have $\Dtau(  c f,0) =1$ for all $c>1$;
	\item \label{tev:4} A sequence of $\spaceD$-valued random elements $X_n, n\ge 1$ defined on an non-atomic probability space $\Pspace$ converges weakly as $n\to \IF$ with respect to the $J_1$ topology to a $\spaceD$-valued random element $X$ defined on the same probability space, if $(X_n(t_1)\ldot X_n(t_k))$ converge in distribution as $n\to \IF$ to $(X(t_1)\ldot X(t_k))$ for all $t_1 \ldot t_k$ in   $ X_{\mathbb{P}}$    and further 
	$$ \lim_{\eta\downarrow 0} \limsup_{n \to \IF} \pk{w'( X_n, K, \eta) \ge \ve  }=0,\quad  \forall \ve>0, \forall K \in \TTTK, $$  
	$$ \limit{m} \limsup_{n \to \IF} \pk*{ \sup_{t\in [-k,k]^l \cap \TO}  \norm[*]{X_n(t)} \ge  m  }=0, \quad \forall k\in \Nset;
	 $$  
\item \label{tev:7} The sets  $X_{\mathbb{P}}$ and $T_\nu$ for $\nu$ a $\sigma$-finite measure on $\borel{D}$ is dense in $\TTT$ 
\cEV{and 
$\nu(disc(\proj_{t}))=0, \forall t\in T_\nu$};
\item \label{tev:8} If $\nu$ is a $\sigma$-finite measure on $\borel{D}$, then for any hypercube $A \subset \R^l$ with corners in $T_\nu$
we have $\nu(Disc(\proj_A))=0$, where  $\proj_A: \spaceD \mapsto D(A, \R^d)=\spaceD_{A}$ with $\proj_A(f)= f_A, f\in \spaceD$  the restriction of $f\in \spaceD$ on $A$. Moreover we can find increasing hypercubes $A_k, k\in \Nset$ with $\nu(Disc(\proj_{A_k}))=0$ such that $[-k,k]^l \subset A_k$ for all $k\in \Nset$;
\item \label{tev:9} \cEV{If  $A_k$ is as in \Cref{tev:8},  the projection map $\proj_{A_n, A_k}: \spaceD_{A_n} \mapsto \spaceD_{A_k} $ with $A_k \subset A_n$ or 
	$A_n=\R^l$   is measurable;}
\item \label{tev:10} \cEV{Let $\nu_n,n\in \Nset, \nu$ be finite   measures on $\borel{\spaceD}$ and let  $A_k,k\in \Nset $ be as specified in \Cref{tev:8} above. If  $\nu_n \circ \proj_{A_k}^{-1}\convweak \nu \circ \proj_{A_k}^{-1}, \forall k\in \Nset$, then $\nu_n \convweak\nu$ as $n\to \IF$.} 
\end{enumerate}	
\label{sam}
\end{theorem}

\prooftheo{sam}
\Cref{tev:0}: This is shown in \cite{PhDTH}[Thm 2.1] for the case $D([0,\IF)^l, \R^d)$ and can be proved with similar arguments for $D(\R^l,\R^d)$.\\
\Cref{tev:1}:  For the case $D([0,\IF)^l, \R^d)$ this is shown in  \cite{PhDTH}[Thm 2.2]. The case $D(\R^l, \R^d)$ follows with the same arguments, 
see also \cite{MR625374}.\\
\Cref{tev:2}: The last two equalities are  shown in \cite{MR625374}[Thm 3.2] for the case $D([0,\IF)^l, \R^d)$ (see also 
\cite{PhDTH}[Thm 2.3]) and can be shown with similar arguments for our setup. 
\COM{It is clear that $\sigma(\proj_t, t\in \TO) \subset \mb{\sigma(\proj_A, A \subset \TO)}$. The converse inclusion follows by showing that 
$\proj_T$ is measurable with respect to \mb{$\sigma(\proj_t, t\in A)$ for all $A\subset \TO$}. This follows by the fact that the product $\sigma$-field 
of a product space (with countable product) of separable metric spaces is given by the corresponding product of $\sigma$-fields, see \cite{Kallenberg}[Lem  1.2]. Note that our spaces are Polish and therefore separable.}\\
\Cref{tev:3}:  We need to show that for all positive sequences  $a_n \to a>0$ as $n\to \IF$ and any $f_n, f\in \spaceD$ such that $\limit{n} \Dtau(f_n,f)$ we have 
$\limit{n} \Dtau(\mb{a_n f_n}, af)=0$. By  the characterisation of the Skorohod topology  there exists $\tau_n$ such that 
\Cref{Si} and \Cref{Sii} hold. Since 
$$\norm[*]{ a_nf_n( \tau_n(t)) - af(t)}\mb{\le} a_n \norm[*]{f_n( \tau_n(t)) - f(t)} +
\abs{ a_n- a} \mb{\norm[*]{f(t)}}$$ and $\sup_{t\in K} \norm[*]{f(t) }$ is finite for any compact $K\subset \R^l$, then the claim follows.\\
\Cref{tev:Mreis} \eEE{By the equivalence of the norms on $\R^d$ and the definition of $\norm[t]{\cdot}$ in the formulation of \Cref{A:2}, we can assume without loss of generality that $\norm[t]{f}=\norm[*]{f(t)}, f\in \spaceD, t\in \TTT$.\\
	  We have that $A \in \mcb_0$ if and only if  there exists $\ve_A>0$ such that for all $f\in A$ we have  $\Dtau(f,\zD)>\ve_A$. Hence for such $A$,  by \eqref{prs} 
there exists $\ve'\in (0,\ve_A)$ and some hypercube $K_A$ such that 
$$ f_{K_A}^*=\sup_{t\in K_A \cap \TO} \norm[*]{f(t)} > \ve'$$
for all $f\in A$. \\  
 Conversely, if for all $f\in A $ we have $f_{K_A}^*> \ve>0$ and thus $f_{[-k,k]^l}^*> \ve$ for all $k$ sufficiently large, since  
$$d_N(f,0) \ge \sup_{t\in [-k,k]^l \cap \TO} \norm[*]{f(t)}=f^*_{[-k,k]^l}, \ \forall N \in \Nset^l\setminus [-k,k]^l, $$
 then by definition of $d_\spaceD$ we have that $d_\spaceD(f,0)\ge f^*_{[-k,k]^l}> \ve'$ for some $\ve'>0$ and all $f\in A$, this means that $A \in \mcb_0$ by the definition of $\mcb_0$ establishing the claim.} \\
\Cref{tev:6}:  By the definition,  for all $c>0, f\in \spaceD$ and  $N(j) \in \Nset^d$  (recall $0$ denotes the zero function in $\spaceD$)
$$d_{N(j)}(cf,0)= (cf)^*_{[-N(j), N(j)]}= cf^*_{[-N(j), N(j)]}\ge c \norm{f(0)}.$$
Hence if $\norm{f(0)}\ge 1$, then 
$$ \Dtau(cf,0)=  \sum_{ j=1}^\IF  2^{-j} \min(1, d_{N(j)}(cf,0))= 1=\Dtau(f,0), \quad c>1.$$
\Cref{tev:4}:  The tightness criteria is given in \cite{PhDTH}[Thm 2.4]. The claim follows now from \cite{MR0402847}[Thm 5.5].\\
\Cref{tev:7}: The fact that $X_\mathbb{P}$ is dense in $\R^l$ is shown in  \cite{MR625374}[p.~182] for $\spaceD=D([0,\IF)^l, \R^d )$ and hence the claim follows for $\spaceD=D(\R^l, \R^d)$.  We can use that result and $\sigma$-finiteness of $\tailmeasure$ to prove that $T_\nu$ is also dense in $\R^l$. Next, for all $t\in T_\nu$ we have that $\proj_{t}$ is continuous for almost all $f\in \spaceD$ with respect to $\nu$ if and only if 
$\nu(\{f\in \spaceD: f_t\not=f_{t-}\})=0$, hence the claim follows.\\
\Cref{tev:8}: \mb{The proof for $A$ is along the lines of} \cite{MR625374}[Lem 4.2] for a probability measure on $\spaceD$ and the argument can be extended to a $\sigma$-finite measure $\nu$. Since
$T_\nu$ is dense in $\R^l$, then $A_k$ with the stated property exists. \\
 \Cref{tev:9}: The case $l=1$ is shown for instance in \cite{sagitov2020weak}[Lem 9.20], where $A_n=\R$. The general case $l$ is  a positive integer  and $A_n$ is a hypercube that includes $A_k$ follows with  similar arguments as therein using further 
\Cref{tev:2} above.\\  
\cEV{\Cref{tev:10}: For probability measures $\nu_n, \nu$ this is the shown in \cite{MR625374}[Thm 4.1] and the remark about proper sequences after the proof of \cite{MR625374}[Thm 4.1]. However,  the proof of the aforementioned theorem as well as the corresponding result \cite{Lindwall}[Thm 3, 3'] have a non-fatal gap, namely the projection map denoted by $r_\alpha$  therein has not been shown to be measurable and therefore the mapping theorem cannot be applied as claimed. The measurability of $r_\alpha$ for $l=1$ is shown in \cite{MR561155}[Lem 2.3] and  the case $l>1$  is claimed in \Cref{tev:9} above. 
 The claim for finite non-null measures $\nu_n,\nu$ follows then,  since the weak convergence implies that $\limit{n}\nu_n(\spaceD)=\nu(\spaceD)\in (0,\IF)$ and hence $\nu_n/\nu(\spaceD), \nu/\nu(\spaceD)$ are probability measures and we have the corresponding weak convergence. }
\QED 
 
Concluding, we present  the  mapping theorem  under the assumption \Cref{A:5} for both $\spaceD$
and $\spaceD'$  equipped with  properly localised boundednesses $\mcb$ and $\mcb'$, respectively. 
\begin{theorem} \label{mapping} (\cite{kulik:soulier:2020}[Thm B.1.21])
		Let $\tailmeasure_z,z>0$ be $\mcb$-boundedly finite  measures on $\borel{\spaceD}$ and let  $\nu$ be a  measure on $\borel{\spaceD'}$. 
		If   $H: \spaceD \mapsto \spaceD'$ is   $\borelD/\borel{\spaceD'}$ measurable  
		and  $\nu_z \convvague\nu$,  then  
		$\nu_z  \circ H^{-1} \stackrel{v,\mcb'}{\longrightarrow} \tailmeasure \circ H^{-1}$, provided that 
		\begin{equation}\label{nGSTT}
			H^{-1}(B) \in \mcb, \ \  \forall B \in \mcb' \cap \borel{D'} \ \text{  and  } \tailmeasure(Disc(H))=0.
		\end{equation}
	\end{theorem}

\end{appendices}

{\bf Acknowledgements}:  \cEG{We are grateful to both reviewers and the Editor for several comments and  suggestions  that lead to significant improvement of the manuscript. Many thanks go to  Richard Davis, Philippe Soulier, Thomas Mikosch and Ilya Molchanov  for   literature suggestions and  discussions. 
MB kindly acknowledges support from the SNSF Grant 200021-191984 and EH,  GS kindly acknowledge support from the SNSF Grant 1200021-196888. }
  
\bibliographystyle{ieeetr}
\bibliography{R.regVarTS}
\end{document}